\documentclass[]{article}
\usepackage{amsmath}
\usepackage{setspace}
\usepackage{graphicx}
\usepackage{nameref}
\usepackage{color}

\usepackage[margin=2cm]{geometry}
\usepackage[toc,page]{appendix}

\begin{document}

\vspace*{.6cm}

\centerline{\LARGE\bf Tipping points near a delayed saddle node bifurcation with}
\centerline{\LARGE \bf periodic forcing}

\vspace*{.6cm}

\centerline{\large Jielin Zhu\footnote{Department of Mathematics, UBC}, Rachel Kuske\footnote{Department of Mathematics, UBC}, Thomas Erneux\footnote{Optique Nonlineaire Theorique, Universite Libre de Bruxelles, BE}}

\vspace*{.6cm}

\noindent
{\bf Keywords:} tipping points, multiple scales, delayed bifurcation, periodic forcing, asymptotic expansions

\vspace*{.3cm}
\noindent
{\bf AMS subject classifications:} 34E13, 37B55, 35B32, 37L10, 34E10 



\begin{abstract}
 We consider the effect on tipping from an additive periodic forcing in a canonical model with a saddle node bifurcation and
a slowly varying bifurcation parameter.
Here tipping refers to the dramatic change in dynamical behavior characterized by a rapid transition away from a previously attracting
state.  In the absence of the periodic forcing, it is well-known that a slowly varying bifurcation parameter produces a delay
in this transition, beyond the bifurcation point for the static case. Using a multiple scales analysis, we consider the
effect of amplitude and frequency of the periodic forcing relative to the drifting rate of the slowly varying bifurcation parameter.
 We show that a high frequency oscillation drives an earlier tipping when the bifurcation parameter varies more slowly,
with the advance of the tipping point proportional to the square of the ratio of amplitude to frequency.
In the low frequency case the position of the tipping point is affected by the frequency, amplitude and phase of
the oscillation. The results are based on an analysis of the
 local concavity of the trajectory, used for low frequencies both of the same order as the
drifting rate of the bifurcation parameter and for low frequencies larger than the
drifting rate. The tipping point location is advanced with increased amplitude of the periodic forcing, with
critical amplitudes where there are jumps in the location, yielding significant advances in the
tipping point. We demonstrate the analysis for two applications with saddle node-type bifurcations.
\end{abstract}

\section{Introduction}

Generally speaking, the term tipping point refers to an abrupt transition in dynamical behavior observed as  the system 
moves to a qualitatively different state, due to small changes in one or more
 factors \cite{Lenton2008}.  Tipping points have been observed
in many different fields, including the start and end of ice-ages \cite{sutera1981}, environmental regime shift  \cite{guttal2008}, synchronized behaviour of neural activity \cite{Meisel2012}, catastrophic collapse in ecology,
 and power systems \cite{Dai2012}.
As many of these transitions are often irreversible or difficult to reverse, it becomes important to develop analytical techniques for predicting the location of tipping points and understanding the key contributing factors to tipping.

A number of
analyses indicate that tipping points can depend on the bifurcation structure of the system, external noise or forcing,
 or the drifting rate of parameters that vary in time \cite{ashwin2012}.
 Of course in applications,
 a combination of these factors may contribute to the system simultaneously. 
A significant number of studies have focused on systems that have a
 slowly varying parameter approaching a saddle-node bifurcation point, where the transition occurs 
for parameter values beyond that of the static bifurcation point, known as a delayed bifurcation. Earlier  recognition of this phenomenon is discussed in the
context of optics \cite{mandel1987} and neuronal dynamics \cite{baer1989}, as well
as general bifurcation theory \cite{haberman1979}.
The structure appears in many applications, including  models of global energy balance \cite{fraedrich1978}, extent of Arctic sea ice \cite{eisenman2009,abbot2011}, and population resilience \cite{Dai2012}, to name a few.    Focal points
for prediction of tipping include
 statistical measures based on historical data
that track the value of the slowly varying parameter
 where an abrupt
transition occurs \cite{sieber2012}. These types of measures have been compared with canonical models, illustrating
that in the presence of noise tipping can occur well before the slowly varying parameter reaches
 the static bifurcation point of the model
\cite{thompson2011}.  
Whether the location of the tipping point is advanced, that is ``early'', or delayed
 relative to the static bifurcation point 
depends on the relationship between the drifting rate 
and the strength of the noise, as studied for different underlying bifurcation structures
 in \cite{berglund2010, kuske1999}.
Predictive techniques using statistical measures and time series analyses 
have been developed to identify signatures of proximity to a tipping point or potential for early
tipping, with examples in 
\cite{ guttal2008,Dai2012,sieber2012,dakos2008} and additional references in \cite{Scheffer2012}.
 However, those measures are not necessarily applicable if the tipping is driven by changes in
the drifting rate or other factors, indicating that a better understanding of the role of
different tipping mechanisms is necessary.

An additional factor appearing in many applications is periodic variation.
Often reduced models, such as those for the extent of Arctic sea ice \cite{eisenman2009, abbot2011,merryfield2008},  
 are developed for analysis by averaging over oscillations with higher frequency, but a closer look
at these types of models indicates different locations of the tipping points when the oscillations are included. 
Periodic variation can occur on different time scales, such as the seasonal variation
considered in \cite{eisenman2009} or regular fluctuations on the time scale of decades \cite{tung2013}.
To explore the effect of periodic forcing with either low or high frequency, we consider the 
canonical model for a saddle node bifurcation with a slowly varying bifurcation parameter
 and  additive periodic forcing. In this setting we
explore how the location of the tipping point depends on  the combined effect of 
frequency, amplitude, and the drifting rate at which bifurcation parameter approaches the static bifurcation point. The different expansions that are necessary to cover the different scenarios provide explicit expressions for the tipping points and capture the key contributions that determine them.

For high-frequency forcing, we use a multiple scales approach combined with layer expansions \cite{Bender1978,Kevorkian1996}.
The outer expansion uses a multiple scales expansion based on the times scales of the oscillations and the drifing rate.
The outer expansion does not provide
the approximation for the tipping point, but indicates the location of the inner expansion 
composed of an oscillatory part plus a correction.
 Combining several multiple scales expansions in terms of fractional powers of the drifting rate to
get the inner expansion, we find that the tipping point is determined from the
 averaged equation for the correction to the oscillatory contribution. 
The form of the averaged equation is similar to that of the case without oscillations, facilitating
 explicit expressions for 
 two main contributions to the location 
of the tipping point relative to the static bifurcation point:
 a delay due to the slowly drifting bifurcation parameter, and an advance
proportional to the square of the ratio of amplitude to frequency for the oscillation.
The shift due to high frequency periodic forcing is also compared with the case of white noise forcing.
An additional rescaling is required in the case of both large amplitude and high frequency.  That rescaling indicates that in 
certain limits, the appropriate expansion for the rescaled system is the low frequency approximation rather than the high frequency case as might be expected.

The multiple scales approach used for the high frequency oscillations is not applicable for low frequency oscillations. Rather, the oscillatory outer or regular solution can be written in terms 
of the slow time only, and where it does not exist, 
we find two different types of tipping points.  For the first type, the outer solution vanishes but its derivative does not. Then
a local approximation for the equation of the trajectory has a singularity corresponding to the tipping point.  The influence
of the oscillation appears through the derivative of the square of the outer solution. In contrast, for the second type of tipping point the derivative of the
square of the outer solution also vanishes. Then a different local approximation is used to consider the concavity of the
trajectory. The concavity indicates whether or not tipping occurs near this location.
 In contrast to the high-frequency forcing, where the location of the 
tipping point varies smoothly with the other parameters, for parameter values where the concavity in the low frequency case changes
 there are discontinuities in the location of the tipping point. 
From the analysis of these two types of tipping points, 
we find that their location is affected
not only by the frequency,  the amplitude, and the drifting rate, but also by  the behavior of
the trajectory as the oscillations cause it to approach zero. This effect can be related to the phase of
the oscillations, given by an initial condition, and can lead to either delayed or advanced tipping. 

We demonstrate the method on two applications, namely a Morris-Lecar model and 
an energy balance model describing the persistence of Arctic sea ice. 
In both cases, we use 
the approaches 
developed for the canonical model directly on normalized versions of these models.  The analysis illustrates the
importance of this normalization, as it is within the rescaled system that we can identify the relationships and
 approximations that lead to quantitative expressions for the tipping point.
We also see that the results can be extended beyond the 
range of asymptotic relationships for drifting rate, frequency, and amplitude 
obtained in the context of the canonical model; for example, we consider
the scenario where the ratio of amplitude to frequency is outside of these limits.
The application of the approaches in these extended ranges is based on recognizing 
the leading order contributions in the outer and local expansions
that provide information about the relevant features for tipping.
 These observations guide the inclusion of the most important terms in 
reduced equations with singularities corresponding to tipping points.
Even without deriving these reduced equations, the normalized systems combined with the results from the canonical model already 
indicate parameter ranges where advanced tipping is predicted, as we discuss in the context of the specific applications.


In Section \ref{sec:model}, we introduce the canonical model with an additive periodic forcing. We briefly review the results when there is no oscillatory term. In Section \ref{sec:hf}, we show that a high frequency oscillation 
changes the position of the bifurcation point for the case of a constant bifurcation parameter, and likewise
  triggers an early tipping when the drifting rate of the bifurcation parameter is slow enough. 
A brief  comparison of the additive noise and the periodic forcing cases is also shown.
In Section \ref{sec:lf}, we show that the position of the tipping point is affected by the frequency, amplitude and phase of
the low frequency oscillation. For the case where the frequency of the oscillatory forcing is the same
order of magnitude as the drifting rate of the bifurcation parameter we develop an analysis based on the
 local concavity of the trajectory to find the threshold of the amplitude at which there is an
abrupt change of the position of the tipping point. Then we adapt this analysis
to consider cases forced by a periodic oscillation with low frequency that is larger than the drifting rate 
of the bifurcation parameter.  For consideration of the system with large amplitude, we study
a rescaled model to which the analysis from the low or high frequency cases can be applied 
for certain parameter combinations. 
In Sections \ref{sec:model_ML} and \ref{sec:sea_ice}  we consider
two models with saddle node-type bifurcations, transformed to be
similar to the canonical model. Then we can directly apply the approaches
from the previous sections to determine the tipping points in these examples.
 
Throughout the paper we compare our analytical results to simulations of the model using a second order
 Runge-Kutta method, except for the Morris-Lecar model in Section 5, where we use a fourth-order Runge-Kutta
method.

\section{The canonical model}
\label{sec:model}

We consider a canonical model with additive periodic forcing,
\begin{eqnarray}
\frac{\mbox{d}x}{\mbox{d}t}&=&a-x^2+A \mbox{sin}(\Omega t), \quad
\frac{\mbox{d}a}{\mbox{d}t}=-\mu, \nonumber\\
 & & x(0)=x_0,\  a(0)=a_0>0,
\label{sys:gen_periodic_forcing}
\end{eqnarray}
where $a(t)$ denotes the slowly varying bifurcation parameter with the drifting rate $0< \mu \ll 1$, $A$ is the amplitude of the oscillation ($A>0$), and $\Omega$ represents the frequency. We assume $a_0$, $A$ are positive constants and $x_0=\sqrt{a_0}$.

\begin{figure}[ht]
\centering 
\scalebox{0.5}
{\includegraphics{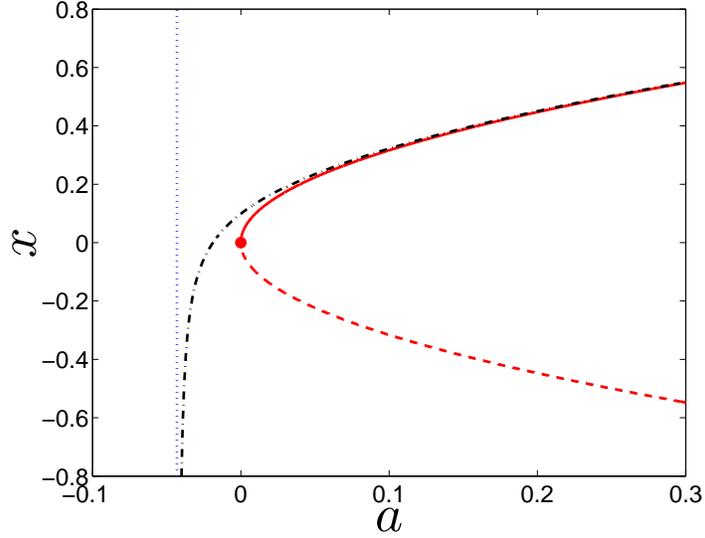}}
\caption{The asymptotic approximation of the slowly varying equilibrium solution (\ref{sol:slow_vary_eq_1}), (\ref{sol:slow_vary_eq_2}) of the system (\ref{sys:gen_periodic_forcing}) without periodic forcing ($A=0$) (black dash-dotted) is compared with the saddle node bifurcation diagram (red solid and dashed lines). The solid (dashed) line 
is the stable (unstable) steady state. The solid circle indicates the static bifurcation point. The vertical blue dotted line indicates the value of $a_d$ (\ref{sol:location_tipping}), the singularity of (\ref{sol:slow_vary_eq_2}). The drifting rate is $\mu= 0.0025$.}
\label{fig:delay_effect}
\end{figure}

The system (\ref{sys:gen_periodic_forcing}) is generalized from the
canonical model for a saddle node bifurcation when $a(t)$ is a constant and 
$A=0$ \cite{strogatz2000}. For any constant $a>0$, there exist a stable equilibrium $x_e^+ = \sqrt{a}$ and an unstable equilibrium $x_e^- = -\sqrt{a}$. 
For $a<0$, no equilibrium exists and the trajectory of (\ref{sys:gen_periodic_forcing}) decreases exponentially for any initial condition $x_0$. 
The bifurcation diagram is shown in Figure \ref{fig:delay_effect}.

For the system  (\ref{sys:gen_periodic_forcing}) without the periodic forcing term ($A=0$), there is no equilibrium because $a(t)$ is slowly varying. With initial conditions $x_0$ and $a_0>0$, the trajectory of the system is exponentially attracted to one specific solution which we call the slowly varying equilibrium solution, analyzed in
 \cite{haberman1979, thompson2011} for $A=0$ which we briefly summarize here.
For $a(t)>0$, the slowly varying equilibrium solution is close to $x=\sqrt{a}$ for $a$ constant which we call the stable branch throughout this paper. As $a(t)$ crosses zero, there is a transition from the slowly varying equilibrium solution to the rapidly decreasing behaviour. This transition occurs for the value $a(t)$ below the static saddle node bifurcation point at $a=x=0$, and thus is called a delayed bifurcation (see \cite{haberman1979} and references therein). When such an abrupt qualitative change occurs, it is often referred to as a tipping point regardless of the value of $a(t)$. The method of matched asymptotic expansions is applied to find the approximation of the slowly varying equilibrium solution for $A=0$
 \cite{haberman1979,berglund2010}, using a slow
time scale $\mu t$ and 
solving for $x(t)$ as a function of $a(\mu t)$,
\begin{equation}
x(a(\mu t)) = \sqrt{a}+\frac{\mu}{4a}-\frac{5}{32}\frac{\mu^2}{a^{5/2}}+O(\mu^3), \label{sol:slow_vary_eq_1}
\end{equation}
for $a(\mu t) \sim O(1)$. Notice that if $a(\mu t)$ is  $O(\mu^{2/3})$, the terms shown in (\ref{sol:slow_vary_eq_1}) are all $O(\mu^{1/3})$ and the
expansion is no longer valid. 
Then  a local approximation is needed, which yields
\begin{equation}
x(a(\mu t)) \sim -\mu^{1/3}\frac{\mbox{Ai}'(a/\mu^{2/3})}{\mbox{Ai}(a/\mu^{2/3})}, \label{sol:slow_vary_eq_2}
\end{equation}
 where $\mbox{Ai}$ denotes the Airy function. The asymptotic approximation of the slowly varying equilibrium solution of (\ref{sys:gen_periodic_forcing}) with $A=0$ is shown in Figure \ref{fig:delay_effect}, and compared with the static saddle node bifurcation diagram. The tipping happens near the singularity in the expression (\ref{sol:slow_vary_eq_2}), which is the first zero of the Airy function,
\begin{equation}
a_d = \mu^{2/3} \cdot (-2.33810...) \, . \label{sol:location_tipping}
\end{equation}
The smaller that $\mu$ is, the slower $a(t)$ varies, and the closer that the tipping point is to the static bifurcation value $a=0$. The $\mu^{2/3}$ law characterizes the delay of the jump transition and has been verified experimentally in various physical systems \cite{erneux1989,Jung1990,Hohl1995}. This scaling law reappears throughout the results in this paper, due to the parabolic nature of the
bifurcation structure near a saddle node bifurcation.

For $A>0$ in (\ref{sys:gen_periodic_forcing}), we explore
the effect of the amplitude and frequency of the oscillation on the location of the tipping point,  
considering both high and low frequency in combination with the slow drift of the parameter $a(t)$.


\section{High Frequency Oscillation $\Omega \gg 1$}
\label{sec:hf}

In order to see how the tipping location is influenced  by high frequency periodic forcing, 
we first consider (\ref{sys:gen_periodic_forcing}) with constant $a$ ($\mu=0$) and determine
the critical value of the bifurcation parameter where the bounded attracting solution loses stability.
We identify this bifurcation point denoted $a_p$, analogous to the saddle node bifurcation point in Section \ref{sec:model}, using the methods of multiple scales and matched asymptotic expansions. This result shows that the
periodic forcing shifts the location of the bifurcation point to a value of $a>0$.  We compare this effect with the shift in the location of the tipping point for combined
 slowly varying $a(t)$ and high frequency forcing. 
 We obtain an expression for this shift in terms of the drifting rate $\mu$ of $a(t)$, the frequency $\Omega$, and the amplitude $A$.  
First we consider values of $A$ smaller than $\Omega$ or of the same order of magnitude of $\Omega$. For larger values of $A$, a different approach is needed, as shown in Section \ref{sec:hf_large_A}.


\subsection{Constant Bifurcation Parameter}
\label{sec:hf_no_drifting}

For $\mu = 0$ ($a(t)$ a constant) and $A \neq 0$, the solution of (\ref{sys:gen_periodic_forcing}) oscillates periodically around $\sqrt{a}$ for $a$ an $O(1)$ constant but it exponentially decreases for parameter values below the bifurcation value of $a=a_p$. We locate $a_p$, using the method of multiple scales \cite{Bender1978}, 
based on both a fast time scale $T = \Omega t$ and the original $O(1)$ time scale $t$, the slow time scale in this case. 
For $x = x(T,t)$, (\ref{sys:gen_periodic_forcing}) becomes
\begin{equation}
x_T +\Omega^{-1} x_t = \Omega^{-1} a -\Omega^{-1} x^2 + \Omega^{-1} A \mbox{sin}(T), \mbox{ where } 
x_T=\frac{\partial x}{\partial T}, \ 
x_t=\frac{\partial x}{\partial t}
 \label{exp:hf_multi_scale}
\end{equation}
 and the expansion
\begin{equation}
x \sim x_0(T,t)+\Omega^{-1} x_1(T,t)+\Omega^{-2} x_2(T,t)+ \cdots, \label{exp:hf_gen_expansion}
\end{equation}
is substituted into (\ref{exp:hf_multi_scale}), yielding
\begin{align}
O(1): x_{0T} &=0, \label{exp:hf_constant_order0} \\
O(\Omega^{-j}): x_{jT} &= R_j(T,t), \,\,\, j=1,2,3 \cdots  \, . \label{exp:hf_constant_general}
\end{align}
The leading order equation (\ref{exp:hf_constant_order0}) indicates that $x_0 = x_0(t)$ so that $x_0$ depends only on the slow time variable $t$.

To find $x_j$ for $j=0,1,2,\cdots$, we use a solvability condition to eliminate secular terms in the solution of (\ref{exp:hf_constant_general}) that monotonously increase in the fast time variable $T$.
The solvability condition is that the right hand side of (\ref{exp:hf_constant_general}) is orthogonal to the homogeneous solution of $x_{jT}=0$ \cite{Bender1978},
\begin{equation}
\lim_{T \to \infty} \frac{1}{T} \int_0^{T} R_j(u,t) du = 0, \,\,\, j=1,2,3 \cdots. \label{eq:solv_condition}
\end{equation}
Under the assumptions of the method of multiple scales, $t$ and $T$ are treated as independent variables.
For $j=1$,
\begin{equation}
R_1(T,t) = a-x_0^2-x_{0t}+A\mbox{sin}(T). \label{exp:hf_constant_order1}
\end{equation}
By applying the solvability condition (\ref{eq:solv_condition}) on (\ref{exp:hf_constant_order1}), we have
\begin{equation}
\begin{array}{lll}
x_{0t} = a-x_0^2, &\Rightarrow &x_{0} =\sqrt{a}, \\
x_{1T} = A \mbox{sin}(T), &\Rightarrow &x_{1} = -A\mbox{cos}(T)+v_1(t),
\end{array}
\label{sol:hf_constant_order1}
\end{equation}
\noindent where $v_1(t)$ is a function that needs to be determined in the higher order analysis.

The higher order corrections $x_j$, $j=2,3$ are obtained in Appendix \ref{append:A}, yielding the asymptotic approximation for the attracting solution of (\ref{sys:gen_periodic_forcing}) for $\mu=0$, $a$ a constant, as
\begin{equation}
x \sim \sqrt{a}+\Omega^{-1} \cdot \left[ -A \mbox{cos}(T) \right] +\Omega^{-2} \cdot \left[ 2\sqrt{a}A\mbox{sin}(T)-\frac{A^2}{4\sqrt{a}} \right]+ \cdots. \label{sol:x_outer}
\end{equation}
Similarly, one can find the approximation of the solution for the case when $x_0 = -\sqrt{a}$ and show it to be unstable.

For $A=0$, (\ref{sol:x_outer}) reduces to $x=\sqrt{a}$ with a saddle node bifurcation point at $a=0$. For $A>0$, (\ref{sol:x_outer}) is valid when both $A$ and $a$ are $O(1)$. For $a=O(\Omega^{-2})$, we substitute $a = \Omega^{-2} b$ and introduce the expansion (\ref{exp:hf_gen_expansion}) into (\ref{exp:hf_multi_scale}) to find $a_p$. In Appendix \ref{append:A}, we find the leading order approximation of the attracting solution for $a=O(\Omega^{-2})$
\begin{equation}
x \sim \sqrt{a- \frac{ A^2}{2\Omega^2}}-\frac{A}{\Omega} \mbox{cos}(T) + \cdots. \label{sol:hf_equili_local}
\end{equation}
For $a$ a constant, (\ref{sol:hf_equili_local}) shows that there are no attracting solutions if $a < a_p = A^2/(2\Omega^2)$.  



\subsection{High Frequency Oscillation in a Slowly Varying System}
\label{sec:hf_slow_drifting}


For $a(t)$ as in (\ref{sys:gen_periodic_forcing}) with  $\mu\ll 1$, we introduce a new parameter $\lambda>0$ by writing $\Omega = \mu^{-\lambda}$ in (\ref{sys:gen_periodic_forcing}). The parameter $\lambda$ captures the relationship between frequency and drifting rate, which we relate to the tipping point below.  
  Here we consider the case $A= o(\Omega)$ for $\Omega\gg 1$, and we consider larger values of $A$ in Section \ref{sec:hf_large_A}.  To obtain an expression for the tipping location, we use a combination of the method of multiple scales together with
 outer and local  expansions.  The outer expansion does not provide the location of the tipping point, but it indicates
a new scaling for a local analysis. In the local analysis we use a different combination of multiple scales,
and obtain an expression for the tipping point in terms of the key parameters.

For the outer approximation, we find it is sufficient to use two scales, the fast time scale $T = \mu^{-\lambda}t$ 
from the oscillatory forcing and the slow time scale $\tau = \mu t$ of the parameter $a(\tau)$.  In terms of these time scales, (\ref{sys:gen_periodic_forcing}) becomes
\begin{align}
\mu^{\lambda+1} x_{\tau} + x_T &= \mu^{\lambda} \left[ a -x^2 +A \mbox{sin}(T) \right], \label{exp:hf_drifting_multi_scale} \\
a_{\tau} &=-1, \nonumber 
\end{align}
By substituting 
the general expansion
\begin{equation}
x \sim x_0(T,\tau)+\mu^{\lambda} x_1(T,\tau)+\mu^{2\lambda} x_2(T,\tau)+ \cdots, \label{exp:hf_outer_expansion}
\end{equation}
into (\ref{exp:hf_drifting_multi_scale}) and solving for each order by applying the solvability condition (\ref{eq:solv_condition}) (details shown in Appendix \ref{append:C}), we get the asymptotic approximation of $x$ for $a(\tau) = O(1)$ 
\begin{equation}
x \sim \sqrt{a}+ \frac{\mu}{4 a} + \mu^{\lambda}\left[ -A \mbox{cos}(T) \right] +\cdots \, . \label{sol:hf_drifting_outer}
\end{equation}
This result is similar to (\ref{sol:slow_vary_eq_1}), with an oscillatory correction term $\mu^{\lambda}\left[ -A \mbox{cos}(T) \right]$ from the high-frequency oscillatory term in (\ref{sys:gen_periodic_forcing}). The asymptotic approximation (\ref{sol:hf_drifting_outer}) describes the attracting solution of (\ref{sys:gen_periodic_forcing}) away from $a=0$, but does not describe tipping. As $a(t)$ approaches zero, specifically for  $a(t) \sim O(\mu^{2/3})$, the first
two terms are both $O(\mu^{1/3})$ so that the expansion (\ref{sol:hf_drifting_outer}) is not valid. In order to determine the influence of the oscillations
on the location of the tipping point compared to  (\ref{sol:location_tipping}), a local analysis for $a(t) \sim O(\mu^{2/3})$ is needed.

For the local analysis, we introduce $a(t) = \mu^{2/3} \alpha(s)$, the fast time variable $T = \mu^{-\lambda} t$ and a new slow time variable $s = \mu^{1/3} t$ and substitute into (\ref{sys:gen_periodic_forcing}). The system (\ref{sys:gen_periodic_forcing}) becomes
\begin{align}
\mu^{1/3+\lambda} x_{s} + x_T &= \mu^{\lambda+2/3}\alpha -\mu^{\lambda} x^2 + \mu^{\lambda} A \mbox{sin}(T), \label{exp:hf_drifting_local_multiple_scale} \\
\alpha_s &=-1, \nonumber
\end{align}
and we  focus on the local equation (\ref{exp:hf_drifting_local_multiple_scale}) to determine 
the location of the tipping point.
Separate analyses for different values of $\lambda$ suggest that it is useful to substitute $x = -\mu^{\lambda} A \mbox{cos}(T) +\mu^{1/3} y(T,s)$ in (\ref{exp:hf_drifting_local_multiple_scale}), yielding the asymptotic result for $\frac{1}{3}\leq \lambda <1$. Then we get

\begin{equation}
\mu^{\lambda+1/3} y_s +y_T = \mu^{\lambda}  [ \mu^{1/3} \alpha - \mu^{2\lambda-1/3}  \frac{A^2}{2} -\mu^{2\lambda-1/3}  \frac{A^2}{2}\mbox{cos}(2T)+ 2\mu^{\lambda}  A y \mbox{cos}(T)-\mu^{1/3} y^2 ], \label{exp:hf_drifting_local_2}
\end{equation}

\noindent and substitute the expansion

\begin{equation}
y \sim y_0(T,s)+ \mu^{\lambda+1/3} y_1(T,s) + \cdots \label{exp:hf_drifting_local_exp1}
\end{equation}

\noindent into (\ref{exp:hf_drifting_local_2}), yielding $y_0 = y_0(s)$,  and

\begin{align}
O(\mu^{\lambda+1/3}): y_{1T} &= -y_{0s} + \alpha -\mu^{2\lambda-2/3} \cdot \frac{A^2}{2} -\mu^{2\lambda-2/3} \cdot \frac{A^2}{2} \mbox{cos}(2T) + \mu^{\lambda-1/3} \cdot 2Ay \mbox{cos}(T)-y_0^2. \label{exp:hf_drifting_in_order1}
\end{align}

\noindent We apply the solvability condition (\ref{eq:solv_condition}) to (\ref{exp:hf_drifting_in_order1}), yielding

\begin{equation}
y_{0s} = \alpha - y_0^2 -\mu^{2\lambda-2/3} \frac{A^2}{2}, \,\,\, \Rightarrow \,\,\, y_0 = -\frac{\mbox{Ai}'(\alpha - \mu^{2\lambda-2/3} \cdot \frac{A^2}{2})}{\mbox{Ai}(\alpha - \mu^{2\lambda-2/3} \cdot \frac{A^2}{2})}. \label{sol:hf_drifting_in_order1}
\end{equation}

\noindent The local approximation of the attracting solution for $a=O(\mu^{2/3})$ is given by

\begin{equation}
x \sim \mu^{\lambda} \cdot [-A \mbox{cos}(\Omega t)]+\mu^{1/3}\cdot \left[-\frac{\mbox{Ai}'[(a-\mu^{2\lambda}\frac{A^2}{2})/\mu^{2/3}]}{\mbox{Ai}[(a-\mu^{2\lambda}\frac{A^2}{2})/\mu^{2/3}]} \right] + \cdots. \label{sol:hf_drifting_local}
\end{equation}

For $\frac{1}{6}<\lambda < \frac{1}{3}$, we substitute the expansion (\ref{exp:hf_outer_expansion}) into (\ref{exp:hf_drifting_local_multiple_scale}) and get the same approximation as (\ref{sol:hf_drifting_local}) in the end. For $\lambda \le \frac{1}{6}$, the frequency is treated as $O(1)$ and is not considered here.

Comparing (\ref{sol:hf_drifting_local}) with (\ref{sol:slow_vary_eq_2}), and recalling the relationship $\Omega=\mu^{-\lambda}$, 
we find that the position of the tipping point is determined by $a$ that satisfies $\mbox{Ai}[(a-\frac{A^2}{2\Omega^2})/\mu^{2/3}]=0$. 
Then the tipping point with high frequency forcing and slowly drifting bifurcation parameter is given by
\begin{eqnarray}
a_{\rm hf} \sim a_d + a_p = \mu^{2/3}K + \frac{A^2}{2\Omega^2}\, , 
\label{ad+ap}
\end{eqnarray}
for $K$ given in (\ref{sol:location_tipping}) and  
 $a_p$ giving the same shift due to the periodic forcing 
as observed for the bifurcation point in (\ref{sol:hf_equili_local}) 
for $\mu=0$.  Recall that in this section we consider the case $A= o(\Omega$), with larger values of $A$ discussed below.

\begin{figure}[h]
\begin{center}$
\begin{array}{cc}
\includegraphics[width=0.45\textwidth]{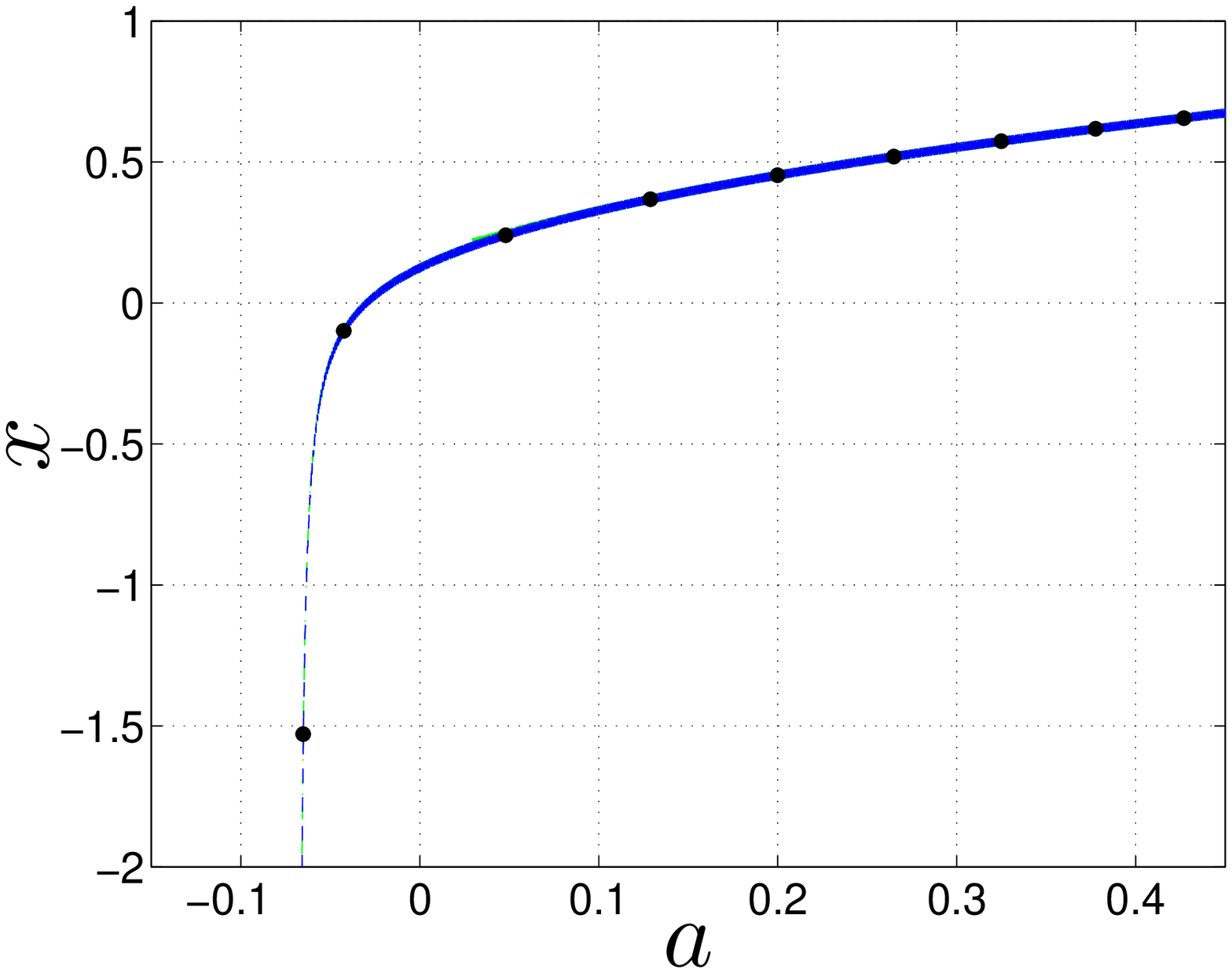} &
\includegraphics[width=0.45\textwidth]{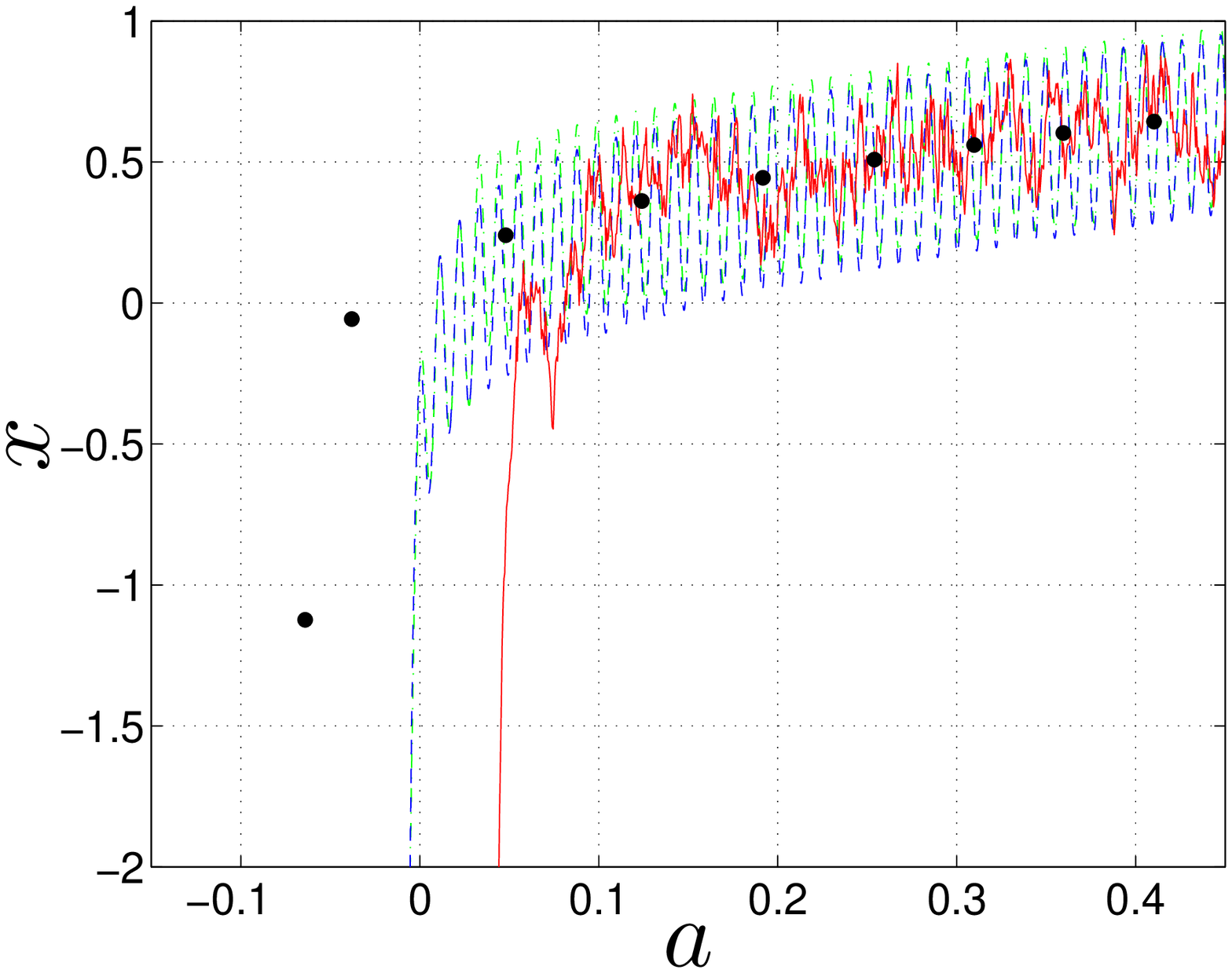}
\end{array}$
\end{center}
\caption {Tipping is shown for $\mu=.005$. LEFT: The asymptotic approximations of the attracting solution for $A=1$ and $\lambda =.8$, (\ref{sol:hf_drifting_outer}) (including higher order corrections up to $O(\mu^{3\lambda})$) and (\ref{sol:hf_drifting_local}) (blue dashed line) are consistent with the trajectory generated numerically (green dash-dotted line). The position of the tipping point is compared with the case of $A=0$ (black heavy dots). RIGHT: Blue dashed and green dash-dotted as on the LEFT with
$A=1$ and $\lambda = .2$, compared with 
the numerical simulation for the stochastic process obtained by replacing the periodic forcing by an additive white noise 
(red solid line). The noise coefficient is $\epsilon=0.2$ for the realization shown.} 
\label{fig:hf_trajectory}
\end{figure}

Figure \ref{fig:hf_trajectory} shows the effect of the high-frequency oscillation
in (\ref{sys:gen_periodic_forcing}) on the location of the tipping point. 
For higher frequency the location is almost the same as the case with no oscillation
but for lower frequency (smaller values of $\lambda$)
there is earlier tipping. We also see that a slower oscillatory forcing in (\ref{sys:gen_periodic_forcing}) results in a larger amplitude oscillation around
the slowly varying equilibrium solution found without oscillatory forcing (\ref{sol:slow_vary_eq_1}), consistent with the expression for the outer approximation (\ref{sol:hf_drifting_outer}).


In \cite{thompson2011, berglund2010}, it was shown that early tipping can be triggered also by white noise, in the system
obtained by replacing the periodic forcing in (\ref{sys:gen_periodic_forcing}) with an additive white noise, specifically,
 $dx = (a-x^2)\,dt + \epsilon dW(t)$ for $W(t)$ a standard
Brownian motion. Early tipping is shown to occur with non-negligible
  probability for noise
amplitude $\epsilon=O(\mu^{1/2})$ or larger. For the shifted 
 tipping points shown in Figure \ref{fig:hf_trajectory} (RIGHT) resulting from noisy and periodic forcings,
the amplitude of the oscillatory forcing ($A=O(1)$) is an order of magnitude larger than the noise coefficient ($\epsilon =o(1)$) that generates a comparable advance of tipping.
The probability of early escape, corresponding to
$x$ crossing the unstable branch, is in the range  $.25<P(x<-\sqrt{a})<.43$ 
for $a$ varying between $a=.05$ and $a=.025$ for the parameters shown here. 
  An asymptotic approach for the corresponding time-dependent
 probability density function is given in \cite{Zhu2014}. 

\begin{figure}[h]
\begin{center}$
\begin{array}{cc}
\includegraphics[width=0.45\textwidth]{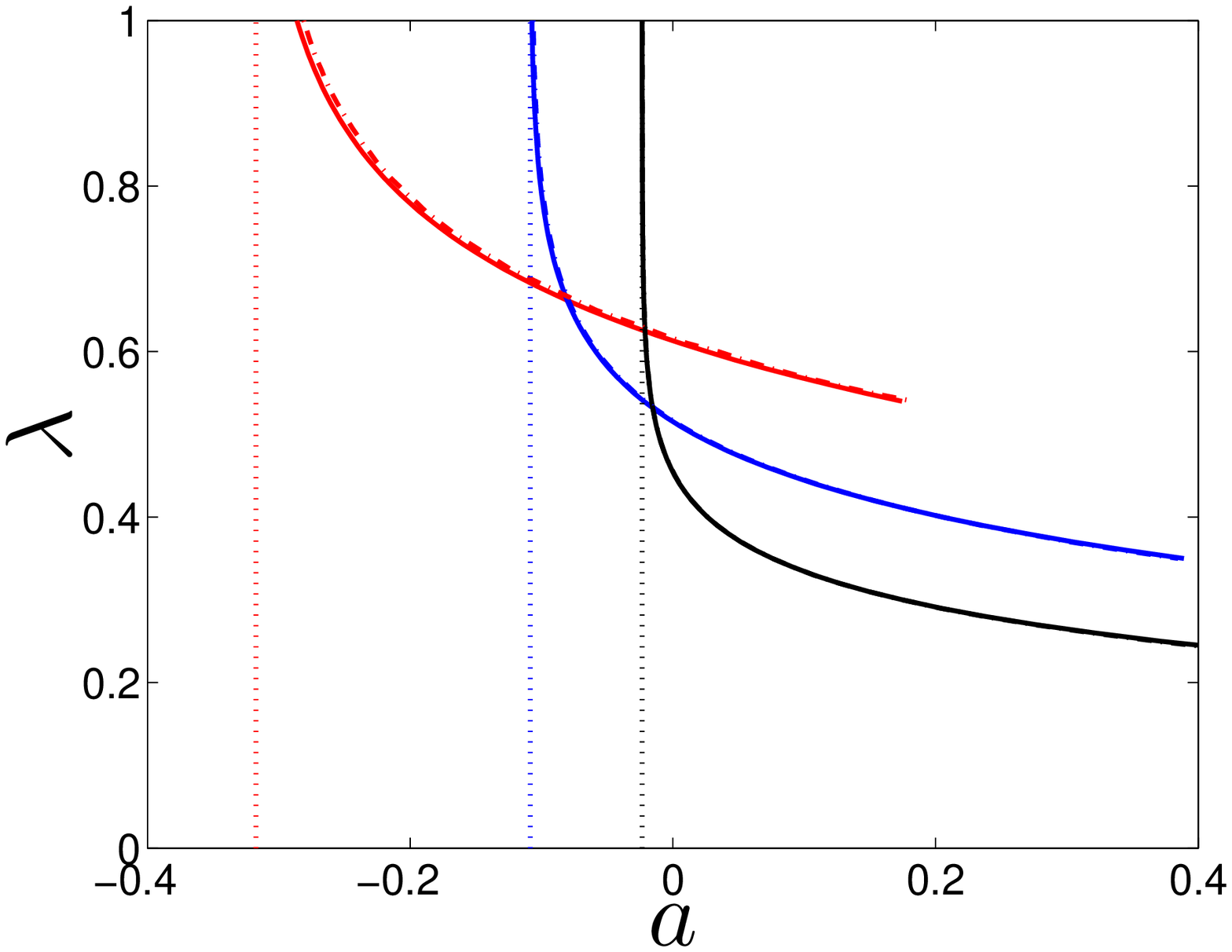} &
\includegraphics[width=0.45\textwidth]{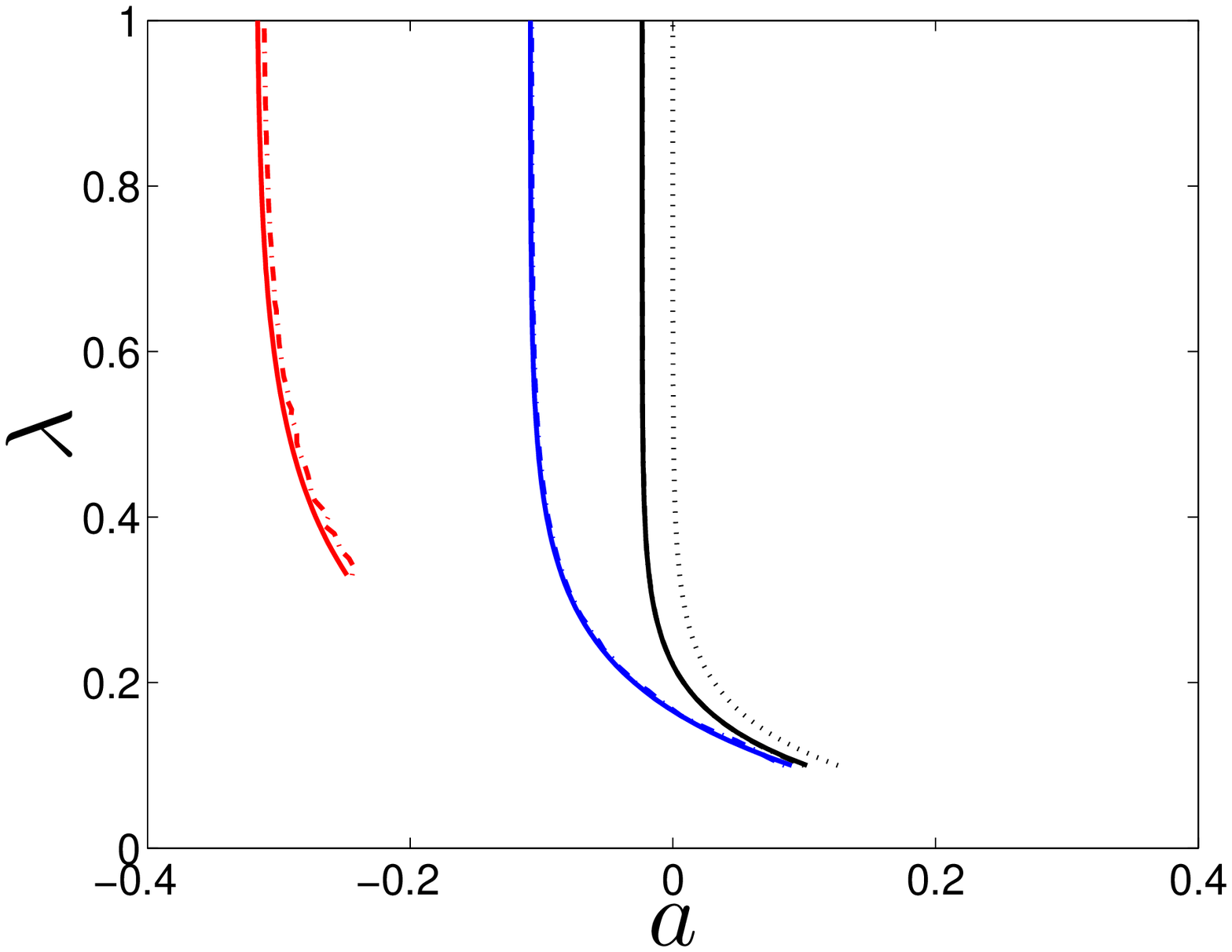}
\end{array}$
\end{center}
\caption{The position of the tipping point for different frequencies $\Omega = \mu^{-\lambda}$, shown for three drifting rates $\mu$. At $\lambda=1$ the values of $\mu$ shown
decrease from left to right in both panels:  $\mu=0.05$ in red; $\mu=0.01$ in blue; $\mu=0.001$ in black. The asymptotic approximation (\ref{ad+ap}) is shown by
the solid line.  
The dash-dotted line shows the location of the tipping point
obtained from numerical simulations, defined as the value of $a$ for $x=-10$.  LEFT: $A=5$; the dotted vertical lines show the position of the tipping point $a_d$ for $A=0$; RIGHT: $A=1$; the black dotted line shows the approximation of the bifurcation point $a_p = A^2/(2\Omega^2)$ for $\mu=0$.}
\label{fig:hf_tipping_position}
\end{figure}

In contrast, for  
high frequency periodic forcing, it is the ratio of amplitude to frequency that
plays a significant role in shifting the tipping point location, rather than 
the magnitude of the amplitude alone.  
The asymptotic approximation
for the tipping point location $a_{\rm hf}$ given in (\ref{ad+ap}) is
compared with numerical simulations in Figure \ref{fig:hf_tipping_position}
for different drifting rates $\mu$ and varying exponent $\lambda$ for $\Omega=\mu^{-\lambda}$.
There we see the competition between the two components in $a_{\rm hf}$.
For $A=1$ and smaller values of $\mu$, that is slower drifting rates, the location
of the tipping point is dominated by $a_d$ corresponding to
  $A=0$ (\ref{sol:location_tipping}). In those cases only for smaller values of $\lambda$, that is lower frequencies, is there a noticeable influence of the periodic
term, as seen in  Figure \ref{fig:hf_tipping_position}(RIGHT).
The advance of $a_{\rm hf}$
 due to the periodic forcing as captured by $a_p$ clearly increases
with $A$. For example, the case $A=5$ is shown in Figure \ref{fig:hf_tipping_position}(LEFT), 
  where there are larger advances associated with $a_p$ for decreasing $\lambda$ 
as compared with negligible advances for $A=1$ and smaller $\lambda$.
  As noted above, for smaller values of $\lambda$ (near $\lambda \approx 1/3$ or below), 
 the frequency $\Omega$ of the oscillatory forcing is
not necessarily large, depending on the value of $\mu$. 
Then the expansion (\ref{exp:hf_outer_expansion}) is no longer a reasonable asymptotic approximation, 
and instead the behavior has some of the features discussed in Section \ref{sec:lf}
 for low frequency forcing.  As seen
in Figure \ref{fig:hf_tipping_position}, the asymptotic approximation is valid
 for a larger range of $\lambda$ for  smaller values of $\mu$.
The approximation $a_{\rm hf}$ and (\ref{sol:hf_drifting_local}) break down 
for $A \sim O(\Omega)$, since the
 coefficient $\mu^{2\lambda} A^2 = A^2/\Omega^2$ that appears in 
(\ref{exp:hf_drifting_in_order1})-(\ref{sol:hf_drifting_local}) 
and in the argument of ${\rm Ai}$ is treated as $o(1)$ to obtain 
$a_{\rm hf}$.  We consider results for larger values of $A$ separately in the next section.

\subsection{High Frequency Oscillation with Large Amplitude ($A \gg 1$)}
\label{sec:hf_large_A}

In Section \ref{sec:hf_no_drifting} and \ref{sec:hf_slow_drifting}, the asymptotic approximations (\ref{sol:hf_drifting_outer}) and (\ref{sol:hf_drifting_local}) are not valid
if $A=O(\Omega)$ or larger. However, for some larger values $A$
and frequency $\Omega$, we can rescale (\ref{sys:gen_periodic_forcing}) to
obtain a similar system forced by an oscillation with amplitude of unity.  In
certain cases we can apply the approximations from 
either Section \ref{sec:hf_slow_drifting} or Section \ref{sec:lf_order1} to the rescaled
system.

Substituting $x = \sqrt{A}z$ and $t = {S}/\sqrt{A}$ into (\ref{sys:gen_periodic_forcing}), we get

\begin{eqnarray}
 & & \frac{\mbox{d}z}{\mbox{d}S}=h-z^2+\mbox{sin}(\omega S), \nonumber\\
 & &\frac{\mbox{d}h}{\mbox{d}S}=-M, \nonumber\\
 & &v(0)=v_0=\frac{x_0}{\sqrt{A}}, \qquad  h(0)=h_0=\frac{a_0}{A},
\label{sys:norm_periodic_forcing}
\end{eqnarray}
where $\omega={\Omega}/{\sqrt{A}}$, $h = \frac{a}{A}$, and
$M = A^{-3/2}\mu$  are the scaled
frequency, bifurcation parameter, and drifting rate, respectively.
In order to consider the case $A \gg 1$, we assume $A = \Omega^P$ with $P \ge 1$. 
Then the relationship between $\omega$ and $M$ is given by 
$\omega = M^{\zeta}$, with $\zeta$ a function of $\lambda$ and $P$
\begin{equation}
\zeta = \frac{{P}-2}{{3}P+{2}/{\lambda}}. \label{exp:hf_drifting_norm_gamma}
\end{equation}

For $1 \le P < \frac{4}{3}$, the high frequency approximation (\ref{sol:hf_drifting_local}) can be applied to (\ref{sys:norm_periodic_forcing}) if $\zeta < -\frac{1}{6}$, which indicates that $\lambda$ in $\Omega = \mu^{-\lambda}$ satisfies $\lambda > \frac{2}{12-9P} \ge \frac{2}{3}$. For larger values of
 $A$ ($P \ge \frac{4}{3}$), the oscillatory forcing and the drifting bifurcation parameter both dominate the dynamics to leading order, so
that the  asymptotic approximation (\ref{sol:hf_drifting_outer}) dominated
by the slowly drifting parameter $a(t)$ is no longer valid.

For $P > 2$ and $\lambda>0$, $\zeta$ is positive which indicates that even though (\ref{sys:gen_periodic_forcing}) has a high-frequency oscillatory forcing, in the normalized system (\ref{sys:norm_periodic_forcing}) the relationship between $\omega$ and $M$
corresponds to a low-frequency case, discussed in the next section. We note that $\zeta$ is a monotone increasing function with respect to $P$ with a horizontal asymptote at $\zeta=\frac{1}{3}$.

\section{Low Frequency Oscillation $\Omega \ll 1$}
\label{sec:lf}

For the case of a low-frequency oscillation ($\Omega \ll 1$), we introduce a new parameter $\nu>0$ by writing $\Omega = \mu^{\nu}$, capturing the relationship between the frequency and the drifting rate. In the following we consider two different cases: 
$\nu \geq 1$ so that $\Omega = O(\mu)$ or smaller, or $\Omega = \mu^{\nu}$ with $0<\nu<1$. Unlike the high frequency case,
for low frequency oscillations we can not use the method of multiple scales, since the key time scales of drifting
and the frequency are both slow. Then the asymptotic approach relies on understanding the local behavior of the trajectory that
is driven by the slow oscillations.

\subsection{$\Omega = O(\mu)$ or smaller}
\label{sec:lf_order1}

Here it is convenient to  write $\Omega = c\cdot \mu$ for $c$ a positive constant
and to write the system in terms of  $\tau = \mu t$, expressing
  $x$ as a function of  $a(\tau)=a_0-\tau$ in (\ref{sys:gen_periodic_forcing}), to get
\begin{equation}
-\mu\frac{dx}{da} = f(a)-x^2=a-x^2+A \mbox{sin}(c \cdot (a_0-a)). \label{exp:lf_multi_scale}
\end{equation}
Substituting $ x \sim x_0 + \mu x_1 + \mu^2 x_2+ \cdots $
into (\ref{exp:lf_multi_scale}), we find the
first two terms of the approximate solution as a function of $a(\tau)$
\begin{align}
x &\sim \pm \sqrt{f(a)}+ \mu \cdot \frac{f'(a)}{4f(a)} +\cdots \nonumber \\
&\sim \pm \sqrt{a+A\mbox{sin}(c \cdot (a_0-a))}  + \mu \cdot \frac{1-c \cdot A\mbox{cos}(c \cdot (a_0-a))}{4(a+A\mbox{sin}(c \cdot (a_0-a)))}+\cdots, \label{sol:lf_outer}
\end{align}
where the positive sign corresponds to the solution that is attracting. 
 Note that the oscillatory term appears in the leading order contribution in  (\ref{sol:lf_outer}), in contrast to the leading order approximation in the case of the high frequency oscillations (\ref{sol:hf_drifting_outer}) where there is no oscillatory term.  Similar to (\ref{sol:slow_vary_eq_1}), the
terms shown in (\ref{sol:lf_outer}) are both $O(\mu^{1/3})$
for $f(a) \sim O(\mu^{2/3})$ and the expansion is no longer valid. Then we expect a tipping point near $a=a_r$, where $f(a)>0$ for $a>a_r$ and $f(a_r)=0$.  A local analysis near $a=a_r$ is used to determine the location of the tipping point.
 In our analysis and simulations we take $a_0>a_r$ for $f(a_r)=0$ with $|a_0 -a_r|= O(1)$, ensuring 
that the system (\ref{sys:gen_periodic_forcing}) attracts to the outer solution (\ref{sol:lf_outer}) 
prior to any tipping.

\begin{figure}[h]
\begin{center}$
\begin{array}{cc}
\includegraphics[width=0.45\textwidth]{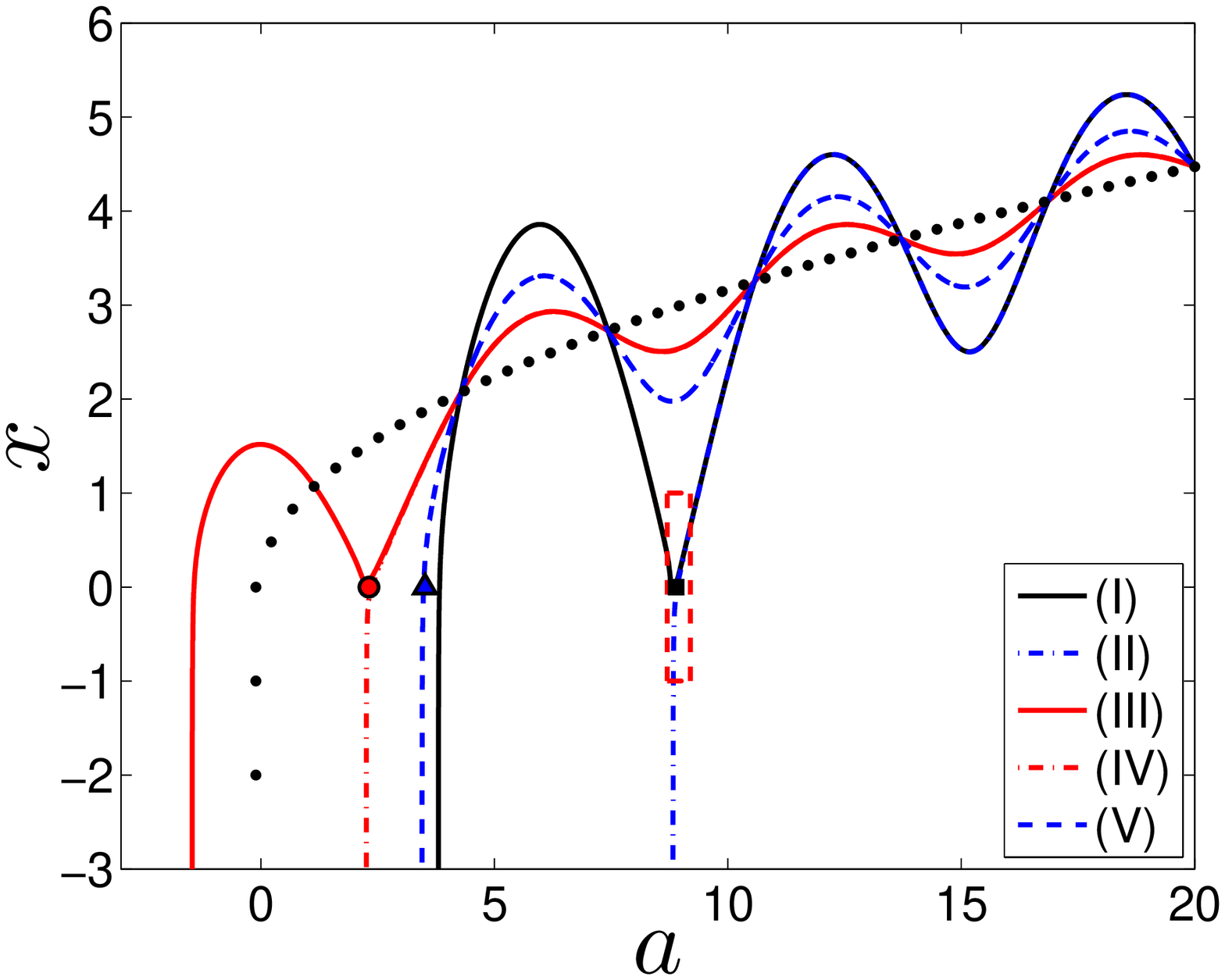}
\includegraphics[width=0.45\textwidth]{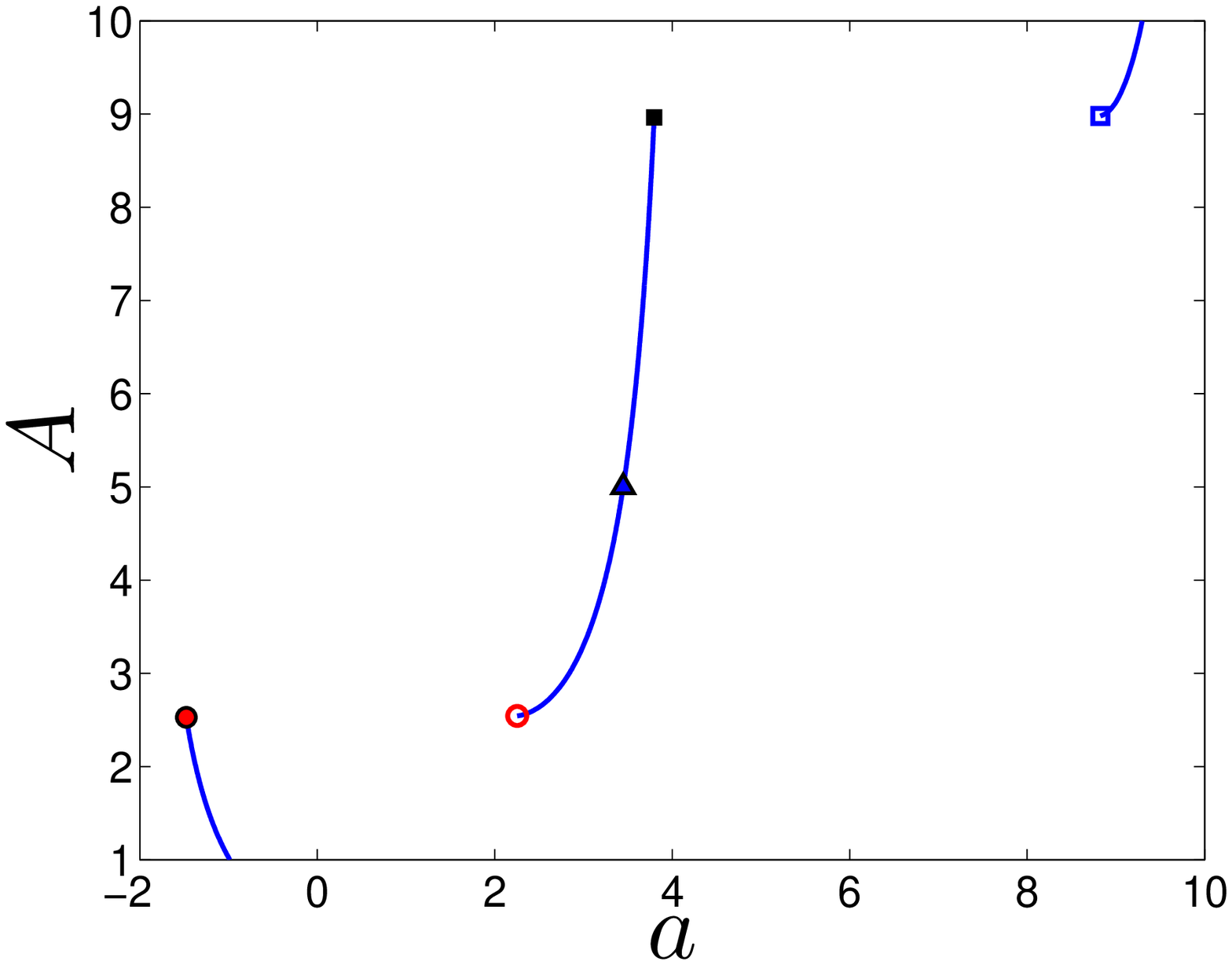}&
\end{array}$
\end{center}
\caption{ Parameter values are $\mu=\Omega=0.01$($c=1$), $a_0=20$ for all graphs. 
LEFT: Numerical simulation of (\ref{sys:gen_periodic_forcing}).  The Roman numerals correspond to cases with different values of
$A$: Case (I), $A=8.9647$; Case (II), $A=8.9797$;  Case (III), $A=2.5268$; Case (IV), $A=2.5418$;  Case (V), $A=5$. With the exception of Case (V), $A$ is near the critical value $A_c$ at which there is
an abrupt change of the tipping position. Heavy black dots: $A=0$.  The solid black square  (red circle) corresponds to parameter values $(a,A) =(a_m,A_m)$,
in Cases (I) and (II) ((III) and (IV)), with $A_m$ taking values slightly less than
the value of $A$ for Cases (I)-(IV).  The dashed rectangle indicates the region 
shown on a larger scale in
  Figure \ref{fig:lf_threshold}(LEFT). 
 The blue triangle marks the value of $a_r$ where $f(a_r)=0$ for $A=5$.
RIGHT: The location of the tipping point as a function of $A$ in (\ref{sys:gen_periodic_forcing}) obtained by numerical simulations, with tipping defined as the value of $a$ where $x=-10$ as in the previous cases. 
 Different markers indicate the value of A and corresponding tipping point value, for the trajectories shown in the LEFT figure. Solid (black) square: (I); Open (blue) square: (II); Solid (red) circle: (III);
Open (red) circle: (IV);  Solid (blue) triangle: (V).}
\label{fig:lf_LO}
\end{figure}
\begin{figure}[h]
\begin{center}$
\begin{array}{cc}
\includegraphics[width=0.45\textwidth]{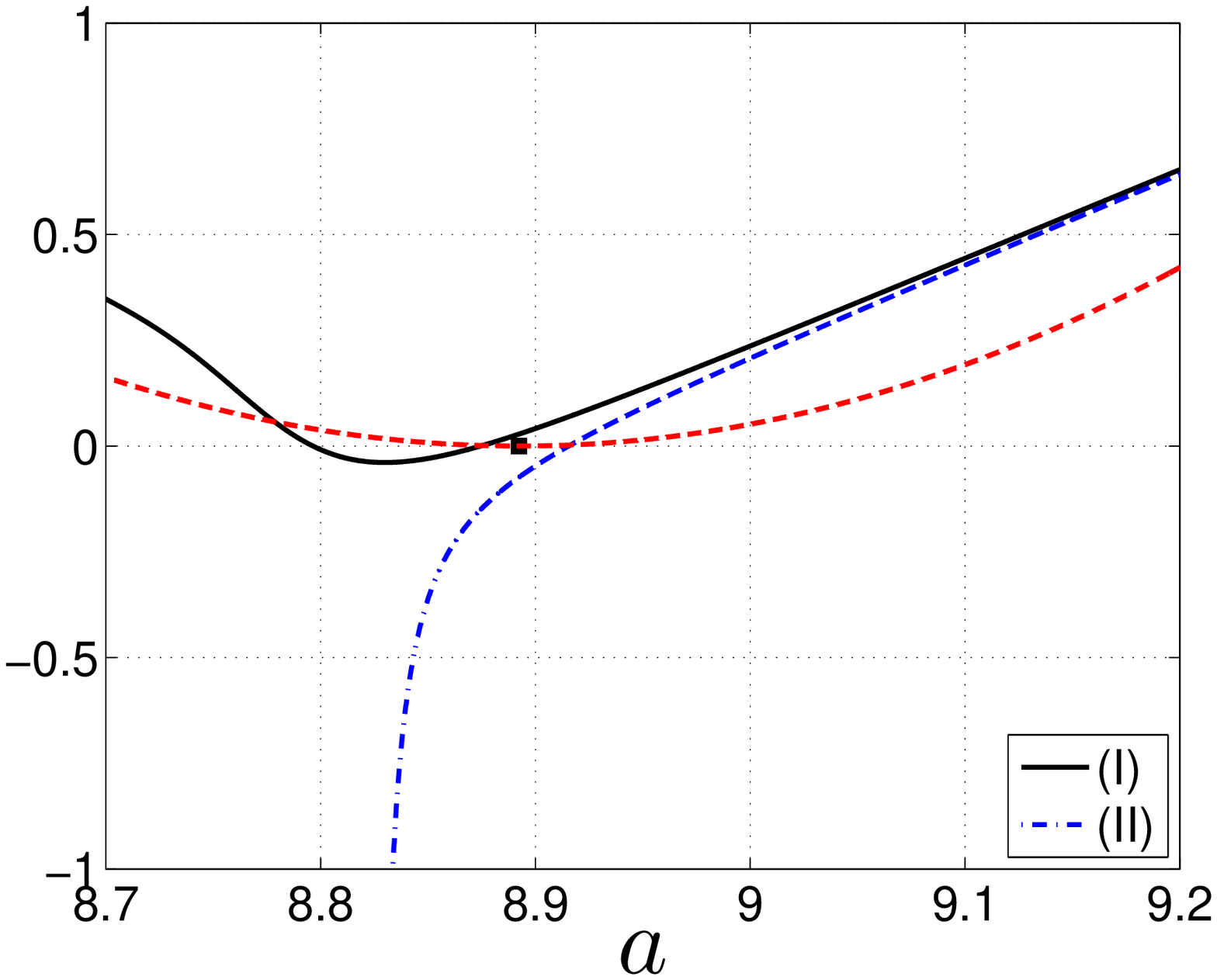}&
\includegraphics[width=0.45\textwidth]{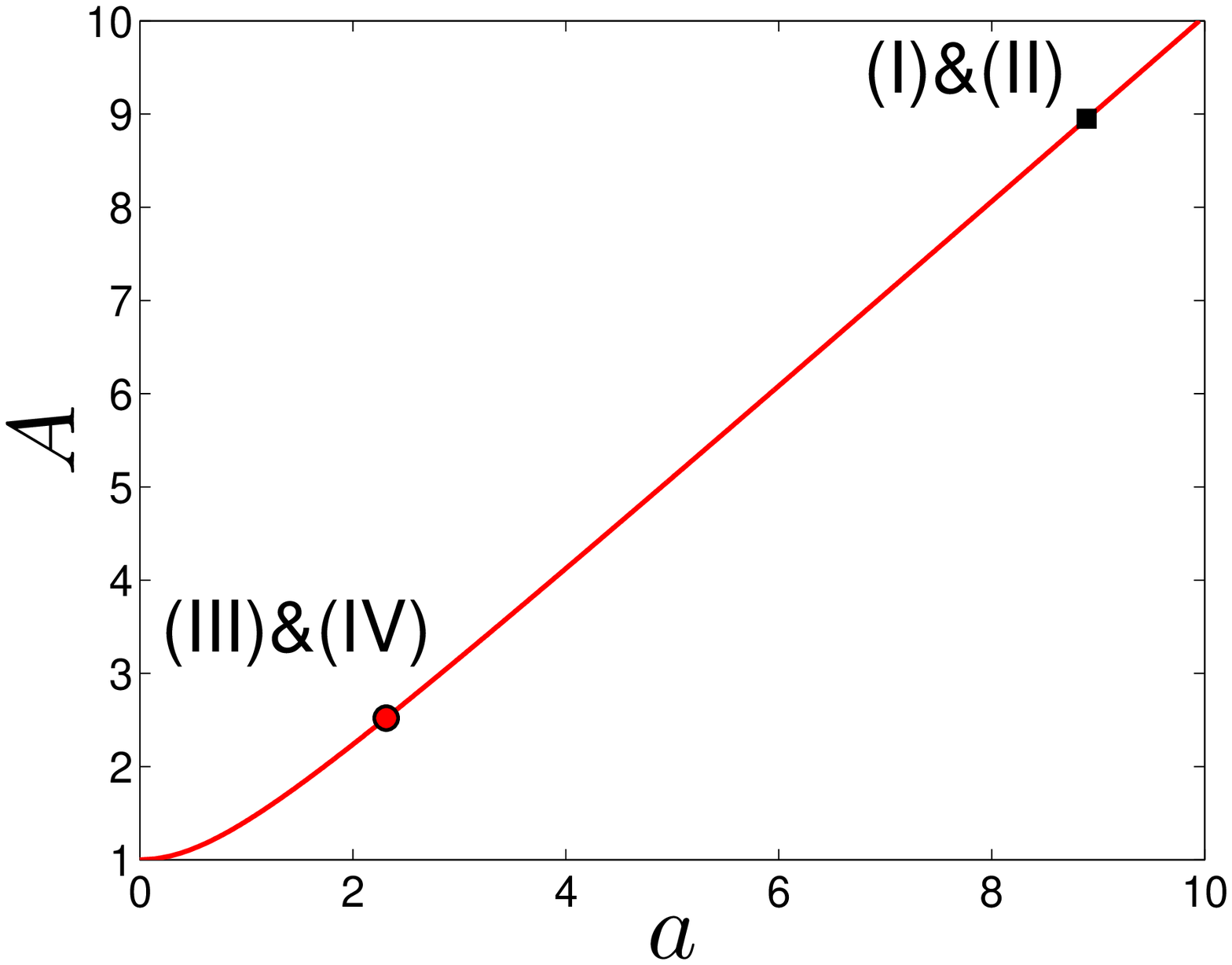}
\end{array}$
\end{center}
\caption{For all graphs, $\mu=\Omega=0.01$ ($c=1$). LEFT: Figure \ref{fig:lf_LO}(LEFT) zoomed-in near 
$a_m$ for $f'(a_m)=0$ (black square), with $f(a)$ for $a_0=20$ and $A=A_m=8.9485$ given by the
   red dashed curve. 
Numerical simulations of Case (I) (solid black line,  $A= 8.9647$) and Case (II) (dash-dotted blue line, $A=8.9797$) are also shown.  RIGHT: 
 $A_m$ as a function of $a_m$ (solid red line), with $a_0$ 
varying along this line as in
(\ref{exp:lf_threshold_1}) and (\ref{exp:lf_threshold_2}). The black square (red circle)  marks 
the value of $a_m$ and $A_m$ for Cases (I) and (II) (Cases (III) and (IV)), with $a_0=20$ as in Figure \ref{fig:lf_LO}.  }
\label{fig:lf_threshold}
\end{figure}

Figure  \ref{fig:lf_LO} shows numerical simulations of 
(\ref{sys:gen_periodic_forcing}) for different values of $A$ and fixed
$a_0$, illustrating important characteristics of tipping that are the basis for our analysis.
 The heavy black dots correspond to $A = 0$ (no modulations). For larger values of $A$,
 the trajectory $x = x(a)$ exhibits larger oscillations and exhibits 
tipping before $a$ reaches or crosses the static bifurcation value at $a=0$. 
As could be expected, the value of $a$ where the tipping occurs increases with $A$ for $a>0$.  
For amplitudes $A$ where the trajectory crosses $x = 0$ for $a(t)>0$, tipping can occur well before the
static bifurcation point is reached.  Abrupt shifts in tipping point location are observed for critical 
amplitudes denoted $A_c$ that are slightly greater than $A_m$, where
the pair $(a,A)=(a_m,A_m)$ satisfies $f(a_m)=f'(a_m)=0$ as discussed below (\ref{exp:lf_threshold_1}). 
For example,  the Cases (I) and (II) shown in Figure
 \ref{fig:lf_LO}(LEFT) have the values of  $A = 8.9647$ and $A = 8.9797$, respectively, with both values slightly greater
than $A_m \approx 8.9485$ for $a_0=20$. However, the tipping for Case (I) occurs at  
$a\approx 3.80$,  (solid black square in  Figure
 \ref{fig:lf_LO}(RIGHT)),
while the tipping for Case (II) occurs at $a \approx  8.82$ (open blue square in  Figure
 \ref{fig:lf_LO}(RIGHT)).  Here the value of $a_m$ is
indicated by the solid black square in 
Figure  \ref{fig:lf_LO}(LEFT), which is near the tipping point for Case (II) but not for Case (I).  This figure illustrates how amplitudes $A>A_c$ cause an early tipping near a value of
$a$ where $x<0$, while for the
same initial condition and $A<A_c$, the periodic forcing is not sufficient to drive tipping at
the same value of $a$ as for $A>A_c$.  Rather, the oscillations continue for another period with
tipping likely to take place near the next crossing of $x=0$.
Figure  \ref{fig:lf_LO}(RIGHT) shows the location of the tipping points for all
trajectories represented in Figure \ref{fig:lf_LO}(LEFT), and indicates the critical values $A_c$ where there
are abrupt changes in the tipping point as a function of $A$. Note that the tipping points for
$A<A_c$ and $A>A_c$ may both correspond to early tipping relative to the static bifurcation
$a=0$ and to the unforced delayed bifurcation $a=a_d$ in (\ref{sol:location_tipping}), as 
is true for Cases (I) and (II).
For smaller values of $A_c$, such as the value near Cases (III) and (IV), 
amplitudes $A>A_c$ drive
early tipping at $a>0$ (Case (IV)), while $A<A_c$ drives tipping at $a<0$, delayed relative to $a=0$ and $a=a_d$ (Case (III)).  As we study separately below, the values $A_c$ are near $A_m$, and are related to the local behavior of the 
trajectory near $a_m$.
  Figure \ref{fig:lf_threshold}(LEFT) 
zooms in on  the trajectories (I) and (II) near $a = 8.9$ and illustrates the obvious change
of concavity near $a_m$, which proves to be an essential part in estimating $A_c$. 

 Before we investigate the solution near $(a_m,A_m)$ where $f(a_m)=f'(a_m)=0$ (and near $A_c$), we first consider
the case where tipping occurs close to $a = a_r$
where $f(a_r)=0$ and $f'(a_r) = O(1)$, as in case (V) in Figure \ref{fig:lf_LO}(LEFT). 
 Near $a=a_r$, specifically $f(a)=O(\mu^{2/3})$, 
 we substitute $a-a_r = \mu^{2/3}B$ and $x = \mu^{1/3}X$ into (\ref{exp:lf_multi_scale}), yielding
\begin{equation}
O(\mu^{2/3}): -X_{B} = f'(a_r)B - X^2 = [1-c A \mbox{cos}(c\cdot (a_0-a_r))]B -X^2, \label{exp:lf_local_IO_A}
\end{equation}
 which has the solution for $f'(a_r)\neq 0$,
\begin{equation}
X = -C^{1/3}\frac{\mbox{Ai}'(C^{1/3}B)}{\mbox{Ai}(C^{1/3}B)}, \,\,\, \Rightarrow \,\,\, x = -(\mu C)^{1/3}\frac{\mbox{Ai}'((C/\mu^2)^{1/3}(a-a_r))}{\mbox{Ai}((C/\mu^2)^{1/3}(a-a_r))}, \qquad  C=f'(a_r). 
\label{sol:lf_local_A}
\end{equation}
Similar to the analysis for (\ref{sol:slow_vary_eq_2}) and (\ref{sol:hf_drifting_local}), the location of the tipping point is determined by the singularity of (\ref{sol:lf_local_A}) where $\mbox{Ai}((C/\mu^2)^{1/3}(a-a_r))=0$, which is
\begin{equation}
a_{lf} = a_r +  a_d/(f'(a_r))^{1/3} \mbox{ {\rm for} } f(a_r)=0 \mbox{ and } f(a)>0 
 \mbox{ {\rm for} } a>a_r,
 \label{sol:lf_tipping_position}
\end{equation}
where $a_d$ is given in (\ref{sol:location_tipping}). 
We differentiate $f(a_r)=0$ with respect to $A$, to obtain
$da_r/dA = (CA)^{-1}a_r$, which indicates that the location of the tipping point increases (decreases) with increasing $A$ for $a_r$ negative (positive), as shown in the Figures \ref{fig:lf_LO}(RIGHT) and
\ref{fig:lf_amplitude_threshold}(RIGHT).

Noting that (\ref{sol:lf_local_A}) is obtained from  (\ref{exp:lf_local_IO_A}) for $f'(a_r)\neq 0$,
  we propose a
separate analysis for $(a,A)$ close to $(a_m, A_m)$, which
 satisfy the conditions,
\begin{align}
f(a_m)&=a_m+A_m \mbox{sin}(c \cdot (a_0-a_m)) = 0, \label{exp:lf_threshold_1} \\
f'(a_m)&=1-c A_m \mbox{cos}(c \cdot (a_0-a_m)) =0, \,\,\, \Rightarrow \,\,\,c^2 a_m^2+1 =c^2 A_m^2. \label{exp:lf_threshold_2} 
\end{align}
The behavior of $f(a)$ near $a_m$ for $A=A_m$ is shown in Figure \ref{fig:lf_threshold}(LEFT).
Figure \ref{fig:lf_threshold}(RIGHT) shows $A_m$ as a function of $a_m$ for different values of $a_0$ 
as in (\ref{exp:lf_threshold_1}) and (\ref{exp:lf_threshold_2}).  From these conditions, we
see that for each $a_0$ there 
is a family of  initial values $a_0+\frac{2k\pi}{c}$ for $k$ an integer that corresponds to
the same pair $(a_m, A_m)$.  Equations (\ref{exp:lf_threshold_1}) and (\ref{exp:lf_threshold_2}) 
indicate that  $A_m > 1/c$ and $0 < a_m<A_m$.  We take $A>0$ in (\ref{sys:gen_periodic_forcing}), noting that similar results can be found for $A<0$.

\begin{figure}[h]
\begin{center}$
\begin{array}{cc}
\includegraphics[width=0.45\textwidth]{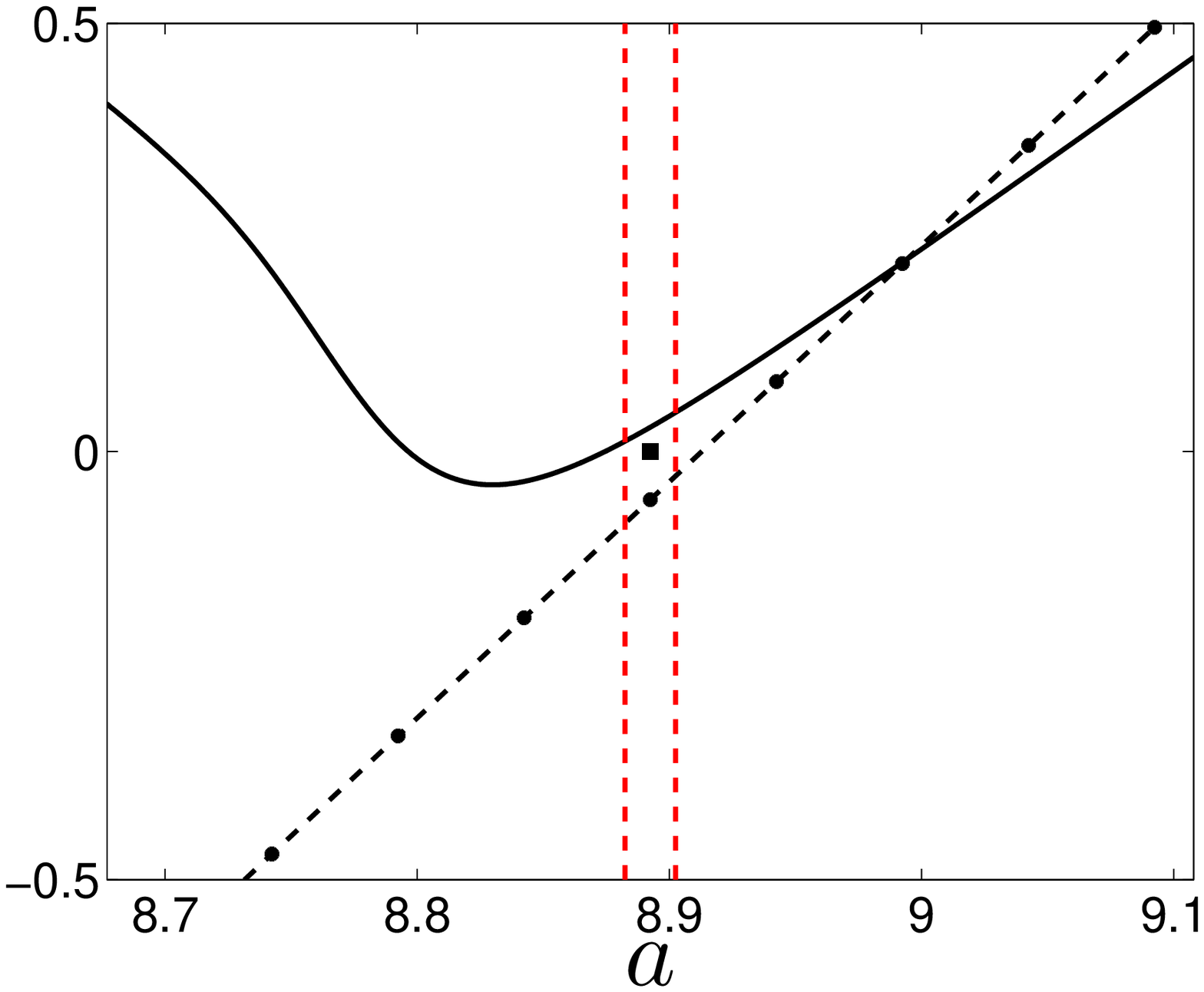}&
\includegraphics[width=0.45\textwidth]{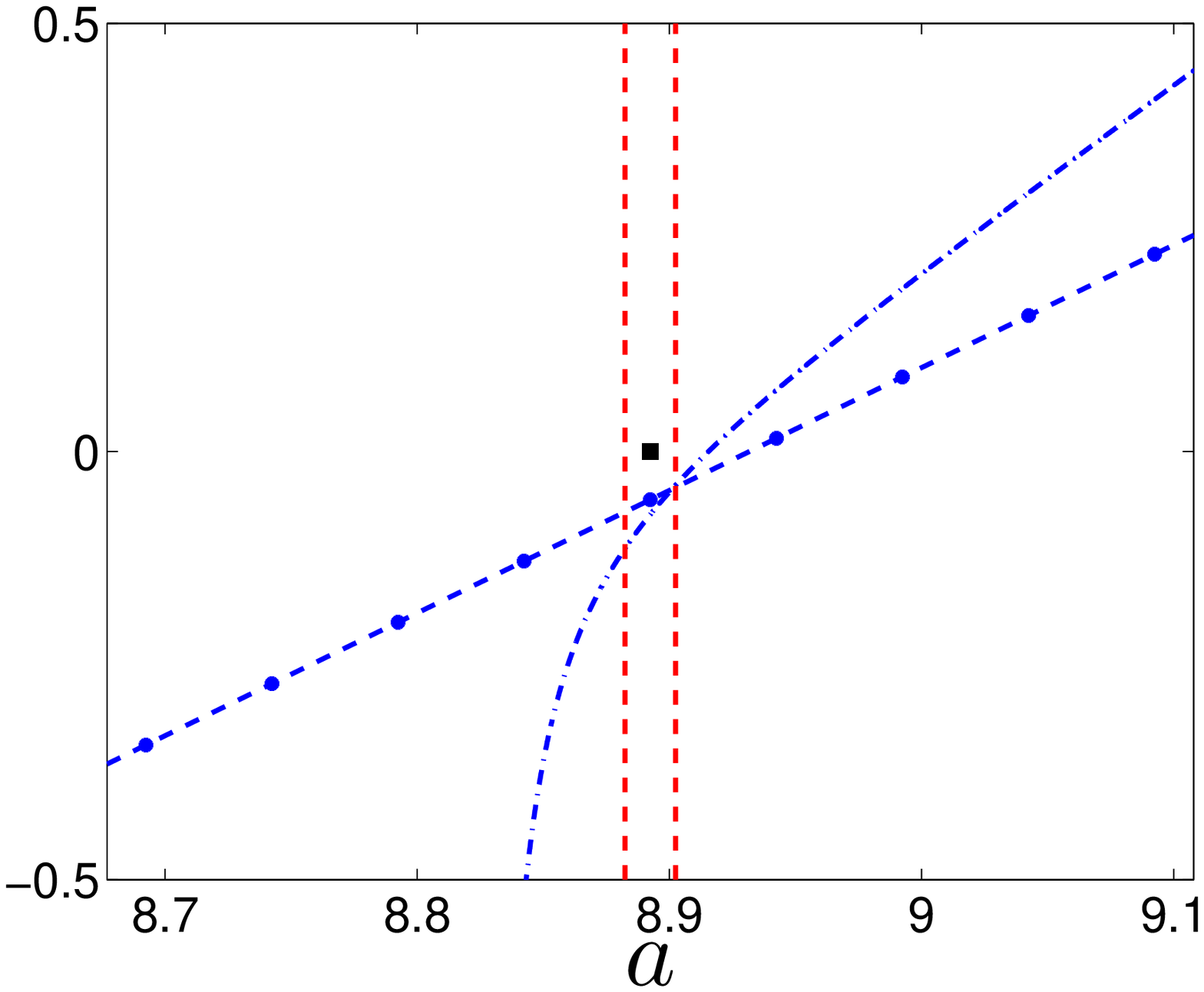}
\end{array}$
\end{center}
\caption{View of Cases (I) and (II) as in Figure \ref{fig:lf_LO}(LEFT), zoomed in on
the scale $|a-a_m|= \mu^{1/3}$, with  
 $|a-a_m| = \mu$ indicated by red dashed lines.
  The black square is the value of $a_m$ where $f(a_m)=f'(a_m)=0$. 
LEFT: The trajectory of Case (I) (solid black line) is compared 
with the concavity threshold (\ref{sol:lf_concavity_threshold}) (black dashed line with markers).
RIGHT: The trajectory of Case (II) (blue dash-dotted line) is compared
with the concavity threshold (\ref{sol:lf_concavity_threshold}) (blue dashed with markers).}
\label{fig:lf_local}
\end{figure}

\begin{figure}[h]
\begin{center}$
\begin{array}{cc}
\includegraphics[width=0.45\textwidth]{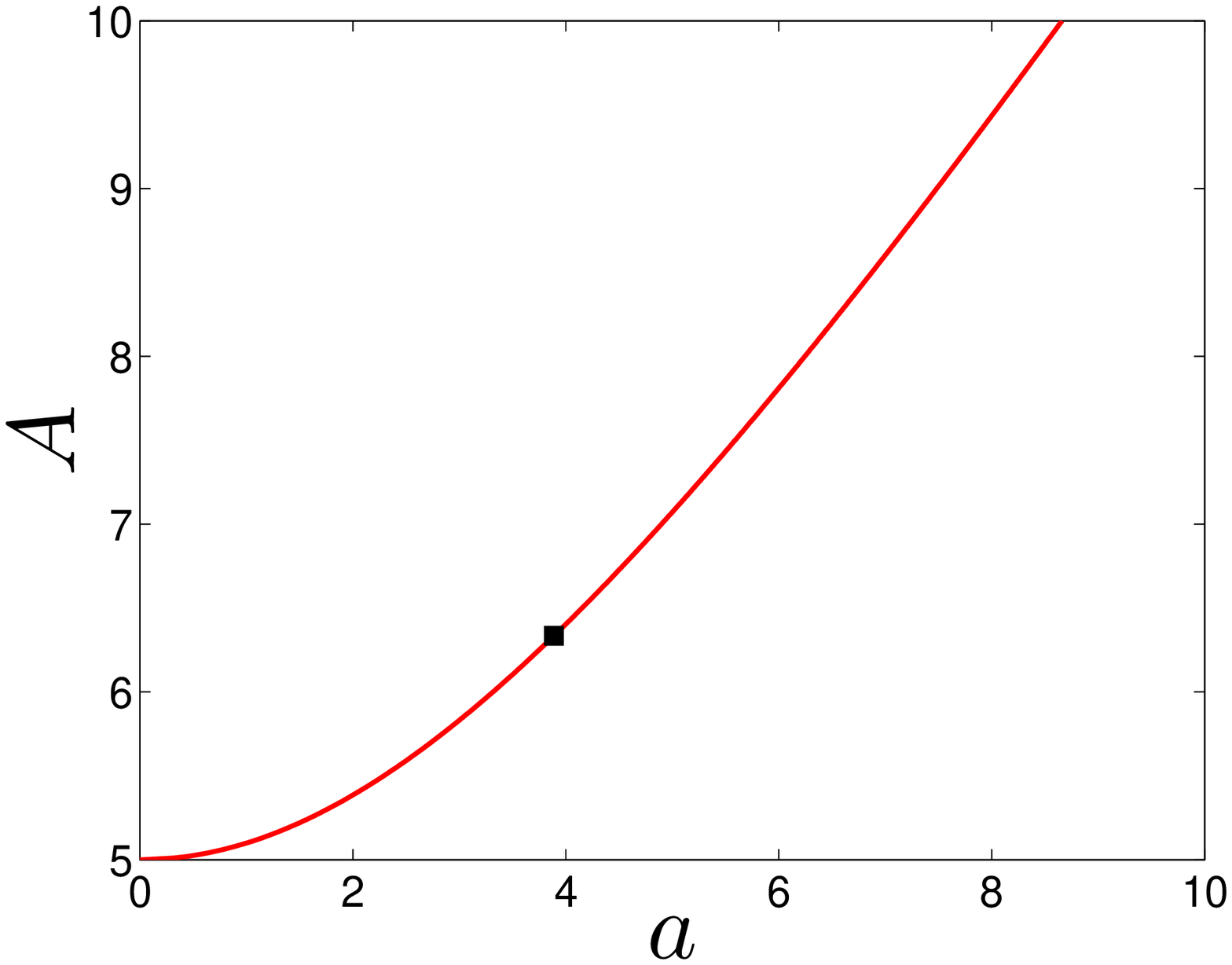} &
\includegraphics[width=0.45\textwidth]{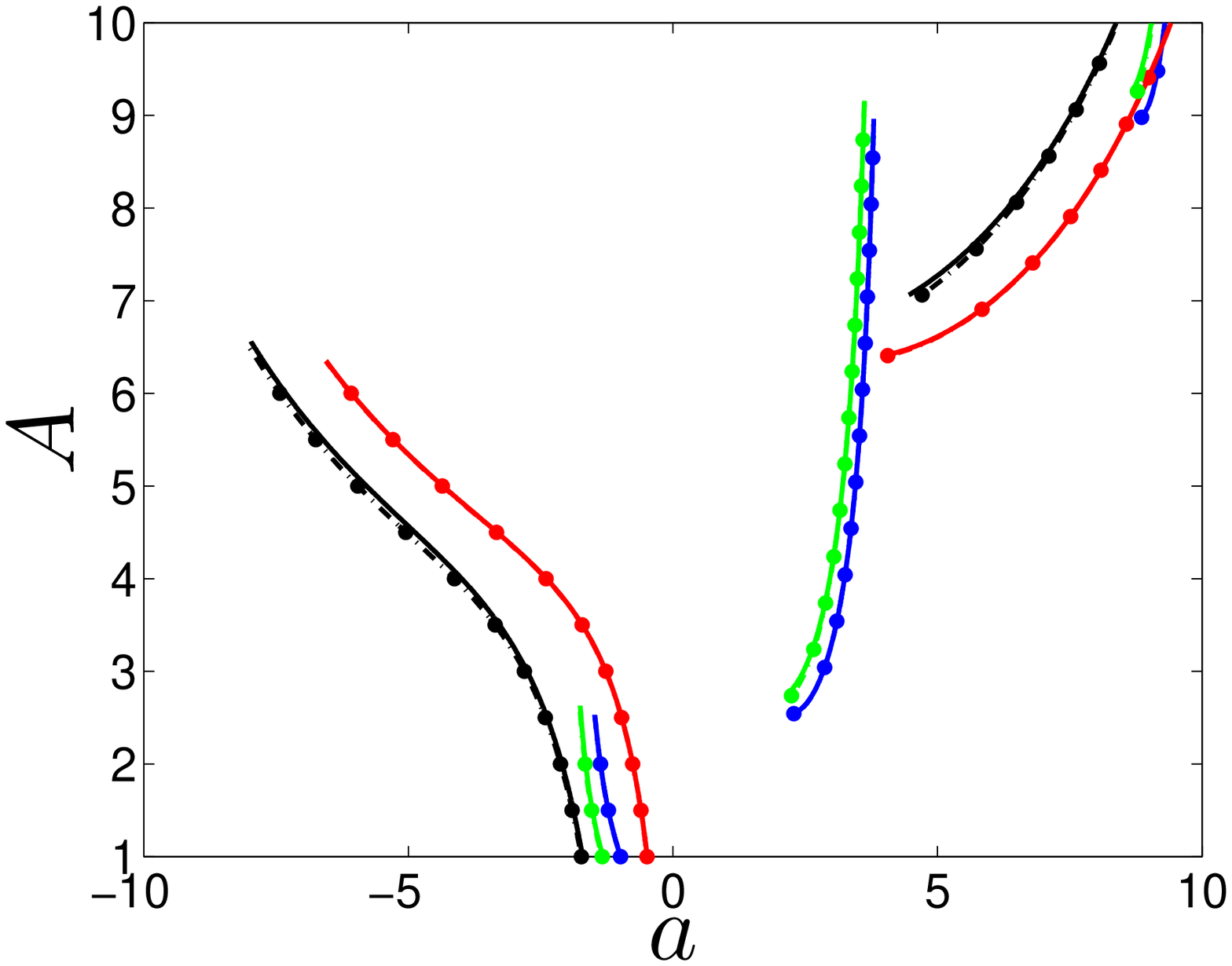}
\end{array}$
\end{center}
\caption{LEFT: Values of $A_m$ (red solid line) from (\ref{exp:lf_threshold_1}) and (\ref{exp:lf_threshold_2})  as a function of $a_m$ for $c = 0.2$. The black square marks the value of $a_m$ and $A_m$ near
parameter values where there is an abrupt change in the tipping point location, shown for the cases with
$c=.2$ in RIGHT. RIGHT: The position of the tipping point as a function of amplitude $A$, comparing the asymptotic approximation (\ref{sol:lf_tipping_position}) with jumps at values (\ref{sol:lf_A_threshold}) (solid lines) with numerical simulation, defining the tipping point as the value of $a$ when $x=-10$ (dash-dotted lines with markers). Different colours indicate different combinations of frequency $\Omega$ and drifting rate $\mu$. Blue: $\Omega=0.01$, $\mu=0.01$, and $a_0=20$; Red: $\Omega=0.01$, $\mu=0.05$, and $a_0=32$; Green: $\Omega=0.1$, $\mu=0.1$, and $a_0=20$; Black: $\Omega=0.1$ and $\mu=0.5$, and $a_0=32$.}
\label{fig:lf_amplitude_threshold}
\end{figure}

 A comparison of Cases (I) and (II) in Figure \ref{fig:lf_LO}(LEFT) with Figure \ref{fig:lf_threshold}(LEFT)
 suggests that the critical amplitude $A_c$, for which there is a jump in the location of the tipping point, is near $A_m$, with $A_c>A_m$.
 Furthermore, the trajectories shown in Figure \ref{fig:lf_threshold}(LEFT) have different concavities near 
$a=a_m$ for $A$ above and below $A_c$. Specifically, for amplitudes $A>A_c$
for which an early tipping occurs, the trajectory is concave down near its crossing of $x=0$.
In contrast, if $A<A_c$, the trajectory is concave up near its crossing of $x=0$ and no tipping
occurs there, yielding oscillations for another period with
tipping likely to take place near the next crossing of $x=0$.
Therefore, to find $A_c$ and the corresponding tipping point, we consider the concavity of $x$ in (\ref{exp:lf_multi_scale}) with respect to $a$ near $a_m$.
Using (\ref{exp:lf_multi_scale}), we write $x''(a)$ 
in terms of the  local variables 
\begin{eqnarray}
a-a_m = \mu \eta, \qquad  x = \mu \xi, \qquad  \mbox{ and }  A = A_m+\mu A_1 \, , \label{local:xietaAm} 
\end{eqnarray}
to focus on the neighbourhood of $a=a_m$ and on amplitude values near $A_m$.  
Substituting the local variables (\ref{local:xietaAm}) in (\ref{exp:lf_threshold_1}) and (\ref{exp:lf_threshold_2}), and expanding
 $\mbox{sin}(c\mu \eta)$ and $1-\mbox{cos}(c\mu \eta)$ for $\mu\ll 1$, we obtain the leading order expression for the local concavity of $\xi(\eta)$,
\begin{eqnarray}
\frac{d^2\xi}{d\eta^2} & =&  \mu\left( \frac{A_1}{A_m}-c^2 a_m \eta +\frac{2A_1 a_m}{A_m}\xi  \right) + O(\mu^2). 
\label{exp:lf_concavity_local_IO}
\end{eqnarray}
 From (\ref{exp:lf_concavity_local_IO})
 we find that the leading order expression of $\xi = x/\mu$ that
 separates trajectories of different concavities is given by
\begin{equation}
\xi_c = \frac{c^2 A_m }{2A_1}\eta-\frac{1}{2a_m} \, .\label{sol:lf_concavity_threshold}
\end{equation}
Then for $\xi > \xi_c $ ($\xi <\xi_c$), 
near $a=a_m$ the trajectory is concave up (down).  To determine the critical value $A_1$
 (and thus $A_c$) for which the trajectories change concavity
 we compare ({\ref{sol:lf_concavity_threshold}}) to the  trajectory itself in 
terms of the local variables $\eta$, $\xi$, and $A_1$. The local expression for the trajectory
$\xi$ is found by substituting (\ref{local:xietaAm}) into (\ref{exp:lf_multi_scale}), 
and using the expansion $\xi \sim \xi_0 + \mu \xi_1 + \mu^2 \xi_2 + \cdots$, 
\begin{align}
\xi \sim  \frac{A_1 a_m}{A_m}\eta+c_0+\mu\left[
 -(\frac{1}{2}c^2a_m-\frac{A_1^2 a_m^2}{A_m^2})\frac{\eta^3}{3}+\frac{A_1}{A_m}(1+2c_0a_m)\frac{\eta^2}{2}+c_0^2\eta+c_1\right]+\cdots,
 \label{sol:lf_local} 
\end{align}
\noindent where $c_0$ and $c_1$ are undetermined constants.
 It is the comparison of the slope of the concavity threshold (\ref{sol:lf_concavity_threshold}),
 $c^2A_m/(2A_1)$, to that of the trajectory (\ref{sol:lf_local}), $A_1a_m/A_m$, that yields the
result for $A_1$.
  If ${A_1a_m}/{A_m}<{c^2A_m}/(2A_1)$ for  $a-a_m= O(\mu)$ and
$x= O(\mu)$, then trajectories that enter the concave up region stay there,  as shown in 
Figure \ref{fig:lf_local}(LEFT), and tipping does not occur near $\eta = 0$ ($a=a_m$). 
Alternatively,  if ${A_1a_m}/{A_m}>{c^2A_m}/(2A_1)$,  then
locally the trajectory has a slope greater than that of the concavity threshold and
trajectories that are near or below
 (\ref{sol:lf_concavity_threshold}) for $a$ near $a_m$ cross into or 
continue in the concave down region, resulting in tipping as shown in Figure \ref{fig:lf_local}(RIGHT). 
Then the critical value of $A_1$ for the change in concavity satisfies
\begin{equation}
A_1^2  = \frac{\Omega c^2A_m^2}{2\mu \sqrt{c^2 A_m^2 - 1}}, \ \ \  c= \frac{\Omega}{\mu}
\qquad \Rightarrow  \qquad A_c \sim A_m + \mu A_1. 
\label{sol:lf_A_threshold}
\end{equation}
where we have used (\ref{exp:lf_threshold_2}) to eliminate $a_m$ and 
 (\ref{local:xietaAm}) to get $A_c$.

  To arrive at (\ref{sol:lf_A_threshold}), we assumed that for $A<A_c$, the trajectory 
  (\ref{sol:lf_local}) is above
 the concavity threshold 
 (\ref{sol:lf_concavity_threshold})
in the local region (\ref{local:xietaAm}). 
To verify this assumption, we compare  the $\xi$-intercept of the
 concavity threshold (\ref{sol:lf_concavity_threshold}),
 $-\frac{1}{2a_m}$,  to that of the trajectory, which to leading order is $c_0$.
Then $c_0$ is determined by matching the
local expansion (\ref{sol:lf_local}) with the outer expansion (\ref{sol:lf_outer}),
as shown in Appendix \ref{append:D},
\[ c_0 = -\frac{1}{2 a_m}+C_1+C_2. \] 
The $\xi$-intercept of the trajectory, relative to $-1/(2a_m)$ is 
obtained from the sign of $C_1 + C_2$.
 As discussed in Appendix \ref{append:D},
 for  $\frac{A_1a_m}{A_m}<\frac{c^2A_m}{2A_1}$ 
   ($\frac{A_1a_m}{A_m}>\frac{c^2A_m}{2A_1}$ ) we find that
$C_1 + C_2 >0$  ($C_1 + C_2 <0$), which is consistent with our assumptions in the derivation
of (\ref{sol:lf_A_threshold}).  Figure \ref{fig:lf_local} illustrates the local behaviour
by comparing the trajectory and the concavity threshold near $a=a_m$. 

The panels in  Figure \ref{fig:lf_amplitude_threshold} illustrate the dependence of these results
on the ratio $c=\Omega/\mu$, and compare the tipping point locations obtained by
 asymptotic approximations
 (\ref{sol:lf_tipping_position}) and (\ref{sol:lf_A_threshold}) and by 
numerical simulations. 
The variation of $A_m$ and $A_c$ with the parameter $c$
is observed by 
comparing the values of $A_m$ shown in Figure \ref{fig:lf_amplitude_threshold}(LEFT) and Figure \ref{fig:lf_threshold}
(RIGHT), and noting the values of $A_c$ where there is a jump in the location of the
tipping point in Figure \ref{fig:lf_amplitude_threshold}(RIGHT). There we see that these
jumps occur near $A=9$ and near $A=2$ for the cases with $c=1$ (blue and green lines), in contrast to the cases
with $c=.2$ where $6<A_c<7$ (red and black lines).   Figure \ref{fig:lf_amplitude_threshold}
 also shows good agreement between
results for the location of the tipping point as a function of $A$ obtained
from
 simulations and from the  asymptotic approximation for $\Omega\ll 1$.

\subsection{$\Omega = \mu^{\nu}$ ($0< \nu <1$)}
\label{sec:lf_order_nu}

\begin{figure}[h]
\begin{center}$
\begin{array}{cc}
\includegraphics[width=0.45\textwidth]{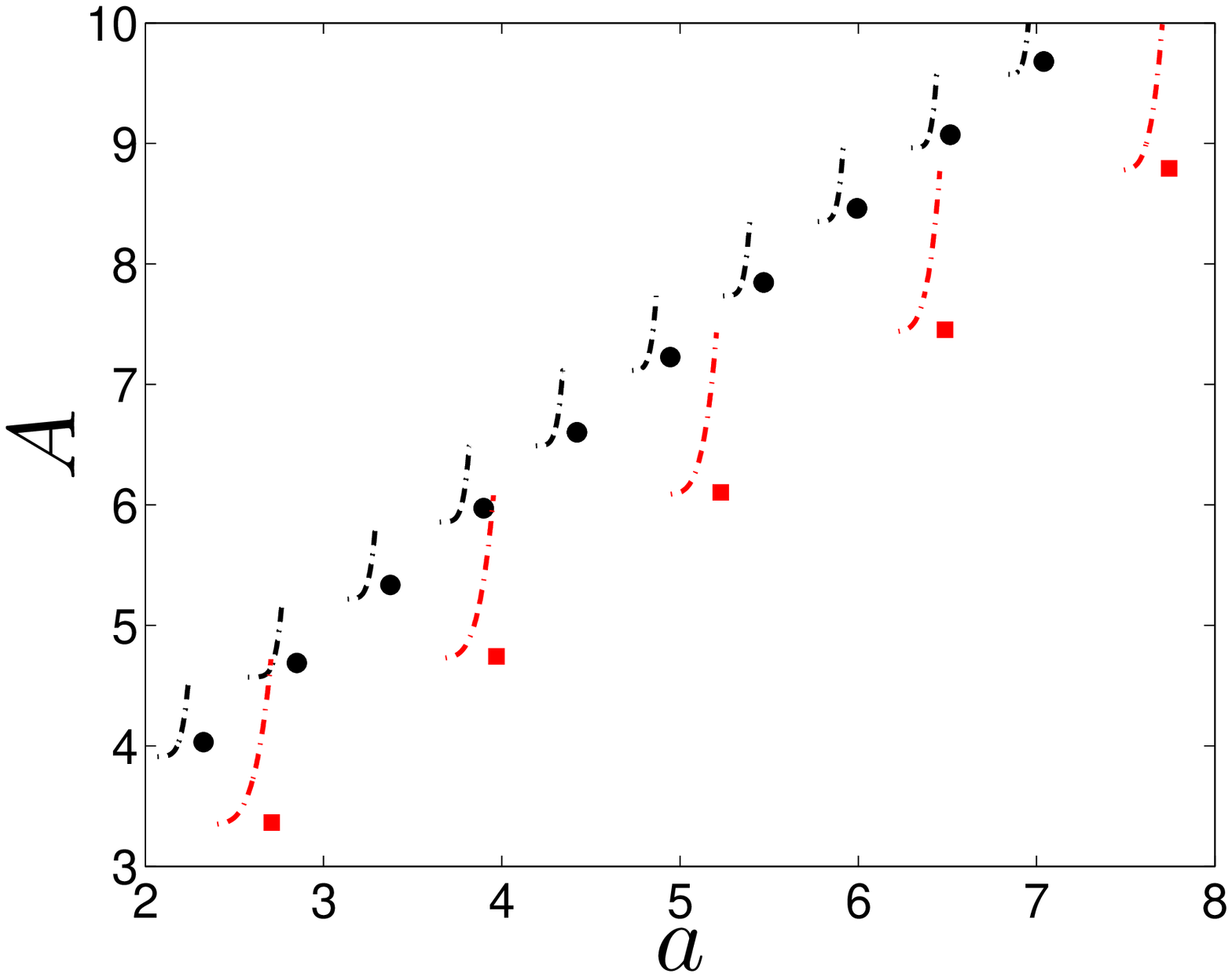} &
\includegraphics[width=0.45\textwidth]{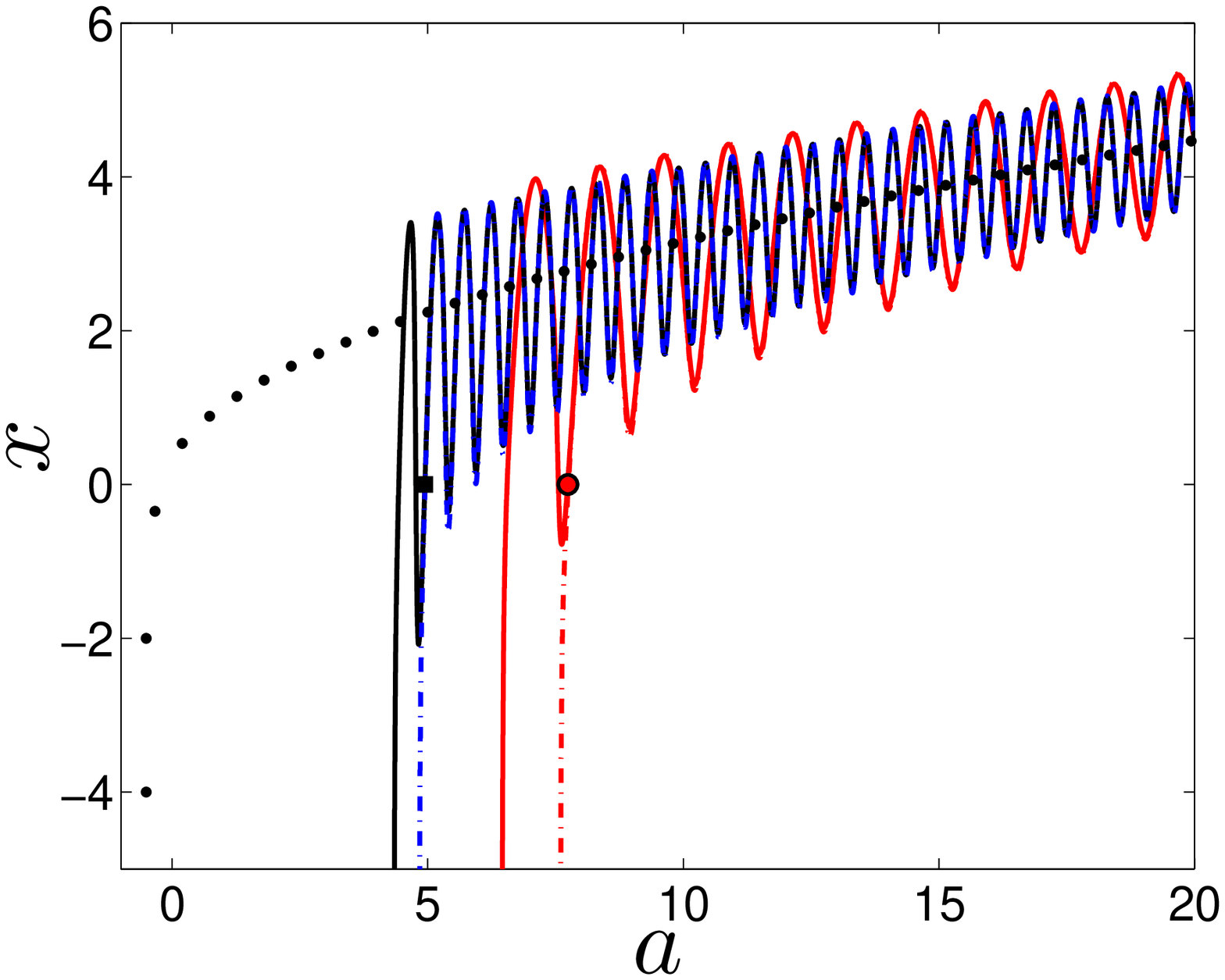}
\end{array}$
\end{center}
\caption{The drifting parameter is $\mu=.1$ for all graphs. LEFT:Numerical simulations of the position of the tipping point (dash-dotted lines),
 defined as the value of $a$ when $x=-10$, compared
with markers indicating $(a_{\rm min},A_c)$ from (\ref{sol:lf_general_concavity1})  and (\ref{sol:lf_general_concavity2}), the critical values of the amplitude at which there is an abrupt change in the location
of the tipping point. Red dash-dotted line and solid squares are for the case $\Omega=0.5$ and $a_0=20$.  Black dash-dotted lines and solid circles are for the case
$\Omega=1.2$ and $a_0 = 20$. 
RIGHT: Numerical simulation of (\ref{sys:gen_periodic_forcing}) with $A$ near $A_c$ in the location of the tipping points in (LEFT):  $A=8.7$ (red solid) and $A=8.9$ (red dash-dotted), both with $\Omega=0.5$; $A=7.1$ (Black solid) and $A=7.3$ (Blue dash-dotted), both with $\Omega=1.2$. The black square and red circle mark the related $a_{\rm min}$
for $\Omega = .5, 1.2$, respectively, satisfying (\ref{sol:lf_general_concavity1})  and (\ref{sol:lf_general_concavity2}).}
\label{fig:lf_tipping_position_gen}
\end{figure}

For $c \gg 1$, so that the low frequency is large compared with the drifting rate, 
   $\mu \ll \Omega \ll 1$,
the approximation (\ref{sol:lf_A_threshold}) is not valid since it is based on the assumption
that $c=O(1)$. Nevertheless,  
simulation results shown in Figure \ref{fig:lf_tipping_position_gen}(LEFT)
for $\Omega$ relatively small illustrate phenomena similar to that shown in Figure \ref{fig:lf_amplitude_threshold}(RIGHT). 
 As in the previous case of $\Omega = O(\mu)$, 
 for a given $a_0$ there are values $A_c$ at which there are abrupt shifts
in the location of the tipping point as a function of the amplitude. These shifts
are related to changes of concavity in the trajectory 
near a local minimum of $f(a)$ at $a=a_{\rm min}$, so that we can again
use the concavity to identify $A_c$.  As above, an early tipping occurs near
a value of $a$ where $x<0$ for $A>A_c$, while for the
same initial condition and $A<A_c$, oscillations continue for another period
with tipping taking place on a future crossing of $x=0$.  There are
also noticeable differences for $\Omega= \mu^{\nu}$, since the frequency is  larger 
than in the previous section.
For example, the consecutive crossings of $x<0$ are closer together for larger frequencies, so 
that the shifts between tipping point location at $A=A_c$ decrease with $\nu$ and the values of $A_c$ are closer together.  Furthermore,
 $f(a_{\rm min})<0$ 
may be $O(1)$ for $a_{\rm min}$ near the tipping point, in contrast to $f(a_m)=0$ in 
Section \ref{sec:lf_order1} 
as shown in Figure \ref{fig:lf_threshold}. As a result, we modify our approach for calculating
$A_c$.

To approximate $A_c$ for $\Omega = \mu^{\nu}$, we first consider the outer approximation to 
(\ref{exp:lf_multi_scale}) given by (\ref{sol:lf_outer}) with $c=\mu^{\nu-1}$.
As in Section \ref{sec:lf_order1} this approximation is
valid only for $f(a)>0$, since as $f(a)$ approaches $0$, 
 the terms shown in (\ref{sol:lf_outer}) are both of the same order for $\mu \ll 1$.
To estimate $A_c$
in this case, we must consider multiple local minima $a_{\rm min}$ of $f(a)$
  for a given $a_0$, satisfying
\begin{eqnarray}
f(a_{\rm min})&=&a_{\rm min}+A \mbox{sin}(\mu^{\nu-1}  (a_0-a_{\rm min}))<0, \label{exp:lf_general_threshold_1}  \\
f'(a_{\rm min})&=&1-\mu^{\nu-1} A \mbox{cos}(\mu^{\nu-1} (a_0-a_{\rm min})) =0, \,\,\, \Rightarrow \,\,\,\mbox{cos}(\mu^{\nu-1} (a_0-a_{\rm min})) =\mu^{1-\nu}\frac{1}{A}. \label{exp:lf_general_threshold_2} 
\end{eqnarray}
Since we are considering larger values of the frequency $\Omega$, the
values of $a_{\rm min}$ are relatively close together and close to values where $f(a)=0$. Then we can not use a local analysis
near values of $a$ where $f(a)=0$,
as was done near $a_r$ in (\ref{exp:lf_local_IO_A}), but instead restrict our analysis near $a_{\rm min}$ for $f'(a_{\rm min})=0$.

Substituting $a=a_{\rm min}+\mu \eta$ and $x = \mu \xi$ into 
$x''(a)$ using (\ref{exp:lf_multi_scale}), 
   (\ref{exp:lf_general_threshold_1}), and (\ref{exp:lf_general_threshold_2}), we get the local expressions for the trajectory and 
its concavity. Keeping the terms
that dominate these expressions for $\mu\ll 1$ yields
\begin{align}
\mu \frac{d\xi}{d\eta} &= -f(a_{\rm min}), \label{exp:lf_local_general_IO}\\
\mu\frac{d^2\xi}{d{\eta}^2} &= \mu^{2\nu}(f(a_{\rm min})-a_{\rm min})\eta-2\mu f(a_{\rm min})\xi. \label{exp:lf_concavity_general_IO} 
\end{align}
As in Section \ref{sec:lf_order1}, we look for the critical value of $A=A_c$ that corresponds to the value where 
the local slope of the trajectory is equal to that of the 
concavity threshold. The latter is obtained from setting $\xi''(\eta)=0$ in 
(\ref{exp:lf_concavity_general_IO})
 yielding
\begin{eqnarray}
\xi_c = \frac{\mu^{2\nu}(f(a_{\rm min})-a_{\rm min})\eta}{2\mu f(a_{\rm min})}. \label{exp:lf_concavity_thresh_gen}
\end{eqnarray} 
Equating the coefficient of $\eta$ from the righthand side of (\ref{exp:lf_concavity_thresh_gen})
with $\xi'(\eta)$ in (\ref{exp:lf_local_general_IO}) gives the relation
\begin{eqnarray}
-2f^2(a_{\rm min})=\mu^{2\nu}(f(a_{\rm min})-a_{\rm min}). \label{exp:lf_equalslope}
\end{eqnarray}
Then, using 
(\ref{exp:lf_general_threshold_1})-(\ref{exp:lf_general_threshold_2}) 
and (\ref{exp:lf_equalslope}), we can eliminate $f(a_{\rm min})$ and obtain a system of equations for
the pair $(a_{\rm min},A_c)$,
\begin{eqnarray}
a_{\rm min} &=& \sqrt{A_c^2-\frac{\mu^2}{\Omega^2}}-\frac{\Omega}{\sqrt{2}} \left[A_c^2-\frac{\mu^2}{\Omega^2}\right]^{1/4}
, \label{sol:lf_general_concavity1} \\
a_{\rm min} &=& a_0 +
 \frac{\mu}{\Omega}\left[\cos^{-1}\left(\frac{\mu}{A_c \Omega}\right) -2k\pi\right], \label{sol:lf_general_concavity2}
\end{eqnarray}
where $k$ is a positive integer. Note that for fixed $a_0$ the system of equations (\ref{sol:lf_general_concavity1})-(\ref{sol:lf_general_concavity2}) is valid for $\mu=o(\Omega)$. It has more than one solution, that is, more than one pair $(a_{\rm min}, A_c)$
at which is there is a change in concavity and thus a jump in the location of the tipping point as a function of $A$.  
Then for small
$\mu$ or sufficiently large $A_c$, $\cos^{-1}(\mu/(A_c\Omega))\sim \pi/2$ in
(\ref{sol:lf_general_concavity2}), so that the jumps in the value of
$a_{\rm min}$ are dominated by $2k \pi {\mu}/{\Omega}$ for integers $k$, related
to the multi-valued function $\cos^{-1}(\cdot)$.

Figure \ref{fig:lf_tipping_position_gen}(LEFT) compares the position of the tipping
points obtained from numerical simulation with the
asymptotic approximation for the critical values $a_{\rm min}$ and $A_c$ given by 
(\ref{sol:lf_general_concavity1})-(\ref{sol:lf_general_concavity2}).
 Figure \ref{fig:lf_tipping_position_gen} illustrates this
approximation for $\mu= .1$, both for $\Omega = \mu^{\nu}$ for $\nu\approx .3$, and for $\Omega = O(1)$.
While the latter case does not appear to be covered by the asymptotic analysis, we recall that 
in Section \ref{sec:hf_large_A}, we showed that the case $A>\Omega^2$ with $\Omega = \mu^{-\lambda}$ for 
$\lambda>0$ can be studied as a low frequency case by rescaling (\ref{sys:gen_periodic_forcing}) 
into the normalized system  (\ref{sys:norm_periodic_forcing}). This rescaled system, with frequency $\omega<1$
and slow drift parameter $M\ll 1$, can be analyzed using expressions
analogous to (\ref{sol:lf_general_concavity1})-(\ref{sol:lf_general_concavity2}). 
For $\frac{\mu}{\Omega} \ll 1$, we see that (\ref{sol:lf_general_concavity1})-(\ref{sol:lf_general_concavity2})
provide a good approximation for the location of the tipping point and the values of $A$ for which there
is a jump in location, even though we
do not have an analytical approximation for the solution $x$ for $\Omega = O(1)$ and $\mu\ll 1$.
For larger values of $\Omega$ with $\mu$ fixed, there are more values of $A_c$ where this
jump occurs, and the corresponding values of $a=a_{\rm min}$ are closer together.
 The increasing number and reduced distance between these jumps indicates that as 
$\Omega$ increases further, the relationship between
the amplitude $A$ and the location of the tipping point is a continuous function,
  as shown in Section \ref{sec:hf_slow_drifting} for large frequency.


\begin{figure}[ht]
\centering 
\scalebox{0.5}
{\includegraphics{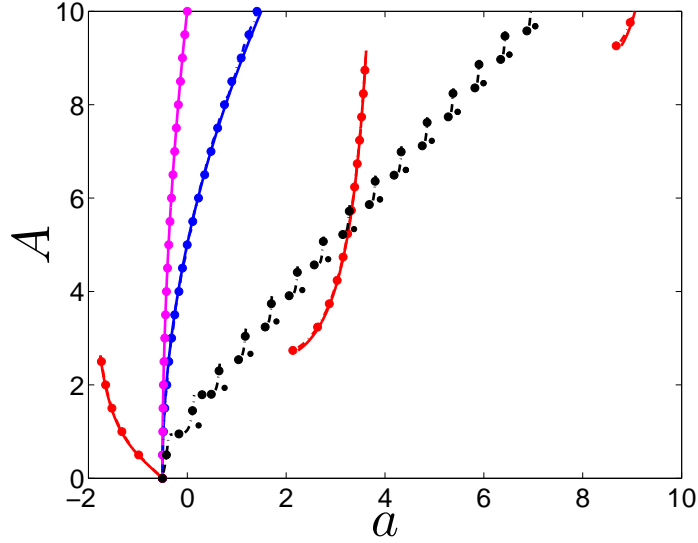}}
\caption{Position of the tipping point as a function of $A$ in (\ref{sys:gen_periodic_forcing}) obtained by numerical simulation, defined as the value of $a$ where $x=-10$ ( dash-dotted
line with solid circles) compared with the asymptotic approximation (solid line) for 
$\Omega=0.1$(red);  $\Omega=5$(blue); and $\Omega=10$(magenta) with $\mu=0.1$ 
and $a_0=20$. The black dash-dotted lines give the location of the tipping point
determined numerically for $\Omega = 1.2$, which are compared with the black circles indicating the value of $A$ where there is jump in the
location of the tipping point based on the asymptotic approximation (\ref{sol:lf_general_concavity1})  and (\ref{sol:lf_general_concavity2}).}
\label{fig:summary}
\end{figure}

The analysis of large amplitude oscillatory forcing in Section \ref{sec:hf_large_A} suggests a 
rescaling of the system,
 which allows the asymptotic approximations to be extended
to some cases where the frequency is not particularly large or small relative to the drifting rate of the bifurcation parameter.
  In Figure \ref{fig:summary} we compare results for low, high, and O(1) frequencies $\Omega$, 
illustrating how 
the location of the tipping point varies with the amplitude and frequency of the periodic forcing.

\section{Example 1: Morris-Lecar Model}
\label{sec:model_ML}

We consider the Morris-Lecar(ML) model \cite{Morris1981,Yu2013} with a slowly increasing
input current and an additive periodic external input. 
The ML model has been used widely as a two-dimensional model, relatively straightforward
to analyze, that captures certain
dynamical behaviour observed for different types of neurons with a variety of states. 
It has the general form
\begin{align}
\gamma \frac{\mbox{d}v}{\mbox{d}\hat{t}} &= -g_{Ca}\hat{m}_{\infty}(v)(v-v_{Ca})-g_{K}(v-v_K) w-g_L(v-v_L)+I_{\rm bias}+ I_{\rm ext}, \label{sys:ML} \\
\frac{\mbox{d}w}{\mbox{d}\hat{t}} &= \hat\kappa(v)(\hat{w}_{\infty}(v)-w),\nonumber \\
\frac{\mbox{d}I_{\rm bias}}{\mbox{d}\hat{t}} &= \hat{\mu}, \nonumber\\
 \hat{m}_{\infty}(v) &=\frac{1}{2}\left(1+\tanh\left(\frac{v-v_1}{v_2}\right)\right), \qquad \hat{w}_{\infty}(v)=\frac{1}{2}\left(1+\tanh(\left(\frac{v-v_3}{v_4}\right)\right), \qquad \hat{\kappa}(v) = \phi \mbox{cosh}\left(\frac{v-v_3}{2v_4}\right), \nonumber
\end{align}
where $v$ represents the membrane voltage potential and $w$ is the probability of opening of $\mbox{K}^{+}$ channels ($0 \leq w \leq 1$). This simplied model of excitability depends on three ionic currents for 
calcium, potassium, and leakage,
described by the terms with the relevant conductances, $g_{Ca}$, $g_K$, and $g_L$.
 The periodic external input is $I_{\rm ext}$, 
$I_{\rm bias}$ represents the bias or base current in the neuron, and we refer the reader to \cite{MLweb} for additional details.  Here we take the commonly
used fixed constants as given in \cite{Yu2013}: 
 $v_K=-84$mV, $v_L=-60$mV, $v_{Ca}=120$mV, $c=20 \mu \mbox{F}/\mbox{cm}^2$, $v_1=-1.2$mV and $v_2=18$mV,
  $g_K=8 \mbox{mS}/\mbox{cm}^2$, $g_L=2\mbox{mS}/\mbox{cm}^2$, $g_{Ca}=4.4\mbox{mS}/\mbox{cm}^2$, $v_3=12 \mbox{mV}$, $v_4=17.4\mbox{mV}$, and $\phi=\frac{1}{15}\mbox{ms}^{-1}$. 
The dynamical behavior of (\ref{sys:ML}) with $\hat\mu=0$
 has been well studied for $I_{\rm ext}=0$, with the current $I_{\rm bias}={\rm constant}$ often playing
 the role of the
bifurcation parameter. 
For that case there is a SNIC (saddle node on an invariant circle) bifurcation \cite{wiki_ML}
at $I_{\rm bias}=I_c \approx 44.09 \mu \mbox{A}/\mbox{cm}^2$,  and corresponding 
voltage $ v=v_{c}\approx -27.14$mV,
so that for $I_{\rm bias}<I_{c}$, 
there exists a stable equilibrium and unstable equilibria. For $I_{\rm bias}>I_{c}$, the system exhibits
oscillations whose frequency  increases from zero with $I_{\rm bias}-I_c$, and there is a Hopf
bifurcation for $I_{\rm bias}$ above $60 \mu \mbox{A}/\mbox{cm}^2$, which 
we do not consider here. The dynamics of
a ML-type neuron with periodic external input and constant
 $I_{\rm bias}$ has been studied in \cite{Lee2007}-\cite{Lee2008}.

\begin{figure}[h]
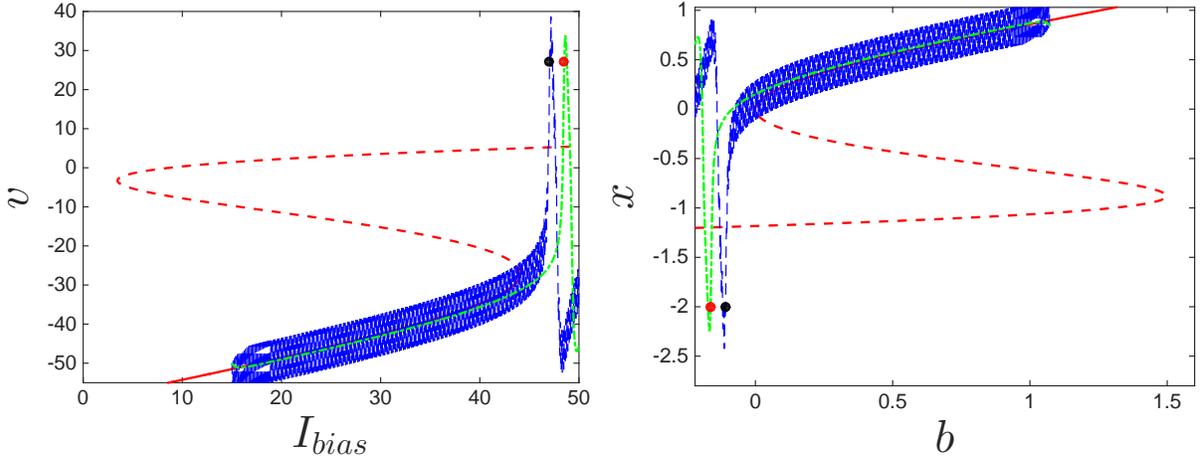

\begin{center}
\includegraphics[width=0.45\textwidth]{figure_ML_bif1.eps} 
\includegraphics[width=0.45\textwidth]{figure_ML_bif2.eps}
\end{center}
\caption[ML Bifurcation]{\label{fig:ML_bifurcation} LEFT: The bifurcation diagrams of 
(\ref{sys:ML}) indicating the stable (red solid) and unstable (red dashed) branches. 
The trajectories for $\hat{A}=0$ (green dash-dotted line) and $\hat{A}\neq 0$ (blue dashed line)
are superimposed on these curves , with the corresponding tipping points indicated by solid circles.
Note that tipping occurs
for a smaller value of $I_{\rm bias}$ for the case $\hat{A}\neq 0$. RIGHT: Same as LEFT, but
for the transformed system (\ref{sys:norm_ML}).}
\end{figure}

We focus on transitions from the stable steady state to oscillations
for a current that is slowly increasing with rate $\hat{\mu}\ll 1$, modeled here as 
slowly varying $I_{\rm bias}$ for simplicity but alternatively it could be an input current,
 and a periodic forcing $I_{\rm ext} = \hat{A}\sin(\hat{\Omega} t)$. 
Sample trajectories superimposed on
the static bifurcation structure are shown in Figure \ref{fig:ML_bifurcation}.
 To determine the tipping points, that is, the value of $I_{\rm bias}$ for the transition, 
we transform (\ref{sys:ML}) to a form similar to the canonical model 
(\ref{sys:gen_periodic_forcing}),
rescaling time 
$\hat{t} = \gamma t$ and normalizing the system about the bifurcation point ($v_{c}$, $I_{c}$), 
\begin{align}
\frac{\mbox{d}x}{\mbox{d}t} &= h(x) -g_K(x+{\cal D})w +b(t) +A\mbox{sin}(\Omega t) \nonumber \\
&= b(t) -{g_{Ca}}m_{\infty}(x)(x+1-\frac{v_{Ca}}{v_c})-{g_L}(x+1-\frac{v_L}{v_c})+\frac{I_{c}}{v_c}-g_{K}(x+ 1-\frac{v_K}{v_c})w+A\mbox{sin}(\Omega t), \label{sys:norm_ML} \\
\frac{\mbox{d}w}{\mbox{d}t} &= \kappa(x)(w_{\infty}(x)-w), \qquad
\kappa(x) = \gamma\hat{\kappa}(v_c(x+1)), \qquad 
w_{\infty}(x) = \hat{w}_{\infty}(v_c(x+1))
\qquad \nonumber\\
\frac{\mbox{d}b}{\mbox{d}t} &= -\mu = -\frac{\gamma}{|v_c|}\hat\mu \nonumber\\
 x & = \frac{v-v_c}{v_c}, \qquad  b=\frac{I_{\rm bias}-I_{c}}{v_c}, \qquad
  \Omega = \gamma\hat\Omega_c, \qquad  v_c A = \hat{A}, \qquad 
m_{\infty}(x)= \hat{m}_{\infty}(v_c(x+1)) \, .
\nonumber
\end{align}
For $\mu=0$ and $A=0$, the bifurcation point ($b_c,x_c$) is at the origin for the transformed
system (\ref{sys:norm_ML}).
Figure \ref{fig:ML_bifurcation} shows the bifurcation diagram for (\ref{sys:ML}) and
 (\ref{sys:norm_ML}), with trajectories for
$A=0$ and $A\neq 0$ obtained from 
 the corresponding
 systems with $\mu\neq 0$. We note here that for the transformed system, the bifurcation
parameter $b$ and dependent variable $x$ take $O(1)$ values. For the system in this form, it 
is then straightforward to identify different scales related to slow drift, high and
low frequency, or large amplitude. Such an identification is important in choosing the correct
approach to approximate the tipping point.



Before considering the case with $A \neq 0$, we give an
analytical expression for the tipping point without oscillatory forcing $A =0$, which is
compared below with the tipping point for $A\neq 0$.  It is useful to find the tipping point for the generic system of the
form
\begin{eqnarray}
\frac{dx}{dt} = Da + k_0 + k_1 x + k_2 x^2 \, ,  \label{hpoly0}
\quad \frac{da}{dt} = -\mu \, .
\end{eqnarray}
Applying the scaling of \cite{haberman1979} described in Section \ref{sec:model},
completing the square in (\ref{hpoly0}) and
using a straightforward transformation again yields an Airy equation, from
which we obtain the singularity corresponding to the tipping point for (\ref{hpoly0}),
\begin{eqnarray}
a_g = (D|k_2|)^{-1/3}a_d- \frac{a_s}{D} \qquad \mbox{for} \qquad a_s = k_0 + \frac{k_1^2}{4|k_2|} , \label{a_g}
\end{eqnarray}
where $a_d$ is the tipping point in (\ref{sol:location_tipping}) for the symmetric case  
and $a_s$ gives the 
shift for a general quadratic polynomial in (\ref{hpoly0}), in contrast
to the symmetric case.

For $A=0$ in (\ref{sys:norm_ML}) we obtain an approximation to
 the tipping point using a local expansion about $x=0$.
 Noting that $0<\kappa(0)=O(1)$ we
expect $w$ to relax to $w_\infty$, so that the local approximation for (\ref{sys:norm_ML}) is
\begin{eqnarray}
 & & \frac{dx}{dt} \approx b + h_0^0 + h_1^0 x + h_2^0 x^2  - g_K(x+{\cal D})[w_{\infty}(0) + w_{\infty}'(0)x + w_{\infty}''(0)x^2/2]
  \label{xpoly0}\\
& & \frac{db}{dt} = -\mu, \qquad h_j^{0} = j! h^{(j)}(x)|_{x=0}  \, . \nonumber
\end{eqnarray}
Then the approximate tipping point for $A=0$ is given by (\ref{a_g}) with $D=1$ and
\begin{eqnarray}
k_0 = h_0^0 -g_K{\cal D}w_{\infty}(0), \qquad
k_1 = h_1^0 -g_K(w_{\infty}(0) +{\cal D}w_{\infty}'(0)), \qquad
k_2 = h_2^0-g_K\left(w_{\infty}'(0) +{\cal D}\frac{w_{\infty}''(0)}{2}\right) 
\, .\qquad \label{b_tip_A0_hf}
\end{eqnarray}
Note that $k_0=k_1=0$ in (\ref{b_tip_A0_hf}) following from the quadratic form
for the saddle node bifurcation at $b=0$ and $x=0$, so that resulting tipping point for
(\ref{xpoly0}) is
$b_{ml} = |k_2|^{-1/3}a_d$.  

\subsection{High Frequency Oscillation $\Omega \gg 1$}
\label{sec:ML_hf}

We use the approach of Section \ref{sec:hf_slow_drifting} to study the case $\mu\ll 1$ and $\Omega\gg 1$ in  (\ref{sys:norm_ML}), again
writing $\Omega=\mu^{-\lambda}$. We restrict our attention to the case where $A=O(1)$.

We introduce the fast time scale $T=\mu^{\lambda} t$ and the slow time scale $\tau=\mu t$, 
and substitute the multiple scale approximations for $x$ and $w$
\begin{align}
x &\sim x_0(T,\tau)+\mu^{\lambda} x_1(T,\tau)+ \cdots, \\
w &\sim w_0(T,\tau)+\mu^{\lambda} w_1(T,\tau) +\cdots, \label{xwexp}
\end{align}
into (\ref{sys:norm_ML}), yielding the sequence of equations
\begin{align}
O(1): \ &x_{0T} = w_{0T}=0, \Rightarrow x_0=x_0(\tau), w_0 = w_0(\tau), \label{exp:ML_outer_LO1} \\
O(\mu^{\lambda}):\  &x_{1T} = b+ h(x_0)-g_k(x_0+D)w_0+A\mbox{sin}(T), \nonumber\\ 
         &w_{1T} = \kappa(x_0)(w_{\infty}(x_0)-w_0) \, .  \label{exp:ML_outer_SO2} 
\end{align}
We apply a solvability condition  similar to (\ref{eq:solv_condition}), namely that the
inner product of the right-hand side of  (\ref{exp:ML_outer_SO2})
with the homogeneous solutions ($x$ and $w$ constant with respect to $T$) 
vanishes. This yields
\begin{align}
 &\lim_{L \to \infty} \frac{1}{L}\int_0^{L} \left[ b + h(x_0)-g_K(x_0+{\cal D})w_0 + A\sin(T)
\right] \, dT= 0 \label{ML_solv}\\ 
 &\lim_{L \to \infty} \frac{1}{L}\int_0^{L} \kappa(x_0)(w_{\infty}(x_0)-w_0) \, dT=0 \, . \nonumber
\end{align}
Then the asymptotic approximation of $x$ for $b(\tau)=O(1)$ is 
\begin{align}
x &\sim x_0+\mu^{\lambda}[-A \mbox{cos}(T)]+\cdots, \ \mbox{ where } \
   -b = h(x_0) - g_K(x_0+ {\cal D})w_{\infty}(x_0) \label{sol:ML_outer_x} \\
w &\sim w_{\infty}(x_0)+\cdots \, ,  \nonumber
\end{align}
 where $x_0$ in (\ref{sol:ML_outer_x}) 
 is the stable equilibrium of (\ref{sys:norm_ML}), 
shown as the solid curve in Figure \ref{fig:ML_bifurcation} (RIGHT). 
The asymptotic approximation describes the attracting solution of (\ref{sys:norm_ML}) away from $b=0$, but does not
give the  tipping point.

 To determine the location of the tipping point, we use
a local approximation for $x\ll 1$ similar to that used in  
 Section \ref{sec:hf_slow_drifting},
\begin{eqnarray}
 x&=&{\cal X}(T)+z(T,s), \qquad {\cal X}(T)=-\Omega^{-1}A\mbox{cos}(T)\label{ML_inner}\\
 w&=& W_0(T,s) + \Omega^{-1} W_1(T,s) \nonumber\\
 z&=& z_0(T,s) + \Omega^{-1} z_1(T,s) , \qquad s = \mu^{1/3} t \, ,
\label{xhf_inner}
\end{eqnarray}
 that is, the solution is given by oscillations on the $T$ scale plus a correction $z$. We could
introduce additional scalings of $z= O(\mu^{1/3})$ and $b = O(\mu^{2/3})$ 
as in (\ref{sec:hf_slow_drifting}), 
but as it does not affect the result, we use $z$ and $b$ for simplicity.
 We also introduce quadratic polynomials given by Taylor expansions of $h(x)$, $\kappa(x)$
and $w_{\infty}(x)$ about $x=0$, 
facilitating
explicit expressions for the tipping point. Then substituting (\ref{ML_inner})-(\ref{xhf_inner}) into (\ref{sys:norm_ML}), 
we find that ${W_0}_T$ and ${z_0}_T$ vanish, so that $W_0 = W_0(s)$, $z_0 = z_0(s)$.
Approximations for $W_0$ and $z_0$ are obtained at the next order, 
 analogous to (\ref{exp:hf_drifting_in_order1}) as shown in Appendix \ref{append:E}.  We apply the solvability condition as in (\ref{ML_solv}), to the equations for $W_1$ and $z_1$ yielding
 for $W_0$ 
\begin{eqnarray}
 W_0 &\sim&  w_\infty(0) +  W_{00} + W_{01}z_0 + W_{02}z_0^2 \, .\label{w0approxmain}
\end{eqnarray}
with the coefficients $W_{00}$,  $W_{01}$, and $W_{02}$ given in (\ref{w0approx}).
Then the equation for $z_0$ written in terms of $t$ is,
\begin{align}
{z_0}_t&=b+h_0^0+h_1^0 z_0+h_2^0 z_0^2+ h_2^0 \frac{A^2}{2\Omega^2} -g_K(z_0 +{\cal D})W_0 
, \label{exp:ML_inner_average} 
\end{align}
Applying (\ref{hpoly0})-(\ref{a_g}) yields the singularity for (\ref{exp:ML_inner_average})
\begin{eqnarray}
b_{hf} &=& |h_2^{\rm hf}|^{-1/3}a_d- b_s \qquad \mbox{for} \qquad b_s = h_0^{\rm hf}
 + \frac{(h_1^{\rm hf})^2}{4|h_2^{\rm hf}|} \qquad \mbox{where}  \label{sol:ML_hf_tipping}\\
  h_0^{\rm hf} & =&-g_k{\cal D}W_{00}+ h_2^0 \frac{A^2}{2\Omega^2} \, , \qquad
  h_1^{\rm hf} = -g_K(W_{00} + {\cal D}(W_{01}-w_\infty'(0))) \, . \nonumber\\
  h_2^{\rm hf} & =&  h_2^0 - g_K(W_{01} + {\cal D}W_{02}) \, , \nonumber
\end{eqnarray}
where we have used the fact that $k_0$ and $k_1$ in (\ref{b_tip_A0_hf}) vanish.
Note that for $\Omega\gg 1$, $W_{01} - w_\infty'(0)$ and $W_{00}$ are proportional
to $\Omega^{-2}$, as are $h_0^{\rm hf}$ and $h_1^{\rm hf}$. Then we recover
(\ref{b_tip_A0_hf}) as $\Omega^{-2}\to 0$.

\begin{figure}[h]
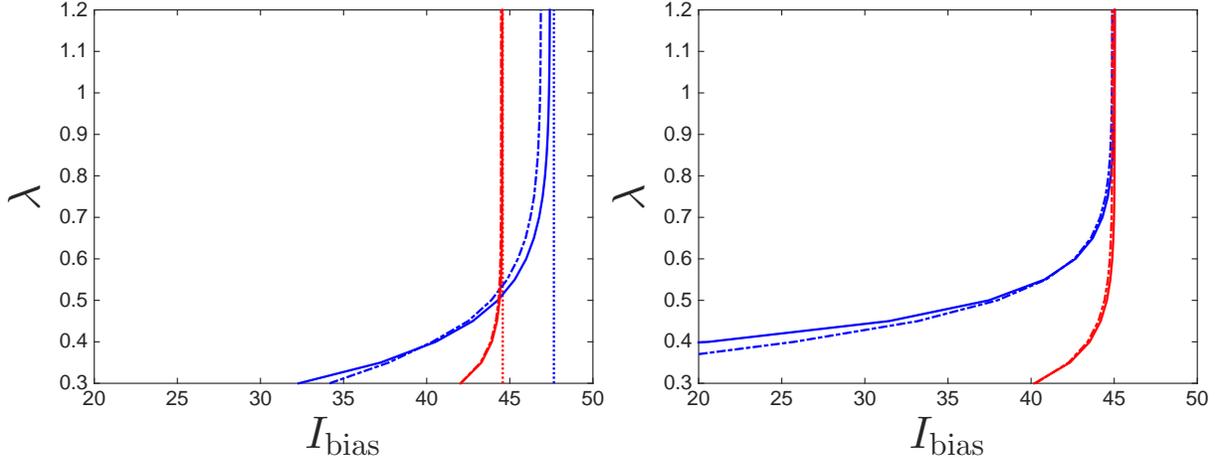

\begin{center}
\includegraphics[width=0.45\textwidth]{figure_ML_hf_tipping1.eps} 
\includegraphics[width=0.45\textwidth]{figure_ML_hf_tipping2.eps}
\end{center}
\caption[ML High Frequency Case]{\label{fig:ML_hf_tipping}  The value $I_{\rm bias}$ 
corresponding to tipping in (\ref{sys:ML}) vs. $\lambda$ for $\Omega=\mu^{-\lambda}$.
The asymptotic results obtained from  (\ref{sol:ML_hf_tipping}) (solid lines) are
compared with numerical simulations (dash-dotted lines), obtained using the
 condition for tipping as
 $v=|v_c|$ in (\ref{sys:ML}).   
LEFT: Results are
shown for amplitude
$A=2$ ($\hat{A} =-54.275$) and
 two drifting rates $\mu = .001$ ($\hat{\mu} = 0.0014$) (red) and $\mu = .02$ 
($\hat{\mu}=0.0271$) (blue)
with asymptotes for large $\Omega$
 approximated by (\ref{b_tip_A0_hf}) (vertical dotted lines)
as $I_{\rm bias} = 44.58$ and $I_{\rm bias}= 47.65$, respectively.
 RIGHT:  Results are shown for $\mu = .003$ ($\hat\mu = .0041$) and amplitude values
$A = 2$ (red) and $A=8$ ($\hat{A} = 217.10$) (blue), with earlier tipping for larger ratios $A/\Omega$.}
\end{figure}

In Figure \ref{fig:ML_hf_tipping}  we compare
the asymptotic approximation for the tipping point for $I_{\rm bias}$  obtained
from (\ref{sol:ML_hf_tipping}) as a function of the
  exponent $\lambda$ for $\Omega=\mu^{-\lambda}$ 
 with numerical simulations 
for different drifting rates $\mu$ and amplitudes $A$.
 For larger values of the ratio $A/\Omega$ we see an increased
advance of the tipping point.  
 For larger $\lambda$ corresponding to higher frequencies $\Omega$, the delay effect caused by slowly varying drifting rate $\hat{\mu}$ dominates the location of the tipping point, regardless of the amplitude of the oscillatory forcing. For smaller $\lambda$, larger amplitude $\hat{A}$ triggers earlier tipping. 

We note the importance of deriving these results within the framework of the transformed system
(\ref{sys:norm_ML})
for $x$.   Figure \ref{fig:ML_hf_tipping} shows good agreement of the asymptotic and computational results for $\Omega = \mu^{\lambda}$ over the range $\lambda >.3$, as predicted by the asymptotic results for
the canonical model.  Considering the variables of the original ML system (\ref{sys:ML}),  we
see that $\hat\Omega$ is smaller than $\Omega$ by a factor of $\gamma = 20$, 
while $\hat\mu \approx 1.36 \mu$. Applying the asymptotic analysis in the original variables  
could suggest that good approximations are available over a smaller range of parameter values, and
could lead to the use of an expansion that is not appropriate for the dynamics.


\subsection{Low Frequency Oscillation $\Omega \ll 1$}
\label{sec:ML_lf}

As in Section \ref{sec:lf}, we write $\Omega = {\cal C}\mu$, rescale time as $\tau=\mu t$ and write $x$ and $w$ as functions of $b(\tau) = b_0-\tau$. Then (\ref{sys:norm_ML}) becomes
\begin{align}
-\mu\frac{\mbox{d}x}{\mbox{d}b} &= b+h(x)
-g_K(x+ {\cal D})w+
A\mbox{sin}({\cal C}(b_0-b))\label{exp:ML_lf_rescale1} \\
-\mu\frac{\mbox{d}w}{\mbox{d}b} &= \kappa(x)(w_{\infty}(x)-w) \, .\label{exp:ML_lf_rescale2} 
\end{align}
Substituting the asymptotic expansion (\ref{xwexp}), we find the leading order approximation 
of $w(x)$ is $w_0 = w_{\infty}(x_0)$ for $x_0 \sim O(1)$ and for $x_0$
\begin{align}
& b+h(x_0)  -g_K(x_0+ {\cal D})w_{\infty}(x_0)+ A\mbox{sin}({\cal C}(b_0-b)) \equiv {\Phi}(x_0,b) = 0 \, .
 \label{f_bdef}
\end{align}
As in Section \ref{sec:lf},
 we expect the tipping to occur in the region where the outer approximation is not valid,
for a pair $(b_r,x_r)$. A Taylor expansion of ${\Phi}$ about $x_r$ evaluated at $b_r$, keeping up to
quadratic terms, yields 
\begin{eqnarray}
 x_0-x_r & =& -\frac{1}{2{\Phi}_2^r}\left[ {\Phi}_1^r\pm \sqrt{({\Phi}_1^r)^2-4{\Phi}_2^r\cdot {\cal F}(b_r)} \right] 
\qquad {\Phi}_j^r = j!\frac{\partial^j \Phi}{\partial x^j}(x_r,b_r),\qquad
 {\cal F}(b_r) ={\Phi}_0^r \, ,
\label{sol:ML_lf_LO}
\end{eqnarray}
where the positive sign corresponds to the solution that is attracting.  It follows that the approximation
for $x_0$ is no longer 
valid for $b<b_r$ where $4{\cal F}(b_r) = (\Phi_1^r)^2/\Phi_2^r$ and 
$2(x_0-x_r) = -\Phi_1^r/\Phi_2^r$.
 As  $x_0 \to x_r$, ${\Phi}_1^r \to 0$, from which we conclude that
$x_r=x_c =0$, the value of $x$ at the saddle node bifurcation for the static case. 
As a result we expect that tipping occurs near $b=b_r$ and $x=0$.
Note that ${\cal F}(b)$ then
has the same functional form as $f(a)$ in (\ref{exp:lf_multi_scale}). 

  As in Section \ref{sec:lf}
 we use a local analysis in terms of $b-b_r=\mu^{2/3} {\cal B}$ and $x=\mu^{1/3}\chi$ 
applied to (\ref{exp:ML_lf_rescale1}), yielding
\begin{equation}
-\mu^{2/3}\chi_{\cal B} = {\cal F}'(b_r) \mu^{2/3}{\cal B}+k_2(\mu^{1/3}\chi)^2
+ o(\mu^{2/3}). \label{exp:ML_lf_local1}
\end{equation}
where ${\cal F}'(b_r)=1-{\cal C}A\cos({\cal C}(b_0-b_r))$ and $k_2$ is given in
(\ref{b_tip_A0_hf}).
Then using the result (\ref{a_g}) yields the tipping point
\begin{align}
b_{\rm lf}&=b_r+\frac{a_d}{({\cal F}'(b_r)|k_2|)^{1/3}} \, .  
\label{sol:ML_lf_tipping}
\end{align}
The result for $b_{\rm lf}$ relies on ${\cal F}'(b_r)=O(1)$, and as in Section \ref{sec:lf}, if ${\cal F}'(b)$ vanishes the approximation
(\ref{exp:ML_lf_local1})
is not valid. Near a value $b=b_m$ where ${\cal F}(b_m)={\cal F}'(b_m)=0$ we again expect to find a jump in the location of the tipping point 
as a function of amplitude for a critical value $A=A_c$.
Following the analysis of Section \ref{sec:lf} we provide a local analysis 
for $(b,A)$ close to $(b_m,A_m)$, which satisfy the conditions
\begin{align}
{\cal F}(b_m) &=b_m+A\mbox{sin}({\cal C}(b_0-b_m))=0, \label{exp:ML_lf_bm1} \\
{\cal F}'(b_m) &= 1-{\cal C} A_m \mbox{cos}({\cal C}(b_0-b_m)) = 0, \; \Rightarrow 1+{\cal C}^2b_m^2
={\cal C}^2 A_m^2\, . \label{exp:ML_lf_bm2}
\end{align}
To find $A_c>A_m$ for which there is a jump in the location of the tipping point, we 
consider the concavity of $x$ in (\ref{exp:ML_lf_rescale1}) 
with respect to $b$ near $b_m$.
Substituting local variables 
\begin{equation}
b-b_m = \mu \eta, \; \qquad x  = \mu \xi , \;  \qquad \mbox{ and } A = A_m+\mu A_1, \label{exp:ML_lf_localvar} 
\end{equation}
into (\ref{exp:ML_lf_rescale1}) and (\ref{exp:ML_lf_rescale2}), we find the leading
 contributions to $x'(b)$ and $x''(b)$ in terms of the local variables,
\begin{align}
 \frac{\mbox{d}^2 \xi}{\mbox{d} \eta^2} &= 
\mu\frac{A_1}{{\cal C} A_m}-\mu({\cal C}^2 \eta b_m)- \mu k_2 \xi\frac{2A_1b_m}{A_m}
 + O(\mu^2) \label{exp:ML_lf_concavity} \\
 \frac{\mbox{d}\xi}{\mbox{d}\eta} &= 
\frac{b_mA_1}{A_m}+O(\mu), \label{exp:ML_lf_trajectory} 
\end{align}
for $k_2$ from (\ref{b_tip_A0_hf}).
From (\ref{exp:ML_lf_concavity}) we see that the line corresponding to the concavity threshold, across which
the trajectory of $x$ changes its concavity, is given by
\begin{eqnarray}
\xi_c^{\rm ML}=- \frac{{\cal C}^2A_m}{2k_2 A_1}\eta + \frac{1}{2{\cal C}k_2b_m} \, . \label{concavity_ML}
\end{eqnarray}
We compare the slope of (\ref{concavity_ML}) to the
slope of $x$  from (\ref{exp:ML_lf_trajectory}),
as in (\ref{sol:lf_concavity_threshold})-(\ref{sol:lf_A_threshold}).
Where these two slopes are equal, we find the change in concavity of the trajectory that corresponds
to tipping. Setting the slope of (\ref{concavity_ML}) equal to (\ref{exp:ML_lf_trajectory}),  and solving for $A_1$ yields

%
\begin{equation}
A_1 = \frac{{\cal C}A_m}{\sqrt{2|k_2|b_m}}\qquad \Rightarrow \qquad
A_c = A_m +\frac{\Omega A_m}{\sqrt{2|k_2|b_m}}
\label{sol:ML_lf_Ac}
\end{equation}
\begin{figure}[h]
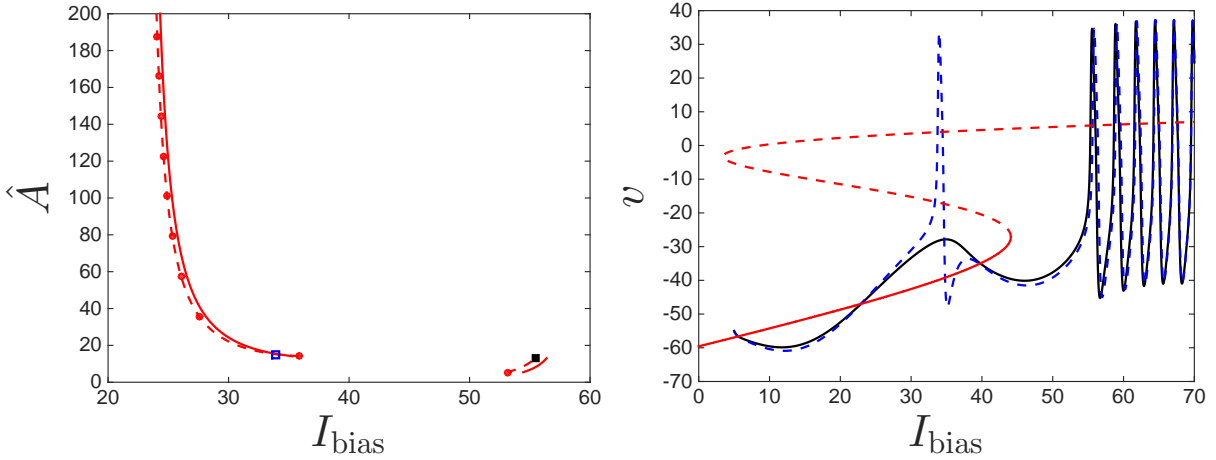

\begin{center}
\includegraphics[width=0.45\textwidth]{figure_ML_lf_tipping.eps} 
\includegraphics[width=0.45\textwidth]{figure_ML_lf_track.eps}
\end{center}
\caption[ML Low Frequency Case]{\label{fig:ML_lf_tipping} Results are shown for $\mu=.04$, $\Omega=0.2$. LEFT: The location of the tipping points 
obtained by numerical simulations of (\ref{sys:ML}) (Red dashed line
 with dots) as a function of
amplitude $\hat{A}$, compared with the asymptotic approximation given by (\ref{sol:ML_lf_tipping}) (Red solid), with tipping defined as $v=|v_c|$. RIGHT: Numerical simulations for different amplitudes $\hat{A}$, $-\hat{A}\approx 13.02$ (black solid) and $-\hat{A}\approx 15.20$ (blue dashed). The related tipping points are shown as the black solid square and blue open square in LEFT.}
\end{figure}

Figure \ref{fig:ML_lf_tipping}(LEFT) compares the asymptotic approximation of
 the location of the tipping points given by (\ref{sol:ML_lf_tipping}) and
(\ref{sol:ML_lf_Ac}) to the results from numerical simulations.
We show the tipping point as a function of the amplitude $\hat{A}$, with the jump in this curve 
occuring at the amplitude threshold $A_c$ given 
by (\ref{sol:ML_lf_Ac}).
 Figure \ref{fig:ML_lf_tipping}(RIGHT) illustrates
the character of the  early transition near $A=A_c$ for (\ref{sys:ML}) with early tipping occuring for
$A>A_c$ but not for $A<A_c$. For the value of $I_{\rm bias}$ where early tipping occurs, there
is not an attracting oscillatory solution nearby, so 
following the tipping event, the system does not remain in the oscillatory state. 
This behavior is in contrast to the system 
(\ref{sys:gen_periodic_forcing}),
for which the tipping is irreversible.  
For some combinations 
of $b_0$ and $A$, the tipping point for (\ref{sys:norm_ML})
 associated with $A_c$ is closer to the SNIC bifurcation, 
in which case after the tipping for $A>A_c$ the system remains in the oscillatory state.



%
%
%


\section{Energy balance model for sea ice}
\label{sec:sea_ice}

We apply the analysis to an energy balance model developed by Eisenman and Wettlaufer \cite{eisenman2009},
as mentioned in the Introduction.    The model describes sea ice dynamics  in terms 
the change of the energy per unit area,  stored either as latent heat in sea ice or as sensible
heat in the ocean mixed layer, 
 \begin{equation}
 E \equiv \left\{ \begin{array}{cc}
  -L_i h_i & E<0  \\
   c_{\rm ml}H_{\rm ml}T_{\rm ml} & E \ge 0,
        \end{array} \right.
        \label{eq:energy}
 \end{equation}
where $L_i$ is the latent heat of fusion for sea ice, $h_i$ is the thickness of sea ice, $c_{\rm ml}$ is the heat capacity of the ocean mixed layer, $H_{\rm ml}$ is the depth of the mixed layer, and $T_{\rm ml}$ is the temperature difference
from the freezing point $0 ^{\circ} \mbox{C}$. 
 The reduced system,  derived from the full heat conduction equations for sea ice thermodynamics of \cite{maykut1971},
expresses the energy dynamics in terms of sensible and latent heat fluxes,
  outgoing longwave and incoming shortwave radiative surface fluxes.
We focus on the partially linearized version model of
 \cite{eisenman2009} that neglects ice export
 which has the form
\begin{equation}
\frac{\mbox{d}E}{\mbox{d}t}  = [ 1-\alpha(E) ]F_S(t) - F_0(t) + \Delta F_0 -F_T(t)T(t,E) +F_B. \label{sys:sea_ice}
\end{equation}
 Here $F_S(t)$ describes the shortwave radiation flux which is seasonally varying, 
with $[1-\alpha(E)]F_S(t)$ describing surface shortwave radiation related to albedo feedback, 
\begin{equation}
\alpha(E) = \frac{\alpha_{\rm ml}+\alpha_i}{2}+\frac{\alpha_{\rm ml}-\alpha_i}{2} \mbox{tanh}\left(\frac{E}{L_i h_{\alpha}}\right), \label{eq:albedo}
\end{equation}
\noindent where $\alpha_{\rm ml}$ and $\alpha_i$ are constants for albedo feedback for the mixed layer and surface ice, respectively. The factor $h_{\alpha}$ is an ice thickness range for smooth transition from $\alpha_i$ to $\alpha_{\rm ml}$.
The temperature dependence of the surface flux is obtained
via a linearization of the Stefan-Boltzmann equation, resulting in the
 outgoing longwave radiation contributions $F_0(t)+F_T(t)T(t,E)$, with
 the seasonal variation in $F_0(t)$ and $F_T(t)$ determined from combined modeling and observational data (see SI Appendix of  \cite{eisenman2009} for details).
 The partial linearization given in  \cite{eisenman2009} follows from taking
the surface temperature $T(t,E)$ as
\begin{equation}
T(t,E) = \frac{E}{c_{\rm ml}H_{\rm ml}} \, . \label{eq:linear_tem}
\end{equation}
The term $F_B$ describes the heat flux into the bottom of the sea ice or the ocean mixed layer. The additional
contribution $\Delta F_0$ captures a change in the surface heat flux, representing warming in the model for $\Delta F_0>0$.

\begin{figure}[h]
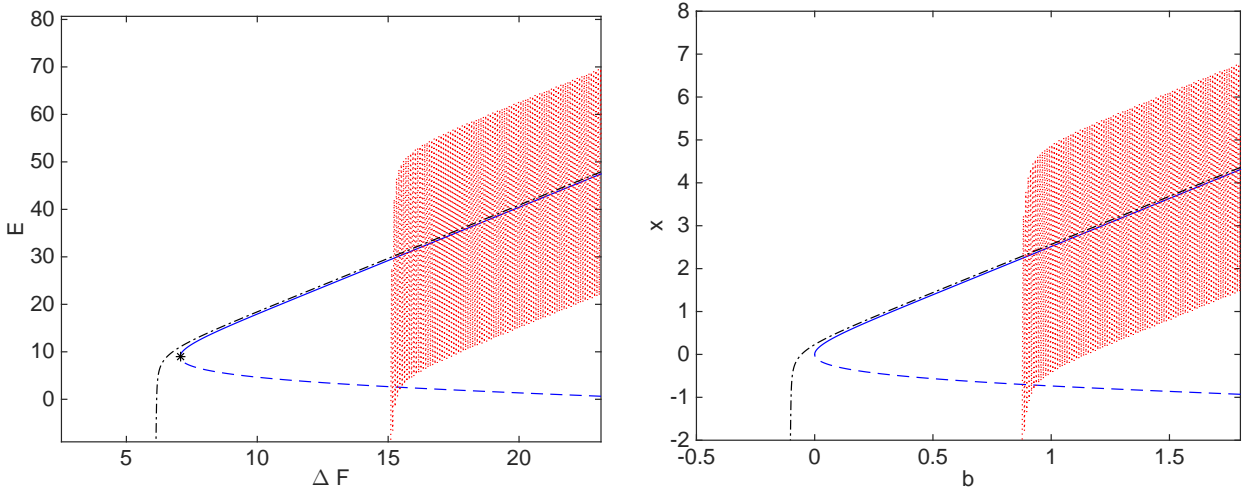

\begin{center}$
\begin{array}{cc}
\hskip -.5in
\includegraphics[width=0.45\textwidth]{E_seaice.eps} &  \ \
\includegraphics[width=0.45\textwidth]{x_seaice.eps}
\end{array}$
\end{center}
\caption{
Time series in both panels are shown for $\mu=.01$. 
LEFT:
Time series for $E$ in (\ref{sys:toy_model}) with oscillatory forcing term $g_j$ (red dotted line)
 and without, that is,
where the oscillatory terms $g_j$
are replaced by their averages $G_j$ 
 (black dash-dotted line). 
 The critical value ($\Delta {F_0}_c, E_c$) is marked with *.
The solid (dashed) blue curve corresponds to stable (unstable) branch of
 $dE/dt=0$ in (\ref{sys:toy_model}) 
 for the averaged model.
RIGHT: 
Time series for $x$ given in (\ref{sys:h_model}) with oscillatory forcing term $g_j$ (red dotted line)
 and without, that is, where the oscillatory terms $g_j$
(black dash-dotted line),
are replaced by their averages $G_j$.
The solid (dashed)  blue curve corresponds
to the stable (unstable) branch of $b=-H(x)$ in (\ref{sys:h_model}).
 When rescaled using the transformation (\ref{xE_transform}), the oscillations for the $x$
approximation are slightly smaller than those for the full $E$ equation, 
although the difference is negligible on this scale.  }
\label{fig:Ex_trajectory}
\end{figure}

In \cite{eisenman2009} the nonlinear threshold behavior of the energy was studied through the bifurcation structure of (\ref{sys:sea_ice}) with the seasonally varying terms $F_0(t)$, $F_S(t)$ and $F_T(t)$. Attracting periodic states for fixed $\Delta F_0$ were determined by simulating the equation, allowing the system to reach the attracting
state, and recording the minimum and maximum value of $E$ for the particular
value of $\Delta F_0$.
For the parameters in \cite{eisenman2009},
  there are two branches of attracting states for $E>0$ and $E<0$ 
  in terms of $\Delta F_0$, the bifurcation parameter describing the change to the surface 
heat flux.  These states are bistable for a range of $\Delta F_0$, with each
branch losing stability to the other via a saddle node bifurcation.
Figure \ref{fig:Ex_trajectory} shows the branch corresponding to the stable 
steady state of the averaged version of (\ref{sys:sea_ice}), that is, 
with the seasonally varying terms $F_0(t)$, $F_S(t)$ and $F_T(t)$ replaced by their
averages, and for
 $E\geq0$ and $\Delta F_0$ a positive constant.
We denote the saddle node bifurcation value on this branch as
  $(\Delta {F_0}_c, E_c)$, and also show the unstable branch for that system.
As would be expected from the analysis of previous sections, the saddle node
bifurcation for the averaged system differs considerably from the tipping point
when there is large amplitude periodic forcing.
The parameters values are as in \cite{eisenman2009}: 
$L_i= 9.4\ {\rm Wm}^{-3}$yr, $c_{\rm ml}H_{\rm ml}= 9.4 {\rm Wm}^{-2}{\rm yr K}^{-1}$, $\alpha_i= .68$, $\alpha_{\rm ml}= .2$,
$F_B = 2 {\rm Wm}^{-2}$, and 
$h_{\alpha} = .5$m.
The seasonally varying quantities take values in the ranges $130\ {\rm Wm}^{-2}>F_0(t)>54\ {\rm Wm}^{-2}$, $3.3\ {\rm Wm}^{-2}{\rm K}^{-1}>F_T(t)>2.5\ {\rm Wm}^{-2}{\rm K}^{-1}$, and $310\ {\rm Wm}^{-2}>F_S(t)>0\ {\rm Wm}^{-2}$.

In considering the case where the surface heat flux varies slowly over time
with a rate $\tilde\mu$ we use the asymptotic approaches described 
in previous sections to study 
tipping from the attracting branch
$E \geq 0$, 
noting that by symmetry
a similar approach applies for the other branch $E<0$.  
We contrast
 the location of the tipping point for the averaged model with 
the tipping location in the case where the oscillatory functions 
$F_0(t)$, $F_S(t)$ and $F_T(t)$ are not replaced with their averages. Since we
obtain analytical expressions for the tipping point, we can easily explore
the parametric dependence on the location of the tipping point, which has implications for the bifurcation structure as a function of the parameters. We note that the 
computational approach used in \cite{eisenman2009} used to explore
bifurcations corresponds to the limit
$\tilde\mu \to 0$.

%

The fluctuations in $F_0(t)$, $F_S(t)$ and $F_T(t)$ are
given as a series of monthly measurements \cite{eisenman2009} (SI Appendix), which
we interpolate  and approximate  using a (finite)
 Fourier series approximation. We assume that the fluctuations are
the same year to year. 
To allow easy identification of the oscillatory forcing terms,
we write (\ref{sys:sea_ice}) in a form analogous to
that of (\ref{sys:gen_periodic_forcing}) and
shift the  bifurcation point for the averaged system from
 $(\Delta {F_0}_c,{E}_c)$ 
 to the origin 
as shown in Figure \ref{fig:Ex_trajectory}.
We rescale the system appropriately, so that the resulting dependent variable $x$
and corresponding slowly varying parameter $b$ are $O(1)$ quantities.  
Substituting
\begin{eqnarray}
   x = \frac{E-{E}_c}{{E}_c}, \quad  b = \frac{\Delta F_0-\Delta {F_0}_c}{{E}_c} \label{xE_transform}
\end{eqnarray}
into (\ref{sys:sea_ice}), and taking
$\Delta F_0$ to vary slowly in time with rate $\tilde\mu$ yields
\begin{align}
\frac{\mbox{d}x}{\mbox{d}t} &= g_1(t)+g_2(t)\mbox{tanh}(g_3(x+1))+g_4(t)x+b,  \label{sys:toy_model} \\
\frac{\mbox{d}b}{\mbox{d}t} &= -\mu = -\tilde\mu/E_c. \nonumber\\
 g_1(t) & = E_c^{-1}\left[(1-\frac{\alpha_{\rm ml}+\alpha_i}{2})F_S(t)-F_0(t)+F_B\right]
+\frac{\Delta {F_0}_c}{{E}_c}+g_4(t),\nonumber\\ 
  g_2(t) & =  -\frac{\alpha_{\rm ml}-\alpha_i}{2E_c} F_S(t), 
\qquad g_3 = 
 \frac{1}{L_i h_{\alpha}}{E}_c, 
\qquad g_4(t) =  -\frac{F_T(t)}{c_{\rm ml}H_{\rm ml}} . \nonumber
\end{align}

In contrast to 
 (\ref{sys:gen_periodic_forcing}), which has only additive periodic forcing,
  (\ref{sys:toy_model}) 
has both multiplicative oscillatory forcing terms with coefficients
 $g_2(t)$ and $g_4(t)$ and an additive oscillatory term $g_1(t)$. So we consider approximations for the terms involving $g_2$ and $g_4$
that replace the multiplicative oscillations with reasonable additive oscillations.
  The oscillations 
in $g_1$ have an amplitude of $\approx 20$, $g_2(t)$ has an amplitude of $\approx 5$, while $g_4(t)$
has oscillations with relatively small amplitude of $\approx 0.5$. We  
anticipate that the fluctuations of $g_4(t)$, as compared with those of $g_1(t)$ 
and $g_2(t)$, have a negligible effect on tipping, generally occurring for
$x=O(1)$ or smaller.
So  
 we neglect the oscillations in the term $g_4(t) x$ by replacing $g_4(t)$ with its average.
Furthermore the results from the previous sections suggest that  large oscillations in
$g_1$ drive  tipping 
for $O(1)$ values of the bifurcation parameter $b$ and thus for $O(1)$ values of $x$,
 for which $\mbox{tanh}(g_3(x+1)) \approx 1$.
Then, it is reasonable to approximate
 $ g_2(t)\mbox{tanh}(g_3(x+1)) \approx  G_2\mbox{tanh}(g_3(x+1)) + g_2(t)-G_2$, so
that (\ref{sys:toy_model}) is approximated by
\begin{align}
\frac{\mbox{d}x}{\mbox{d}t} &= b+H(x)+q(T),  \label{sys:h_model} \\
\frac{\mbox{d}b}{\mbox{d}t} &= -\mu. \nonumber\\
H(x) &= G_1+G_2 \mbox{tanh}(g_3(x+1))+G_4 x, \quad  q(T) = g_1(t)-G_1 + g_2(t)-G_2 \, ,
\end{align}
 where  $G_j$ is the average value of $g_j(t)$
 over one period for $j=1,2,4$ and 
  $q(T)$ is an oscillatory function with mean zero and
period of $t=1$. With the oscillatory contributions represented as Fourier series, 
the time scale of $q(T)$ is $T=\Omega t$ with  $\Omega = 6.28$. 
The computed trajectories shown in
 Figure \ref{fig:Ex_trajectory} illustrate the validity of
 (\ref{sys:h_model}) as an approximation for
(\ref{sys:toy_model}).


We note here that for the transformed system, as in Section \ref{sec:model_ML}, the bifurcation
parameter $b$ and dependent variable $x$ take $O(1)$ values. Then
we identify different scales related to slow drift, high and
low frequency, or large amplitude in this setting, which again
  is important in choosing the appropriate asymptotic approach.

An analytical expression for the tipping point in the case $q(T) =0$ is found
as described in (\ref{hpoly0})-(\ref{a_g}), based on 
a Taylor series of $H(x)$ in (\ref{sys:h_model}) about $x=0$.
Then
\begin{eqnarray}
b_d = |H_2|^{-1/3}a_d- b_p \qquad \mbox{for} \qquad b_p = H_0 + \frac{H_1^2}{4|H_2|}, \label{hpoly0_tip} \qquad
H_j = j! H^{(j)}(0) \nonumber 
\end{eqnarray}
is the tipping point for (\ref{sys:h_model}) for $q(T)=0$, and is
 compared below with the tipping point for $q(T)\neq 0$.   We note that $H_0$ and $H_1$
are small, but not identically zero as are the analogous terms in Section 
\ref{sec:model_ML}, since we are using an approximate
system (\ref{sys:h_model}) rather than the full model from which the saddle node bifurcation point
 $(\Delta {F_0}_c,E_c)$
was determined.

Next, we consider whether (\ref{sys:h_model}) corresponds to 
high frequency forcing 
as in Section \ref{sec:hf_slow_drifting}, or low frequency forcing as in
Section \ref{sec:lf}.  For $\mu\ll 1$
and $\Omega = 6.28$,
 we would expect this system to be in the high frequency case. 
However, since the amplitude of the oscillations is large, 
Section \ref{sec:hf_large_A} indicates that
$\log(A)/\log(\Omega)$ must be considered  to determine if oscillations  with amplitude $A$
 correspond to high or low frequency forcing in a rescaled system. We find that 
$4/3<\log(\max(q(T)))/\log(\Omega)<2$ for
 (\ref{sys:h_model}) falls in between the asymptotic ranges identified 
in Section \ref{sec:hf_large_A} for the high and low frequency cases.
 Nevertheless, we adapt the concepts from
the previous sections to approximate the shift in the tipping point 
driven by the oscillations in (\ref{sys:h_model}).

We begin with an outer approximation for $x$ based on the approach from Section \ref{sec:hf_slow_drifting} and
Appendix \ref{append:C}.
Substituting a multiple scale approximation,
$x(\mu t,T) = x_0(\mu t,T)
+ \Omega^{-1} x_1(\mu t,T) + \ldots$ for
$T= \Omega t$
 into (\ref{sys:h_model})
we get
\begin{eqnarray}
 & {x_0}_T = 0, \Rightarrow  x_0 = x_0(t)
 & {x_1}_T =  b+H(x_0) + q(T) \, .
\end{eqnarray}   
The solvability condition (\ref{eq:solv_condition}) for $x_1$ yields $b+H(x_0)=0$ so that the leading order contributions to $x$ are
\begin{eqnarray}
 x \sim  H_+^{-1}(-b) + \Omega^{-1}Q(T)\, , \quad Q'(T) = q(T), \label{sea_ice_outer}
\end{eqnarray}
where $H_+^{-1}(-b)>0$ corresponds to the stable branch for $x>0$
and $Q(T)$ has zero average.
The leading order approximation for the outer solution (\ref{sea_ice_outer}) is  composed of oscillations $\Omega^{-1}Q(T)$ 
about the branch of the averaged system $b=-H(x)$.

 Where the outer approximation is no longer valid, we construct a local  expansion. To find an
appropriate local approximation,  we appeal to the results shown for the canonical model in 
Section \ref{sec:hf_slow_drifting}, noting that we must adapt that approach 
since the amplitude of the oscillatory forcing is large. 
 If the amplitude was $O(1)$, we would write the inner solution in the multiple scale form
 $x = \Omega^{-1}Q(T) + Y(t,T)$ analogous to the local expansion in Section \ref{sec:hf_slow_drifting}. 
This form is
an approximation composed of oscillations near $x=0$ with $Y$ a small correction to be determined.
However, since the oscillations are large for this application, we expect that the outer approximation 
 breaks down away from $x=0$. In fact, if we examine the outer approximation (\ref{sea_ice_outer}), we see that 
it reaches values of $x$ that  are on the unstable lower branch shown in Figure \ref{fig:Ex_trajectory}. 
We can approximate the value of $x$ where this occurs by identifying a value $b^*$ and corresponding $x^* = H^{-1}_+(-b^*)$
on the stable  branch of $b= -H(x)$ where
\begin{eqnarray}
  H_+^{-1}(-b^*) + \Omega^{-1}\min(Q(T)) = H_-^{-1}(-b^*) \, , \label{bc_def}
\end{eqnarray}
where $H_-^{-1}(-b^*)$ is on the unstable branch shown in Figure \ref{fig:Ex_trajectory}. 
Then we adapt the approach from  Section \ref{sec:hf_slow_drifting}
 writing the solution in the multiple scale form that has oscillations about $x^*$ with $Y$ a correction,
\begin{eqnarray}
 x = x^*+\Omega^{-1}Q(T) + Y(T,t). \label{xc_local}
\end{eqnarray}
We could also rescale $b$, $Y$ and the time
scale $t$ with a power of $\mu$.  However, we obtain the same result without this rescaling, 
so for simplicity we do not introduce such rescaled variables here.
Substituting the  multiple scale expression (\ref{xc_local}) into (\ref{sys:h_model}) yields 
\begin{eqnarray}
\Omega Y_T &=& -Y_t + b + H(x^* +\Omega^{-1} Q + Y) \label{Yeq1}\\
 &\sim & -Y_t + b + H(x^* +\Omega^{-1} Q) +H'(x^* +\Omega^{-1}Q) Y + \frac{H''(x^*+\Omega^{-1}Q)}{2}Y^2 \, . \label{Yeq2}
\end{eqnarray}
 In contrast to the approximation used in  (\ref{exp:ML_inner_LO1})
 for the Morris-Lecar model or in (\ref{hpoly0}),
 we have not replaced $H(x)$ with a  polynomial in $x$ before substituting 
(\ref{xc_local}).  The reason that we avoid this step is that
the forcing $q(T)$ is relatively large in this case, so we must consider
$H(x)$ over a larger range of $x$ values. Over this range, it is not possible to
approximate $H(x)$ accurately with a Taylor series about a single value of $x$. 
Rather we have substituted (\ref{xc_local}) directly into $H(x)$ and then expanded about $Y=0$,
assuming that $Y$ is a small correction. 
Then the solvability condition (\ref{eq:solv_condition}) for (\ref{Yeq2}) yields the equation for  $Y$ 
\begin{eqnarray}
 Y_t &=&  b + \left(\lim_{L \to \infty}\frac{1}{L}\int_0^{L}H(x^*+\Omega^{-1} Q(T)) dT\right) 
+\left(\lim_{L \to \infty}\frac{1}{L}\int_{0}^{L} H'(x^*+\Omega^{-1}Q(T)) dT\right) 
\, Y \nonumber\\
  & &+ \left(\lim_{L\to \infty}\frac{1}{L}\int_0^{L}\frac{H''(x^*+\Omega^{-1}Q(T))}{2} dT \right) \, Y^2 \nonumber\\
  &=&  b + {\cal H}_0 + {\cal H}_1 Y + {\cal H}_2 Y^2 \, .\label{Ysolv}
\end{eqnarray}
 From  (\ref{bc_def}) we get
$x^*\approx 2.2$.
Then we obtain the tipping point from the singularity for (\ref{Ysolv}), using the results
from (\ref{hpoly0}) and 
\begin{eqnarray}
b_{\rm tip} &=& |{\cal H}_2|^{-1/3}a_d- b_Q \quad \mbox{for} \qquad b_Q = {\cal H}_0 + \frac{{\cal H}_1^2}{4|{\cal H}_2|} \, . \label{btip} 
\end{eqnarray}
In Figure 14 we compared this approximation 
 with the numerical results in terms of the original
variables $E$ and $\Delta F_0$. As would be expected for a system forced
by oscillations with a large amplitude to frequency ratio, the advance of the tipping point
is dominated by the oscillations. As $\tilde\mu \to 0$, the tipping point approaches the
end of the branch of attracting oscillatory solutions computed in \cite{eisenman2009}.

\begin{figure}[h]
\begin{center}$
\begin{array}{cc}
\includegraphics[width=0.45\textwidth]{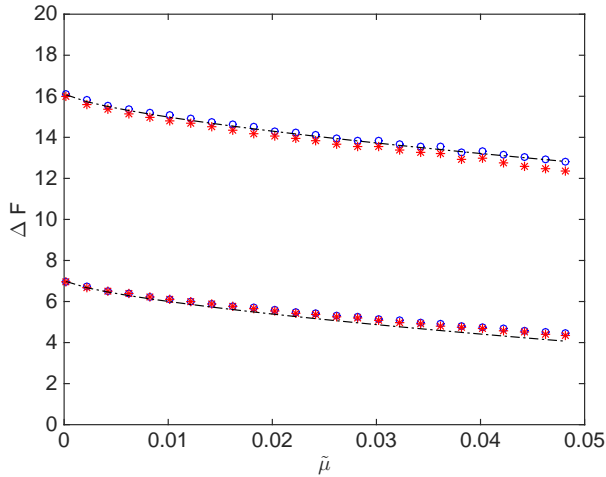}
\end{array}$
\end{center}
\caption{Comparison of the tipping points $\Delta F_0$ obtained
from the asymptotic approximation (\ref{btip}) (black dash-dotted lines),
compared with the numerically obtained values from the full model for
$E$ (\ref{sys:toy_model}) (blue o's) and from the approximate model for $x$ (\ref{sys:h_model}) (red *'s).
The lower values correspond to the tipping point for the averaged
system with $g_j(t)$ replaced with $G_j$, while the upper curve
corresponds to the case where the oscillations are not replaced
with their averages.}
\label{fig:DeltaFtip}
\end{figure}

The fact that we are able to obtain an analytical expression for the tipping
point allows us to explore the impact of parameters values
 on the bifurcation structure of the model 
 and the location of transitions between states. 
Specifically, we can explore the
 existence or prominence of hysteresis between different 
states as a function of the parameters of a periodically forced system. 
Within the context of a large scale
 atmosphere-ocean global climate model studied in \cite{armour2011}, the 
question was considered whether or not parameter variation yields
  hysteresis in the 
transitions between states analogous to $E>0$ and $E<0$ in (\ref{sys:sea_ice}).
 There it was seen that
 the hysteresis obtained for time-varying parameters, predicted by lower dimensional models where there
is bistability for different states, is less prominent and 
in certain contexts it is not observed in the large scale computations.
An analogous question of how the bifurcation structure 
could change with parameters
was explored computationally in \cite{eisenman2012}
for a simplified version of (\ref{sys:sea_ice}).
For example,  recomputing the bifurcation diagram for different parameter
values, \cite{eisenman2012} shows a loss of bistability with  
increased contributions to longwave radiation from $F_T(t)$,
and thus a loss of the hysteresis loop observed for a varying bifurcation
parameter.
In the context of tipping, the loss of
hysteresis corresponds
to shifts in the tipping points so that regions of bistability are no 
longer observed.
The value of the analytical expression for the tipping point (\ref{btip})
 is in determining shifts in the tipping location 
 as a function of a parameter, rather than having to recompute the 
entire bifurcation structure through repeated
 simulation of the full model over a potentially large range
 of parameter values.

\begin{figure}[h]
\begin{center}$
\begin{array}{cc}
\includegraphics[width=0.60\textwidth]{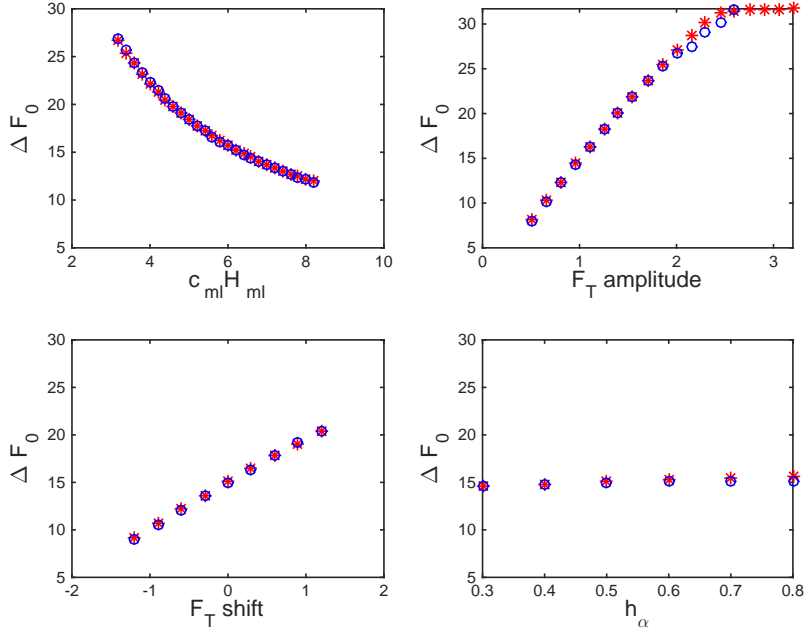}
\end{array}$
\end{center}
\caption{Analytical approximation (blue circles) and
numerical approximation (red *'s) for the value $\Delta F_0$ corresponding to
tipping as a function of different parameters: the ratio of $E$ to $T(t,E)$
in (\ref{eq:linear_tem}) (UPPER 
LEFT), a multiplicative factor for the amplitude of the oscillations $F_T$
(UPPER RIGHT), a shift in the mean of $F_T$ (LOWER LEFT), and the ice thickness parameter describing transitions in the albedo (LOWER RIGHT).
 In the UPPER RIGHT figure,
the approximation is not valid for tipping points with values $\Delta F_0 > 31$. This corresponds to larger amplitudes of $F_T$, for which  
there is no overlap in the range of $\Delta F_0$ for upper and lower attracting 
branches of $E$, indicating no hysteresis for varying $\Delta
F_0$.  Then determining  $x^*$
from (\ref{bc_def}) is no longer valid for the asymptotic approximation.
In that case the numerical approximation asymptotes to $\Delta F_0 \approx 31$. }
\label{tip_parvar}
\end{figure}

In Figure \ref{tip_parvar} we demonstrate how this loss of hysteresis
is obtained directly from the analytical expression (\ref{btip})
for the tipping point.  
 The UPPER RIGHT panel  shows that
increasing the amplitude of $F_T(t)$ causes earlier tipping.
 For tipping points increased to $\Delta F_0\approx 31$, there is
a loss of hysteresis, since tipping  occurs for approximately the
same value of $\Delta F_0$ for which there is a transition from the
lower solution branch for $E$ (not shown) to the upper branch in 
Figure \ref{fig:Ex_trajectory}.   
The approximation (\ref{btip}),  based on (\ref{bc_def}),
indicates this behavior for the amplitude of $F_T$ increased above
a factor of about 2.7, as also observed in \cite{eisenman2012}. 
 For this increase in the oscillations there is no overlap in $b$ for the upper and lower attracting 
branches for the solution $x$, and no intermediate unstable branch.
 Then determining  $x^*$
from (\ref{bc_def}) is no longer valid, signaling the loss of bistability and of hysteresis.
For the same reason, the numerical approximation asymptotes to the value 
of $\Delta F_0 \approx 31$, as shown in Figure \ref{tip_parvar}(UPPER RIGHT).
From the normalized model (\ref{sys:h_model}) we can also identify other
parameters $c_{\rm ml}H_{\rm ml}$ and the mean of $F_T$ that can drive
a considerable shift in the tipping point.  The LEFT panels of Figure
\ref{tip_parvar} show the effect of varying these parameters, which is not as large as for the increased
amplitude of $F_T$, as also observed in \cite{eisenman2012}.

The variation of the tipping point with $h_{\alpha}$ is also shown. For
this model, approximating the tipping location using 
 (\ref{bc_def})-(\ref{btip}) is based on the existence of an
upper solution branch.  
Then Figure \ref{tip_parvar}(LOWER LEFT) shows the approximation for
values $h_\alpha$ where $x^*$ in (\ref{bc_def}) can be determined.
The upper branch used in (\ref{bc_def}) no longer exists for $h_\alpha$ sufficiently small, as also seen in \cite{eisenman2012}.

\vspace{1cm}

%
%
%
%


\section{Summary and Future Work}

 In this paper we have analyzed the factors that determine the location of a tipping point in the
canonical model for a saddle node bifurcation with  bifurcation parameter $a$ that varies at a slow rate $\mu\ll 1$ and
   forced by a periodic oscillation with amplitude $A$ and frequency $\Omega$.
  We see that the main contributions that determine the tipping point location are different for
 low and high frequency. In order to get analytical expressions to identify the key factors in the location of the tipping point,
  we apply the methods
of multiple scales and matched asymptotic expansions to obtain analytical expressions for the combined effect of $\Omega$, $A$, and $\mu$.
 The cases of low and high frequency in the oscillatory forcing naturally require the consideration of different relationships
between the time scales 
of the oscillation and the slow drifting rate.  The size of the  amplitude also plays an important role.

For a high-frequency forcing, there are two competing contributions to the location
of the tipping point relative to the static bifurcation point at $a=0$. 
 The first contribution is a delay due to the slow drifting rate $\mu$  of the bifurcation parameter $a$, well-known in a variety of studies without
periodic forcing \cite{mandel1987, haberman1979}. The second contribution is an advance proportional to the square of the ratio of the amplitude to frequency for the oscillation, observed also in the case where the bifurcation parameter is constant ($\mu =0$). These contributions reflect that the trajectory in this case is 
approximated by the linear combination of the slowly varying equilibrium solution found for $A=0$ in \cite{haberman1979} and an additive contribution that 
takes the form
of the oscillatory forcing with a scale factor. Near the tipping point, the trajectory can be expressed in terms of an oscillatory part plus a correction that satisfies an Airy function on the slow time
scale, with the oscillations averaged on this time scale shifting the singularity that corresponds to the tipping point. The resulting
expressions for these contributions then give an expression for the location of the 
tipping point that varies continuously with the parameters.  Whether 
the combined effect results in an advance or delay compared to
$a=0$ depends on the magnitude of $\mu^{2/3}$ relative to $A^2/\Omega^2$. 
In particular, we see noticeable advances of the tipping point for combinations
of larger amplitudes with frequencies that are not very large. 

Advances in the tipping point due to periodic forcing are even more apparent for lower frequencies. For $\Omega\ll 1$, the
trajectory is dominated by a nonlinear combination of the slowly drifting bifurcation
parameter and the oscillatory forcing, since they give contributions that are of the same order of magnitude. Then the location of the tipping point
depends not only on the amplitude, frequency, and drifting rate, but also on the phase of the oscillation. 
Tipping occurs near values of $a$ where the trajectory is near $x\sim 0$, where a local analysis reveals that
 the change of concavity for $x \sim 0$ provides an analytical expression for the tipping point. This expression captures
 the discontinuous dependence of the tipping location on parameters such as the amplitude. 
These discontinuities follow directly from the fact that tipping occurs only for values of $a$ where the trajectory is near zero: a 
slight change in parameters can change the approach to zero of the trajectory, postponing tipping to a later point at which the trajectory crosses
zero. This shift in the tipping point is on the order of one period of the low frequency oscillation.

The development of this suite of approaches provides a foundation for capturing the behavior
of the tipping point in terms of the key parameters.
Within the context of the applications in Sections \ref{sec:model_ML} and \ref{sec:sea_ice}, we have demonstrated how the approximations developed for the canonical model can be used directly in other models. The influence of the oscillations, and thus the choice of the appropriate approximation, depends on
certain scaling relationships between drifting rate, amplitude, and frequency. These are
identified
within the context of the canonical model where the dependent variable and bifurcation parameters
 take $O(1)$ values
and the saddle node bifurcation at the origin. Thus an important step is normalization of the model to this standard form, from 
which the parameter relationships can be extracted.  Within the transformed model one can already identify parameter ranges where
advanced tipping is likely to occur, by relating the results of the canonical model to the rescaled parameters.  

The results for both the canonical model and the applications illustrate how the concepts developed within the asymptotic approaches can be extended to 
wider ranges of parameters.  For example, the results of Section \ref{sec:lf_order1}  for low frequency, where the frequency is of the same order of magnitude as the drifting rate, are extended in Section \ref{sec:lf_order_nu} for larger values of frequency that are $O(1)$, by using the values where the regular asymptotic 
expansion no longer exists as a basis for a local analysis of the change of concavity.  Another example
of the extension of the method is in the application
of Section \ref{sec:sea_ice}, where the ratio of amplitude to frequency for the oscillation is outside of the asymptotic range identified for large
amplitude oscillations in Section \ref{sec:hf_large_A}.  Nevertheless, as for the canonical model, we apply a regular or outer expansion, and from conditions where the outer expansion loses its validity
we identify the region for a local multiple scale analysis.  The solvability condition for
the local approximation again provides an averaged equation with
quadratic nonlinearity, whose singularity provides the approximation of the location of the tipping point.

We have also compared the advance of the tipping point for high frequency forcing with the shift due to additive white noise forcing, where a significant 
advance of the tipping point occurs if
the magnitude of the noisy forcing scales as the square root of the drifting rate. In order to obtain a comparable shift with high frequency forcing 
alone, a much larger amplitude is required.  In contrast, for low frequency forcing, it is possible to obtain a larger advance than in the
case of high frequency forcing.  Since white noise forcing includes both high and low frequencies, the conditions for advancing the tipping point as indicated
in the analysis of low frequency forcing are captured in the noisy forcing.

 There are a number of remaining questions about the influence of oscillatory and noisy forcing on a system with a slowly drifting bifurcation.
  Here and in \cite{sieber2012,berglund2010}  the case of
additive forcing has been considered, but in general the oscillatory and noisy forcing can be multiplicative or parametric, as can be seen
explicitly in the model considered in \cite{eisenman2009} without averaging over seasonal periodicity.
In addition, the combination of noisy and oscillatory forcing should be of interest, as they appear together in applications. For example,
for low frequency oscillations, we see that the location of the tipping point depends on the phase, and we would expect that the noise causes
phase shifts.  The question then remains when and if stochastic phase shifts, in combination with periodic forcing, play a significant role
in the tipping location.  
A related area is the phenomenon of rate induced tipping, where exceeding a critical drifting rate can result in early tipping \cite{ashwin2012}. 
The effect of noisy 
or oscillatory forcing on these types of transitions has not been analyzed.

The shift of the tipping point observed for different types of oscillations indicates that
we must consider their effect on transitions via a saddle node bifurcation.
  For example, a common underlying structure for a  hysteresis loop consists of two branches for 
stable states with a bistability region and loss of stability via saddle node bifurcations. 
If oscillations are included that advance the tipping points,
 the size of the hysteresis loop can be reduced, or the hysteresis can be removed
completely. 

More generally in a slow-fast system where the slow variable plays the role of
a bifurcation parameter for the rest of the system, any oscillatory forcing can change 
the time to transition when a saddle node bifurcation is embedded in the system's structure.
Complex phenomena such
 as bursting dynamics or mixed-mode oscillations typically
appear in slow-fast systems and involve combinations of transitions via
 different types of bifurcation with a 
slow variable viewed as a bifurcation parameter. For phenomena of this type in which
a saddle node bifurcation is present,
the results in this paper can quantify the influence of oscillations on the length of 
the intervals where steady
states are attracting. 
In addition to the saddle node bifurcation considered here, there are other types of bifurcations that are potentially sensitive to 
combinations of slowly drifting parameters and oscillatory and noisy forcings.  This
sensitivity has been studied in other contexts, including
different types of bifurcations \cite{berglund2010,kuske1999, borowski2010}, although
not necessarily from the point of view of considering competitive or combined effects of both oscillations and noise on transitions or bifurcations.

\begin{appendices}

\section{The Higher Order Analysis for $\mu=0$ and $\Omega \gg 1$}
\label{append:A}

We continue the asymptotic analysis in Section \ref{sec:hf_no_drifting} after (\ref{sol:hf_constant_order1}) for $O(\Omega^{-2})$ and $O(\Omega^{-3})$, yielding
\begin{align}
O(\Omega^{-2}):  x_{2T}+x_{1t}& = -2x_0 x_1,  \Rightarrow x_{2T} =R_2(T,t)= -v_{1t}-2x_0 v_1+2 A x_0 \mbox{cos}(T) \label{exp:hf_constant_outer2} \\
O(\Omega^{-3}):  x_{3T}+x_{2t} &=-x_1^2-2x_0 x_2, \Rightarrow x_{3T} =R_3(T,t)=-x_{2t}-x_1^2-2x_0 x_2, \label{exp:hf_constant_outer3}
\end{align}
\noindent By applying the solvability condition (\ref{eq:solv_condition}) to (\ref{exp:hf_constant_outer2}) and (\ref{exp:hf_constant_outer3}) and solving the resulting equations, we have
\begin{align}
O(\Omega^{-2}): v_{1t} &= -2x_0 v_1,  \label{exp:hf_constant_local_2} \\
x_{2T} &= 2A x_0 \mbox{cos}(T), \,\,\, \Rightarrow x_2 = 2\sqrt{a}A \mbox{sin}(T)+v_2(t), \label{sol:hf_constant_outer2} \\
O(\Omega^{-3}): v_{2t} &= -\frac{A^2}{2}-2\sqrt{a}v_2, \label{exp:hf_constant_local_3}
\end{align}
\noindent The attracting solutions of (\ref{exp:hf_constant_local_2}) and (\ref{exp:hf_constant_local_3}) are their stable equilibria, which are $ v_1 = 0$ 
and $ v_2 = -\frac{A^2}{4\sqrt{a}}$.
After substituting these equilibria and (\ref{sol:hf_constant_order1}) and (\ref{sol:hf_constant_outer2}) into (\ref{exp:hf_gen_expansion}), we get (\ref{sol:x_outer}) in Section \ref{sec:hf_no_drifting}.

In order to find the local behaviour of (\ref{sys:gen_periodic_forcing}) with $\mu=0$, we substitute $a = \Omega^{-2}b$ into (\ref{exp:hf_multi_scale}) and get
\begin{equation}
x_T + \Omega^{-1} x_t = \Omega^{-3}b -\Omega^{-1}x^2+\Omega^{-1} A \mbox{sin}(T). \label{exp:hf_constant_local}
\end{equation}
The contributions at each order are
\begin{eqnarray}
O(\Omega^0):& & \qquad x_{0T}=0, \,\,\, \Rightarrow \,\,\, x_0 = x_0(t), \label{exp:hf_constant_local_order0} \\
O(\Omega^{-1}):& & \qquad x_{1T}+x_{0t}=-x_0^2+ A \mbox{sin}(T), \,\,\, \Rightarrow x_{1T} = -x_{0t}-x_0^2 + A \mbox{sin}(T), \label{exp:hf_constant_local_order1} \\
O(\Omega^{-2}):& & \qquad x_{2T}+x_{1t} = 0, \,\,\, \Rightarrow x_{2T}=-x_{1t}, \label{exp:hf_constant_local_order2} \\
O(\Omega^{-3}):& & \qquad  x_{3T}+x_{2t} = b-x_1^2, \,\,\, \Rightarrow x_{3T} = -x_{2t}+b-x_1^2. \label{exp:hf_constant_local_order3}
\end{eqnarray}
Substituting (\ref{exp:hf_constant_local_order0}) into (\ref{exp:hf_constant_local_order1}),  applying the solvability condition (\ref{eq:solv_condition}), and solving the resulting equations,  we have
\begin{align}
x_{0t} &= -x_0^2,\,\,\, \Rightarrow \,\,\, x_0 = 0, \label{sol:hf_constant_local_0} \\
x_{1T} &= A\mbox{sin}(T), \,\,\, \Rightarrow \,\,\, x_1 = -A\mbox{cos}(T)+v_1(t). \label{sol:hf_constant_local_1}
\end{align}
Substituting (\ref{sol:hf_constant_local_1}) into (\ref{exp:hf_constant_local_order2}) and applying the solvability condition (\ref{eq:solv_condition}), we have
\begin{equation}
v_{1t} = 0, \,\,\, \Rightarrow \,\,\, v_1 = d, \,\,\, \Rightarrow \,\,\, x_{2T} = 0, \,\,\, \Rightarrow \,\,\, x_2 = x_2(t), \label{sol:hf_constant_local_11}
\end{equation}
\noindent where $d$ is constant and needs to be determined. Substituting (\ref{sol:hf_constant_local_1}) and (\ref{sol:hf_constant_local_11}) into (\ref{exp:hf_constant_local_order3}) and applying the solvability condition (\ref{eq:solv_condition}), we have
\begin{equation}
x_{2t} = b-\frac{A^2}{2}-d^2, \,\,\, \Rightarrow d^2 = b-\frac{A^2}{2}, \label{sol:hf_constant_local_2}
\end{equation}
\noindent which provides the equilibrium solution of $x_2$. If $b-\frac{A^2}{2}<0$, $x_{2t}$ is always negative and
there is no attracting solution.

By substituting (\ref{sol:hf_constant_local_0}), (\ref{sol:hf_constant_local_1}), (\ref{sol:hf_constant_local_11}) and (\ref{sol:hf_constant_local_2}) into (\ref{exp:hf_gen_expansion}), we have the leading order approximation (\ref{sol:hf_equili_local}) in Section \ref{sec:hf_no_drifting}.

\section{The asymptotic approximation of the solution for $0<\mu \ll 1$ and $\Omega \gg 1$ in the outer region}
\label{append:C}

We substitute (\ref{exp:hf_outer_expansion}) into (\ref{exp:hf_drifting_multi_scale}), yielding
\begin{align}
O(1): x_{0T}&=0, \,\,\, \Rightarrow \,\,\, x_0 = x_0(\tau), \nonumber \\
O(\mu^{\lambda}):  x_{1T} &= a -x_0^2 +A \mbox{sin}(T).   \label{exp:hf_drifting_outer_order1}
\end{align}
\noindent Applying the solvability condition (\ref{eq:solv_condition}) of (\ref{exp:hf_drifting_outer_order1}) and
solving the resulting equations yields
\begin{equation}
\begin{array}{lll}
x_0^2 = a, &\Rightarrow & x_0 = \sqrt{a},\nonumber \\
x_{1T} = A \mbox{sin}(T), &\Rightarrow & x_1 = -A \mbox{cos}(T)+v_1(\tau). \nonumber
\end{array}
\end{equation}
The higher order contributions of (\ref{exp:hf_drifting_multi_scale}) depend on the value of $\lambda$. The two cases are:
\begin{eqnarray}
 & &\mbox{For }0<\lambda \le 1, \qquad
O(\mu^{2\lambda}): \mu^{1-\lambda} x_{0\tau} + x_{2T} = -2x_0x_1. 
\label{exp:hf_drifting_outer_order2}\\
 & &\mbox{For } \lambda > 1, \qquad
O(\mu^{1+\lambda}): x_{0\tau} + x_{2T} = \mu^{\lambda-1} \cdot ( -2x_0x_1). 
\label{exp:hf_drifting_outer_order2b}
\end{eqnarray}
In either case, we substitute $x_0 = \sqrt{a}$ in (\ref{exp:hf_drifting_outer_order2}) 
or (\ref{exp:hf_drifting_outer_order2b}) 
and apply the solvability condition (\ref{eq:solv_condition}), to get
\[ x_{0\tau} = \mu^{\lambda-1} (-2x_0 v_1) \,\,\, \Rightarrow \,\,\, v_1 = \mu^{1-\lambda} \cdot \frac{1}{4 a}. \]
Then the asymptotic approximation of the attracting solution is
the result (\ref{sol:hf_drifting_outer}) of Section \ref{sec:hf_slow_drifting}.

\section{The local approximation for $\Omega \ll 1$, near $a_m$}
\label{append:D}

Writing (\ref{exp:lf_multi_scale}) in terms of the local variables  $\eta$, $\xi$, and $A_m$ by (\ref{local:xietaAm}), we have

\begin{equation}
-\mu \xi_{\eta} = a_m [1-\mbox{cos} (c \mu \eta)]+\mu \eta - \frac{1}{c}\mbox{sin}( c\mu \eta)- \mu \frac{A_1 a_m}{A_m}\mbox{cos}(c \mu \eta)- \mu \frac{A_1}{c A_m}\mbox{sin}(c \mu \eta)-\mu^2\xi^2. \label{exp:lf_slope_local}
\end{equation}
Substituting the expansion $\xi \sim \xi_0 + \mu \xi_1 + \mu^2 \xi_2 + \cdots$, into (\ref{exp:lf_slope_local}), yields
\begin{align}
O(\mu): -\xi_{0\eta} &= -\frac{A_1 a_m}{A_m}, \,\,\, \Rightarrow \,\,\, \xi_0= \frac{A_1 a_m}{A_m}\eta+c_0, \label{sol:lf_local_1} \\
O(\mu^2): -\xi_{1\eta} &= \frac{1}{2}a_mc^2 \eta^2 - \frac{A_1}{A_m}\eta -\xi_0^2, \nonumber \\
\Rightarrow \xi_1 &= -(\frac{1}{2}c^2a_m-\frac{A_1^2 a_m^2}{A_m^2})\frac{\eta^3}{3}+\frac{A_1}{A_m}(1+2c_0a_m)\frac{\eta^2}{2}+c_0^2\eta+c_1, \label{sol:lf_local_2}
\end{align}

\noindent which provides the local approximation (\ref{sol:lf_local}). 
We find $c_0$ by writing (\ref{sol:lf_local}) in terms of $x$ and an intermediate scaled
 variable $a- a_m = {\cal K}\mu^{q}$, for $q<1$ and ${\cal K}$ an $O(1)$ positive constant, to get
\begin{equation}
x = \frac{A_1 a_m}{A_m}\cdot {\cal K}\mu^{q}+\mu c_0 -(\frac{1}{2}c^2a_m-\frac{A_1^2 a_m^2}{A_m^2})\frac{{\cal K}^3\mu^{3q-1}}{3} +\frac{A_1}{2A_m}(1+2c_0 a_m){\cal K}^2\mu^{2q}+\mu^{1+q} c^2{\cal K}+\mu^2 c_1 +\cdots \, , \label{exp:lf_match_in}
\end{equation}
  and match (\ref{exp:lf_match_in}) with the Taylor expansion of the outer approximation (\ref{sol:lf_outer}) around $a_m$ written in terms of ${\cal K}$
 \begin{align}
 x = \mu^{q} \cdot \sqrt{\frac{c^2 a_m}{2}}{\cal K}+\mu^{1-q} \cdot \frac{1}{2{\cal K}}(1-\frac{A_1 a_m}{A_m c \sqrt{2a_m}})+O(\mu^{2q},\mu^{2-3q}). \label{exp:lf_match_out}
 \end{align}
We take $\frac{1}{3} < q < \frac{1}{2}$ so that the outer approximation (\ref{exp:lf_match_out}) is valid. By matching the lower order terms in (\ref{exp:lf_match_in}) and (\ref{exp:lf_match_out}), we find $c_0$ in the form of $c_0 = -\frac{1}{2 a_m}+C_1+C_2$.  This form is convenient for comparing with the $\xi$-intercept $-1/(2a_m)$ of
the concavity threshold (\ref{sol:lf_concavity_threshold}), reducing this comparison to
determining the sign of $C_1+C_2$, given by
\begin{align*}
C_1 &= \mu^{-q} \frac{A_m}{{\cal K} a_m A_1}(\sqrt{\frac{c^2 a_m}{2}}-\frac{A_1 a_m}{A_m})+\mu^{1-3q} \frac{A_m}{2{\cal K}^3 a_m A_1}(1-\frac{A_1 a_m}{A_m c\sqrt{2a_m}}), \nonumber \\
C_2 &= \mu^{q-1} \frac{A_m}{3a_m A_1}(\frac{1}{2}c^2a_m-\frac{A_1^2 a_m^2}{A_m^2}).
\end{align*}
We find $c_0$ for $A>A_m$ only, since $A_c>A_m$ by the definition of $A_m$ and $a_m$.
If $\frac{A_1a_m}{A_m}<\frac{c^2A_m}{2A_1}$ and $A>A_m$,
both $C_1$ and $C_2$ are positive and $c_0>-\frac{1}{2a_m}$. Then the local trajectory (\ref{sol:lf_local}) is in the concave up region as assumed for $A_c> A>A_m$ and  $a\leq a_m$ ($\eta\leq 0$).
If $\frac{A_1a_m}{A_m}>\frac{c^2A_m}{2A_1}$, both $C_1$ and $C_2$ are negative. 
In that case the local trajectory (\ref{sol:lf_local}) is in the concave down region for $a\leq a_m$ $(\eta \leq 0)$, resulting in tipping near $a=a_m$.


\section{The local approximation for $w$ in the ML model with $\Omega\gg 1$}

\label{append:E}

 The equations for the corrections $z_1$ and $W_1$ are
\begin{align}
&{z_1}_T \sim -\mu^{1/3} {z_0}_s +  \left[b + (h_0^0+h_1^0({\cal X}+z_0)+h_2^0 ({\cal X}+z_0)^2)-
g_K({\cal X}+z_0+{\cal D})W_0 \right], \label{exp:ML_inner_LO1} \\
 &{W_1}_T = -\mu^{1/3}{W_0}_s+ \nonumber\\
 &\left(\kappa(0) + \kappa'(0)({\cal X}+z_0) + \frac{\kappa''(0)}{2} ({\cal X}+z_0)^2\right)
 \left(w_{\infty}(0)+ w_{\infty}'(0)({\cal X}+z_0)+ \frac{w_{\infty}''(0)}{2}({\cal X}+z_0)^2
-W_0\right), \label{exp:ML_inner_LO2}
\end{align}
where we have kept linear and quadratic terms in ${\cal X}$ and $z_0$, analogous to (\ref{exp:hf_drifting_in_order1}).  We apply the solvability condition as in (\ref{ML_solv}), yielding for $W_0$
\begin{eqnarray}
 W_0 &=& \frac{\lim_{L \to \infty} L^{-1} \int_{0}^{L}
(\kappa(0) + \kappa'(0)({\cal X}+z_0) + \kappa''(0)({\cal X}+z_0)^2/2)
 (w_{\infty}(0)+ w_{\infty}'(0)({\cal X}+z_0)+ {w_{\infty}''(0)}({\cal X}+z_0)^2/2) \, dT}
 { \lim_{L\to \infty} L^{-1} \int_{0}^{L}
 (\kappa(0) + \kappa'(0)({\cal X}+z_0) + {\kappa''(0)}({\cal X}+z_0)^2/2) \, dT} \nonumber\\
  &\approx&  w_\infty(0) +  W_{00} + W_{01}z_0 + W_{02}z_0^2\label{w0approx}\\
W_{00} &=& \frac{w_\infty'(0)\kappa'(0) + w_\infty''(0)\kappa(0)/2}{W_{0d}}\frac{A^2}{2\Omega^2} + O(A^4/\Omega^4) \nonumber\\
 W_{01} &=& \frac{w_\infty'(0)\kappa(0)}{W_{0d}} - \frac{\kappa'(0)W_{00}}{W_{0d}} + O(A^2/\Omega^2)\nonumber\\
 W_{02} &=& W_{00}\left(\frac{2\Omega^2}{A^2} -\frac{\kappa''(0)}{2W_{0d}} +
\frac{2\kappa'(0)^2}{W_{0d}^2}\right) - \frac{w_\infty'(0)\kappa'(0)\kappa(0)}{W_{0d}^2} + O(A^2/\Omega^2)\nonumber\\
W_{0d} &=& \kappa(0) + \frac{A^2\kappa''(0)}{4\Omega^2} \nonumber
\end{eqnarray}
where we have used a Taylor series about $z_0$ to obtain a second order polynomial in
$z_0$ as an approximation for $W_0$.  Substituting (\ref{w0approx})
in the equation for $z_0$ written in terms of the original time scale yields
(\ref{exp:ML_inner_average}).  We note that the main contributions in the equation for $z_0$, used to determine the shift in the tipping point, are the terms
involving $h_j^0$ and $g_kw_\infty(0)$.
The additional terms with 
coefficients involving the derivatives of $w_\infty$ and $A^2/\Omega^2$
 give small corrections, since
these coefficients are small relative to $g_K$ and the terms with $h_j^0$.


\end{appendices}


\begin{thebibliography}{9}
\bibitem{Lenton2008}
        Lenton, T. M., H. Held, E. Kriegler, J. W. Hall, W. Lucht, S. Rahmstorf, and H. J. Schellnhuber (2008), {\em Tipping elements in the earth's climate system}, Proc. Nat. Acad. Sci. USA, 105(6) 1786-1793.
\bibitem{sutera1981}
        Sutera A. (1981), {\em On stochastic perturbation and long-term climate behaviour. Quarterly Journal of the Royal Meteorological Society}, 107(451):137-151.
\bibitem{guttal2008}
        Guttal V. and C. Jayaprakash (2008), {\em Changing skewness: an early warning signal of regime shifts in ecosystems}. Ecol Lett. 2008 May, 11(5):450-60. doi: 10.1111/j.1461-0248.2008.01160.x.Epub 2008 Feb 12.
\bibitem{Meisel2012}
        Meisel C, and C. Kuehn (2012), {\em  Scaling effects and spatio-temporal multilevel dynamics in epileptic seizures}. PLos ONE 7(2):e30371. doi:10.1371/journal.pone.0030371
 \bibitem{Dai2012}
Dai, L., D. Vorselen, K. S. Korolev and J. Gore (2012), {\em  Generic indicators for loss of resilience before a tipping point leading to population collapse}, Science 336, 1175, doi:10.1126/science.1219805.
\bibitem{ashwin2012}
        Ashwin P., S. Wieczorek, R. Vitolo and P. Cox (2012), {\em Tipping points in open systems: bifurcation, noise-induced and rate-dependent Examples in the Climate System}, Phil. Trans. R. Soc. A 370, 1166-1184, doi: 10.1098/rsta.2011.0306.
\bibitem{mandel1987}Mandel P and Erneux T (1987){\em  The slow passage through a steady bifurcatoin: delay and memory effects}, J. Stat. Phys., 48 1059-1070.
\bibitem{baer1989}Baer S.M., Erneux T. and Rinzel J. (1989) {\em The slow passage through a Hopf bifurcation: Delay, memory
effects, and resonance}, SIAM J. Applied Math., 49, 55-71.
\bibitem{haberman1979}
        Haberman R. (1979), {\em Slowly varying jump and transition phenomena associated with algebraic bifurcation}, SIAM J. Appl. Math., 37(1), 69-106.
\bibitem{fraedrich1978}
Fraedrich K. (1978), {\em Structural and stochastic analysis of a zero-dimensional climate system}, Quart. J. R. Met. Soc., 104, pp.461-474.
\bibitem{eisenman2009} Eisenman I., and J. S. Wettlaufer (2009), {\em Nonlinear threshold behavior during the loss of Arctic sea ice}, Proc. Natl. Acad. Sci. U. S. A., 106(1), 28-32, doi:10.1073/pnas. 0806887106.
\bibitem{abbot2011} Abbot, D.S., Silber M., and Pierrehumbert R.T. (2011), {\em Bifurcations leading to
summer Arctic sea ice loss}, J. Geophysical Res., Vol. 116, D19120, doi:10.1029/2011JD015653.
 \bibitem{merryfield2008}
         Merryfield, W. J., M. M. Holland, A. H. Monahan (2008), {\em Multiple Equilibria and Abrupt Transition in Arctic Summer Sea Ice extent,} Arctic Sea Ice Decline: Observations, projections, mechanisms, and implications. Bitz, C. M., Deweaver, E. T., Tremblay, L. B., Eds. American Geophysical Union, 151-174.
\bibitem{sieber2012}
        Sieber J. and J. M. Thompson (2012), {\em Nonlinear softening as a predictive precursor to climate tipping}. Philosophical Transactions of the Royal Society A: Mathematical, Physical and Engineering Sciences, vol. 370, issue 1962, pp. 1205-1227.
\bibitem{thompson2011}
        Thompson J. M. T. and J. Sieber (2011), {\em Climate tipping as a noisy bifurcation: a predictive technique}, IMA Journal of Applied Mathematics 76, 27-46
\bibitem{berglund2010}
        Berglund N. and Barbara Gentz, in C. Laing and G. J. Lord (Edts), (2010) {\em Stochastic Methods in Neuroscience}, Oxford University press. pp. 65-93.
\bibitem{Zhu2014} Zhu, J. and R. Kuske, {\em Multiple scale WKB-type approximations for probability densities in stochastic delayed saddle node bifurcations},  
preprint.
\bibitem{kuske1999}
        Kuske R. (1999), {\em Probability densities for noisy delay bifurcations}, Journal of Statistical Physics, Vol. 96, Issue 3-4, pp 797-816.
\bibitem{dakos2008}
        Dakos, V., M. Scheffer, E. H. van Nes, V. Brovkin, V. Petoukhov, and H. Held (2008), {\em Slowling down as an early warning signal for abrupt climate change}, Proc. Nat. Acad. Sci. USA, 105(38), 14308-14312.
\bibitem{Scheffer2012}Scheffer M, Carpenter S. R., Lenton T.M., Bascopmpte J., Brock W., Dakos V., van de Koppel J., 
van de Leemput I. A., Levin S. A., van Nes E. H., Pascual M., Vandermeer J. (2012) {\em Anticipating critical transitions}, Science,
338,  344-348, DOI: 10.1126/science.1225244.
\bibitem{tung2013} K.K. Tung  and J. Zhou  (2013), {\em Using data to attribute episodes of warming and cooling in instrumental records},  Proc. Nat. Acad. Sci. USA, 110(6) 2058-2063.
\bibitem{Bender1978}
 C.M. Bender and S.A. Orszag, (1978)  {\em  Advanced Mathematical Methods for Scientists
and Engineers}, section 11.3, McGraw-Hill, New York (1978).
\bibitem{Kevorkian1996}J. Kevorkian and J.D. Cole, (1981) Perturbation Methods in Applied Mathematics,
Appl. Math. Sciences 34, Springer, New York; (1996) Multiple Scale and
Singular Perturbation Methods, Appl. Math. Sciences 114, Springer, New
York.
\bibitem{strogatz2000}
        Strogatz S. H. (2000), Nonlinear Dynamics and Chaos, (1995), Addison-Wesley, Reading, MA. 
\bibitem{erneux1989} T. Erneux and J.-P. Laplante, (1989), {\em Jump transition due to a time-dependent
bifurcation parameter in the bistable ioadate-arsenous acid reaction}, J. Chem.
Phys. 90, 6129-6134.
\bibitem{Jung1990}
P. Jung, G. Gray, and R. Roy, (1990), {\em  Scaling law for dynamical hysteresis}, Phys.
Rev. Lett. 65, 1873-1876.
\bibitem{Hohl1995}A. Hohl, H.J.C. van der Linden, R. Roy, G. Goldstein, F. Broner, and S.H.
Strogatz, (1995), {\em Scaling laws for dynamical hysteresis in a multidimensional laser
system}, Phys. Rev. Lett. 74, 2220-2223.
        nd Chaos: With Applications to Physics, Biology, Chemistry, and Engineering. Cambridge, MA: Perseus.
\bibitem{maykut1971}
Maykut G.A. and N. Untersteiner (1971).{\em Some results from a time-dependent thermodynamic model of sea ice}, J. Geophys. Res. 76, 1550-1575. 
\bibitem{armour2011} Armour, K.C., I. Eisenman, E. Blanchard-Wrigglesworth, K.E. McCusker, and C.M. Bitz (2011), {\em The reversibility of sea ice loss in a state-of-the-art cliamte model}, Geophysical Research Letters 38, L16705.
\bibitem{eisenman2012}
        Eisenman I. (2012), {\em Factors controlling the bifurcation structure of sea ice retreat}, Journal of Geophysical Research, VOL. 117, D01111, doi:10.1029/2011JD016164.
\bibitem{Morris1981}
        Morris C. and H. Lecar (1981),{\em Voltage oscillations in the barnacle giant 
muscle.}  Biophys J. 35, 193-213.
\bibitem{MLweb} \mbox{http://www.scholarpedia.org/article/Morris-Lecar$\_$model.}
\bibitem{Yu2013}
        Yu N., Y. Li and R. Kuske (2013), A computational study of spike time reliability in two types of threshold dynamics, J. Math. Neuroscience, 3, doi:10.1186/2190-8567-3-11.
\bibitem{wiki_ML} \mbox{http://en.wikipedia.org/wiki/Morris-Lecar$\_$model.}
\bibitem{Lee2007} Lee S., {\em Study on the onset bifurcations of a Morris-Lecar neuron under a periodic current},  J. Korean Physical Society, 50:346-350, 2007.
\bibitem{Lee2008} Lee S., {\em  Bifurcation analysis of a mode-locking structure in a strongly forced Morris-Lecar neuron},  J. Korean Physical Society, 
52:11-16, 2008.
\bibitem{borowski2010}
P. Borowksi and R. Kuske (2010), {\em  Characterizing noisy mixed mode oscillations in neuronal models}, Chaos, 20,  043117.
\end{thebibliography}
\end{document}